\documentclass[11pt]{article}
\usepackage{amsmath,amssymb,amsthm,mathtools}
\usepackage[left=1in,right=1in,top=1in,bottom=1in]{geometry}
\usepackage{xcolor}
\usepackage{graphicx}
\usepackage{xurl}
\usepackage{microtype}
\usepackage{hyperref}
\newcommand{\originalrevisioncolor}{black}
\newcommand{\revisioncolor}{black}
\newcommand{\rd}[1]{#1}
\newcommand{\rev}[1]{{\color{\revisioncolor}#1}}
\hypersetup{hidelinks}
\DeclareRobustCommand{\changed}[1]{#1}
\DeclareRobustCommand{\mathred}[1]{#1}

\newcommand{\doi}[1]{\href{https://doi.org/#1}{\nolinkurl{#1}}}

\newtheorem{theorem}{Theorem}[section]
\newtheorem{lemma}[theorem]{Lemma}
\newtheorem{proposition}[theorem]{Proposition}
\newtheorem{corollary}[theorem]{Corollary}
\theoremstyle{definition}
\newtheorem{definition}[theorem]{Definition}
\newtheorem{remark}[theorem]{Remark}
\numberwithin{equation}{section}

\newcommand{\pl}{(-\Delta_p)^s}
\newcommand{\Cw}{\mathcal C_\omega}
\newcommand{\R}{\mathbb R}
\newcommand{\Sp}{\mathbb S^{N-1}}
\newcommand{\Ls}{L^{p-1}_{sp}(\R^N)}
\newcommand{\dist}{\operatorname{dist}}

\title{\large\bf Fundamental solutions of the fractional $p$-Laplacian in cone-like domains, with applications to Liouville-type theorems and global H\"older regularity}
\author{Reshmi Biswas\thanks{Universidad del B\'io-B\'io, Departamento de Matem\'atica, Concepci\'on, Chile. \texttt{rbiswas@ubiobio.cl}}
\and Alexander Quaas\thanks{Universidad T\'ecnica Federico Santa Mar\'ia, Departamento de Matem\'atica, Avda. Espa\~na 1680, Valpara\'iso, Chile. \texttt{alexander.quaas@usm.cl}}}
\date{}

\begin{document}
\maketitle

\begin{abstract}
Let $\mathcal C_\omega$ be the cone generated by a proper $C^2$ domain $\omega\subset\mathbb S^{N-1}$. We study positive homogeneous solutions of the fractional $p$-Laplace equation in $\mathcal C_\omega$ with zero exterior values. For $p\ne2$, the operator is nonlinear and not uniformly elliptic. For every $p>1$, we identify two extremal homogeneities $\beta^+>0>\beta^-$, prove that every fundamental solution, if one exists, has one of these exponents, and construct homogeneous eigenfunctions whose eigenvalues have signs determined by the position of the homogeneity relative to $\beta^\pm$.

As an application, we obtain the Lane-Emden threshold $q^*=p-1+sp/\beta^+$ for every $p>1$: nontrivial nonnegative supersolutions of $(-\Delta_p)^s u\ge u^q$ in $\mathcal C_\omega$ do not exist for $0<q<q^*$, whereas homogeneous supersolutions exist for $q>q^*$. We also establish Hopf-type lower bounds at conical boundary points with the optimal exponent. When $p-1>sp,$ we construct both fundamental solutions, one singular at the vertex and the other vanishing there. For $p=2$, these solutions are known; our method recovers them for every $s\in(0,1)$.

The homogeneous eigenfunctions also serve as barriers, yielding boundary decay and global H\"older regularity for $(-\Delta_p)^s u=f$ in bounded Lipschitz domains, with zero exterior values and $f$ continuous up to the boundary. The admissible H\"older exponents are determined by the exterior cones. All results hold for every $N\ge2$; in particular, $N>sp$ is not needed.


\end{abstract}

\section{Introduction}

Positive homogeneous solutions in cones describe \changed{the behaviour of solutions near} conical boundary \changed{points and dilation-invariant singularities}. Their degree of homogeneity \changed{reflects} the interaction between the operator and the opening of the cone and governs Liouville theorems, boundary decay and Hopf\changed{-type} estimates. \changed{This theory is well developed for} local operators and uniformly elliptic nonlocal operators\changed{; see} \cite{KM,Tolksdorf,PV,ASS,AS,BB,Michalik,gab}. \changed{Here we} develop it for the fractional $p$-Laplacian, which is \changed{nonlocal and, for $p\ne2$,} nonlinear and not uniformly elliptic. \changed{Our results identify} the extremal \changed{homogeneities and} {the signs of the corresponding homogeneous eigenvalues}, \changed{determine} the Liouville threshold, and \changed{give} optimal Hopf \changed{bounds} and decay rates \changed{at infinity}. \changed{Our} method \changed{combines} comparison in bounded domains, interior and boundary regularity, and the homogeneity of the operator.

Let $N\ge 2$, $s\in(0,1)$ \changed{and} $1<p<\infty$. For \changed{a} measurable \changed{function} $u:\R^N\to\R$, the fractional $p$-Laplacian is \changed{given, whenever the principal value exists, by}
\begin{equation}\label{eq:pl}
\pl u(x)=\mathrm{P.V.}\int_{\R^N}\frac{\Psi_p(u(x)-u(y))}{|x-y|^{N+sp}}\,dy,\qquad \Psi_p(t):=|t|^{p-2}t,
\end{equation}
{\mathred{Here $\Psi_p(0):=0$,} the \changed{{normalization}} constant is \changed{set} equal to $1$\changed{, and P}.\changed{V.\ denotes the Cauchy principal value.}}

Throughout this article, \mathred{$N\ge2$, $0<s<1$, $p>1$,} and $\omega\subset\Sp$ is a domain of class $C^2$ with $\overline\omega\neq\Sp$\changed{. We} set $\Cw:=\{x\in\R^N\setminus\{0\}: x/|x|\in\omega\}$ and study the \changed{Dirichlet} problem
\begin{equation}\label{mp}
\pl u=0 \ \text{ in } \Cw,\qquad u=0\ \text{ in }\R^N\setminus\Cw .
\end{equation}
A {fundamental solution} of \eqref{mp} is a viscosity solution of \eqref{mp} (in particular\changed{,} $u\in\Ls$\changed{;} see Section~\ref{sec:prel}) \changed{that is} positive in $\Cw$, continuous in $\R^N\setminus\{0\}$ \changed{and} homogeneous: $u(tx)=t^{-\beta}u(x)$ for $t>0$ \changed{and} some $\beta\ne0$.

Following \cite{ASS,gab}\changed{,} we define the extremal exponents
\begin{equation}\label{def:beta}
\beta^+:=\sup\{\beta>0:\ \mathcal S_\beta(\omega)\neq\emptyset\},\qquad
\beta^-:=\inf\Big\{\beta\in\big(-\tfrac{sp}{p-1},0\big]:\ \mathcal S_\beta(\omega)\neq\emptyset\Big\},\qquad \alpha^-:=-\beta^-,
\end{equation}
{where $\mathcal S_\beta(\omega)$ is \changed{defined} in \changed{Section~\ref{sec:prel}}.}

A {homogeneous eigenpair} with exponent $\beta$ is a pair $(\mu,U)$, where $\mu\in\R$ and \mathred{$U\in C(\R^N\setminus\{0\})\cap\Ls$} is $(-\beta)$-homogeneous \changed{and} nonnegative, \changed{is} positive in $\Cw$ \changed{and} zero outside $\Cw$, and \changed{solves} $\pl U=\mu|x|^{-sp}U^{p-1}$ in $\Cw$.

Our first \changed{main} result \changed{concerns positive exponents.}
\begin{theorem}\label{mt}
\begin{itemize}
\item[(a)] $\beta^+$ is well defined and $0<\beta^+<\frac N{p-1}$; if $N>sp$, then $\beta^+\ge\frac{N-sp}{p-1}$.
\item[(b)] For every \mathred{$p>1$ and} $\beta\in(0,\frac N{p-1})$\changed{,} there is a homogeneous eigenpair with exponent $\beta$, and every such eigenpair has $\mu>0$ if $\beta<\beta^+$ and $\mu<0$ if $\beta>\beta^+$.
\item[(c)] If \mathred{$p>1$ and $\Phi^+$} is a fundamental solution of \eqref{mp} with exponent $\beta>0$, then $\beta=\beta^+$.
\item[(d)] If $p>\frac1{1-s}$ \textup{(}\changed{equivalently}, $p-1>sp$\textup{)} or $p=2$, then \eqref{mp} has a fundamental solution $\Phi^+$ with exponent $\beta^+$\changed{.}
\end{itemize}
\end{theorem}

\mathred{The case $p=2$ in Theorem~\ref{mt}}(\changed{d)} is known\changed{;} see \cite[Theorem~1.1]{gab} and \cite{BB}. For general $p$, \changed{the existence of} $\Phi^+$ \changed{follows from} a separation property for eigenpairs with exponent $\beta^+$ (Section~\ref{sec:sep})\changed{. The} restriction $p-1>sp$ enters only in the error estimate of Lemma~\ref{lem:error}\changed{,} through the choice of $\alpha$ \changed{in \eqref{eq:alpha}} when $p\ge2$ and a dyadic sum when $p<2$\changed{,} and in \changed{the analogous estimate} in Theorem~\ref{thm:phiminus}. {We do not establish uniqueness of the profile up to a multiplicative constant.} No restriction of the form $N>sp$ is needed. When $N>sp$, the radial function $|x|^{-\frac{N-sp}{p-1}}$ is a homogeneous supersolution, which gives the lower bound in Theorem~\ref{mt}(a). When $N<sp$, every $|x|^{-\beta}$ with {$0<\beta<\frac N{p-1}$} is a subsolution by Theorem~\ref{thm:qp} (and for $N=sp$ the radial fundamental solution is logarithmic \cite{qp}), {so a radial model is unavailable; positive supersolutions} with $\beta>0$ are then produced by the exterior condition, through the positivity of the eigenvalue at $\beta=0$ (Lemma~\ref{lem:zero} and Corollary~\ref{cor:betaplus}). {This is consistent with the local case, in which separable solutions also exist in every smooth cone when} $N\le p$ \cite{Tolksdorf,PV}.

{The \changed{principal} difficulty in \changed{studying} the fractional $p$-Laplacian \changed{for $p\ne2$ stems from} its strong \changed{nonlinearity} and \changed{lack of} uniform \changed{ellipticity}. \changed{For} the Pucci and Isaacs operators studied in \cite{gab}, differences of solutions are controlled by extremal operators with fixed ellipticity constants, and this structure drives both the construction and the comparison of homogeneous profiles. \changed{For} the fractional $p$-Laplacian \changed{with $p\ne2$}, \changed{this} property \changed{is lost. Indeed,} writing $\Psi_p(a)-\Psi_p(b)=A_p(a,b)(a-b)$ with $A_p(a,b)=(p-1)\int_0^1|b+t(a-b)|^{p-2}\,dt$ \mathred{for $(a,b)\ne(0,0)$}, \mathred{with the convention at $(0,0)$ specified in Section~\ref{sec:prel}, we find that the} difference of two solutions \changed{satisfies}, at best, a linear equation whose kernel $A_p(\delta u,\delta v)|x-y|^{-N-sp}$, \changed{where} $\delta u=u(x)-u(y)$ \changed{and $\delta v=v(x)-v(y)$}, depends on the solutions and, for $p\ne2$, \changed{admits} no two-sided bounds \changed{independent of them}.

\changed{The geometry introduces {an additional difficulty}. Unlike} the radial \changed{problem} studied in \cite{qp} (see also \cite{BMS} \changed{and} Theorem~\ref{thm:qp}), \mathred{and unlike the local case, problem~\eqref{mp} does not reduce to an eigenvalue problem on $\omega$: the exterior condition couples the angular profile to every radial scale, so we study a weighted homogeneous eigenvalue problem in the whole cone.}\\

\changed{As a first application} of Theorem~\ref{mt}\changed{, we obtain a Liouville-type theorem for the Lane--Emden inequality}. \changed{Its threshold} is \mathred{determined by $\beta^+$, and the argument applies for every $p>1$, independently of the restriction used to construct fundamental solutions.}}

\begin{theorem}\label{thm:liou}
Let $q>0$\changed{,} and \changed{set} $q^*:=p-1+\frac{sp}{\beta^+}$.
\begin{itemize}
\item[(a)] If $q<q^*$ and $u\in C(\R^N\setminus\{0\})\cap\Ls$, $u\ge0$ in $\R^N$, is a viscosity supersolution of $\pl u\ge u^q$ in $\Cw$ (with $u^q:=(u^+)^q$), then $u\equiv0$ in $\R^N\setminus\{0\}$.
\item[(b)] If $q>q^*$, there exists a nonnegative $u\in C(\R^N\setminus\{0\})\cap\Ls$, positive in $\Cw$ and homogeneous, which is a viscosity supersolution of $\pl u\ge u^q$ in $\Cw$.
\end{itemize}
\end{theorem}
\changed{Together}, \changed{parts}~(a) and~(b) of \changed{Theorem~\ref{thm:liou} show that} $q^*$ is sharp\changed{,} except possibly \changed{at} $q=q^*$\changed{. The critical} case \changed{remains} open for $p\neq2$\changed{;} for $p=2$\changed{,} it belongs to the nonexistence range \cite[Theorem~1.2]{gab} (\changed{see} Remark~\ref{rem:crit}). \changed{Formally, when} $\omega=\Sp$\changed{, we recover} $q^*=\frac{N(p-1)}{N-sp}$, the Serrin-type exponent { (for $N>sp$)}. \changed{In} Theorem~\ref{thm:liou}(a)\changed{,} $u$ \changed{need} not vanish outside $\Cw$.

\changed{We next} state the \changed{counterpart} of Theorem~\ref{mt} for the \changed{negative} exponent $\beta^-$ and its application to a \changed{Hopf-}type \changed{estimate}. The profile $\Phi^-$ (see Section~\ref{sec:neg}\changed{, in particular Subsection}~\ref{7.4}) vanishes at the vertex and grows at infinity\changed{; in a} half-space, \changed{it is} $(x\cdot e)_+^s$\changed{. Following} \cite{ASS,gab}\changed{,} we \changed{also} call \changed{this profile} a fundamental solution\changed{. It} plays the role of a Martin-type kernel with pole at infinity.

\begin{theorem}\label{thm:neg}
\begin{itemize}
\item[(a)] \mathred{Let $p>1$. We have} $-\frac{sp}{p-1}<\beta^-<0$; $\mathcal S_{\beta^-}(\omega)\ne\emptyset$ and $\mathcal S_\beta(\omega)=\emptyset$ for $-\frac{sp}{p-1}<\beta<\beta^-$.
\item[(b)] For \mathred{$p>1$ and} every $\beta\in(-\frac{sp}{p-1},0]$\changed{,} there is a homogeneous eigenpair with exponent $\beta$, and every such eigenpair has $\mu>0$ if $\beta>\beta^-$ and $\mu<0$ if $\beta<\beta^-$.
\item[(c)] If \mathred{$p>1$ and} $\Phi$ is a fundamental solution of \eqref{mp} with exponent $\beta<0$, then $\beta=\beta^-$; moreover\changed{,} \eqref{mp} has no positive $0$-homogeneous solution.
\item[(d)] If $p>\frac1{1-s}$ or $p=2$, \changed{then} \eqref{mp} has a fundamental solution $\Phi^-$ with exponent $\beta^-$, \changed{that is,} $\Phi^-(tx)=t^{\alpha^-}\Phi^-(x)$.
\item[(e)] \textup{(\changed{Hopf-}type result at conical points)} Let $\Omega\subset\R^N$ be open, $x_0\in\partial\Omega$, and $x_0+\mathcal R(\Cw\cap B_\rho)\subset\Omega$ for a rotation $\mathcal R$ and some $\rho>0$. If $v\in C(\R^N)\cap\Ls$, $v\ge0$ in $\R^N$, $v>0$ in $\Omega$, is a viscosity supersolution of $\pl v\ge0$ in $\Omega\cap B_\rho(x_0)$, then for every compact $K\subset\omega$ there is $c>0$ with $v(x_0+\mathcal Rx)\ge c|x|^{\alpha^-}$ for $x\in\Cw\cap B_{\rho/12}$, $x/|x|\in K$. This holds for every $p>1$, and the exponent $\alpha^-$ cannot be replaced by any $\alpha<\alpha^-$.
\end{itemize}
\end{theorem}
For a hemisphere\changed{,} $\alpha^-=s$ for every $p>1$, \changed{which gives the power in} the classical \changed{Hopf lemma} (Remark~\ref{rem:halfspace}). \mathred{For a general cone, Corollary~\ref{cor:hopfcone} identifies the optimal lower-bound exponent at the vertex. Remark~\ref{rem:decay} gives the corresponding} optimal decay $|x|^{-\beta^+}$ at infinity\mathred{. When $p=2$, the fundamental solutions are known; see \cite[Theorem~1.1]{gab}} for \changed{cones contained in a hemisphere and \cite{BB,Michalik}}.

Using the eigenfunctions as barriers\changed{,} {we also establish boundary decay and global} H\"older regularity for $\pl u=f$ in \changed{bounded} Lipschitz domains, with \mathred{zero exterior values and $f$ continuous up to} the \changed{boundary} (Theorem~\ref{cor:lipreg}). \mathred{For every admissible narrower exterior aperture $\theta$, the same exponent $0<\alpha<\alpha^-(\omega^c_\theta)$ is available both in the boundary estimate $|u|\le C\dist(x,\partial\Omega)^\alpha$ and in the global conclusion $u\in C^\alpha(\R^N)$. Here $\omega^c_\theta$ is the spherical complement of the narrower exterior cap, as defined in that theorem. These conclusions hold for every $p>1$.}\\

\emph{\bf Novelty and related work.} To the best of our knowledge, {for $p\ne2$,} Theorems~\ref{mt}(d) and \ref{thm:neg}(d) give the first positive homogeneous solutions vanishing outside a general cone for a nonlocal operator that is not uniformly elliptic, and Theorem~\ref{thm:neg}(e) is the first Hopf\changed{-}type bound at conical boundary points, in particular in Lipschitz domains, for the fractional $p$-Laplacian. \changed{{Moreover, a further important application is} the global H\"older regularity in bounded Lipschitz domains (Theorem~\ref{cor:lipreg}).} {A genuinely nonlocal feature is that $\beta^+$ exists for every $N\ge2$: when $N<sp$ every decaying radial power {$|x|^{-\beta}$, $0<\beta<\frac N{p-1}$,} is a subsolution, and when $N\le sp$ the decaying supersolutions {arise from the exterior condition} (Corollary~\ref{cor:betaplus}); in \cite{gab} the positive exponent is obtained under $\tilde N^+>2\alpha$.} The known Hopf\changed{-}type lemmas for this operator assume an interior ball condition \cite{hopf,OS,Ataei}; the rate $\dist(x,\partial\Omega)^s$ fails at inward corners already for linear stable operators \cite{Grube}; and for uniformly elliptic Isaacs operators \cite[Theorem~1.4]{gab} gives every exponent larger than $\alpha^-$ when $\omega$ is contained in a hemisphere. H\"older continuity up to the boundary under an exterior measure density condition is due to \cite{KKP}, with an unspecified exponent, and global H\"older regularity with the boundary bound $|u|\le C\dist(x,\partial\Omega)^s$ in $C^{1,1}$ domains \changed{is established in} \cite{reg}. For $p=2$ the profiles and the relation $\beta^+=N-2s-\beta^-$ are known \cite{BB,Michalik,gab}; our separation argument recovers them for every $s\in(0,1)$. Homogeneous $p$-harmonic functions in cones were studied in \cite{KM,Tolksdorf,PV,ALV}, fundamental solutions of fully nonlinear local equations in \cite{ASS,ASS2,AS}, Liouville theorems for Pucci operators in \cite{CL}, and fundamental solutions and Liouville theorems for nonlinear integral operators in $\R^N$ \changed{in} \cite{FQ}. For the fractional Laplacian, cone profiles, Green functions, Martin kernels and heat kernels appear in \cite{BB,Michalik,BSS,Terracini,BGcone,BKLP}.\\

Regarding the tools and methodology \changed{adopted} in this article, {we highlight} \changed{a few facts}.  {In addition to standard comparison in bounded domains, our proofs use} the following results from the regularity theory of fractional $p$-Laplace equations, \changed{stated} in Section~\ref{sec:prel}. The equivalence of weak and viscosity solutions \cite{lindgren,bm} allows us to use pointwise test functions, also in the singular range $1<p\le\frac2{2-s}$, together with weak comparison. Local boundedness, interior H\"older estimates, the weak Harnack inequality and Caccioppoli estimates with tail \cite{dkp,dkp2,Cozzi} give compactness of eigenfunctions and stability. The Dirichlet problem in $C^{1,1}$ domains, with global H\"older regularity and the boundary bound $|u|\le C\dist(x,\partial\Omega)^s$ \cite{reg} (see also \cite{IMS20,IM}), and the Hopf lemma \cite{hopf}, used in a quantitative form, give the two-sided boundary behaviour $c\,d^s\le U\le C\,d^s$ of homogeneous eigenfunctions on the unit sphere. For $p\ge2$, the higher H\"older regularity of \cite{BLS} provides the interior estimate behind the error bound of Section~\ref{sec:exist}. \mathred{Together with its counterpart for $1<p<2$ in \cite{GL}, it also allows the boundary decay exponent to be retained in the global H\"older estimate.} {For further regularity results, see} \cite{BP,LindgrenHolder,BL,GL,Diening,KMSdata,KMSself,BT2025}; {and for gradient estimates, see also}  \cite{GJS,BRS,GJ2026}.
% and their bearing on the remaining range $p-1\le sp$ are discussed in Remark~\ref{rem:sepopen}.

The eigenpairs come from a continuation theorem applied to a compact map built from minimal cone solutions (Section~\ref{sec:eig}). Proposition~\ref{prop:nonex} compares a modified homogeneous subsolution with a positive supersolution on bounded annuli; the different radial rates determine the eigenvalue signs, without strong comparison for functions of the same homogeneity. At $\beta^+$, if a limiting eigenvalue were positive, the perturbation $|x|^{-\varepsilon}\underline U^{\,a}$, with $a<1$ close to $1$ and $\varepsilon=2\beta^+(1-a)$, would be a supersolution with homogeneity larger than $\beta^+$: the radial error is of order $\varepsilon$ uniformly up to the lateral boundary when $p-1>sp$ (Lemma~\ref{lem:error}), and the concavity of $t\mapsto t^a$ gives a gain of order $(1-a)\log(1/\underline U)$ there (Lemma~\ref{lem:gain}). The negative exponent is treated in the same way.\\

{The article is organised as follows.}  In Section~\ref{sec:prel} we gather preliminaries, in  Section~\ref{sec:fund} we study $\beta^+$ and the nonexistence principle, in  Section~\ref{sec:eig} we construct the eigenfunctions, and finally in  Sections~\ref{sec:exist}, \ref{sec:liou} and \ref{sec:neg} we prove Theorems~\ref{mt}, \ref{thm:liou} and \ref{thm:neg}. {Section~}\ref{sec:sep} contains the separation principle and the case $p=2$.

\section{Preliminaries}\label{sec:prel}

\subsection{Notation and notions of solution}
We write $B_r(x)$ for balls, $B_r=B_r(0)$, $d(x)=\dist(x,\R^N\setminus\Cw)$ for $x\in\changed{\R^N}$, $\mathcal H^{N-1}$ for the surface measure on $\Sp$, and
\[
\Ls:=\Big\{u\in L^{p-1}_{\rm loc}(\R^N):\ \int_{\R^N}\frac{|u(x)|^{p-1}}{(1+|x|)^{N+sp}}\,dx<\infty\Big\}.
\]
\mathred{For $A\subset\Sp$, we use $\mathcal C_A:=\{r\theta:r>0,\ \theta\in A\}$, also when $A$ is compact. Values at the origin do not affect any integral; a homogeneous profile is extended continuously there when possible, and otherwise assigned the value $0$. Convergence $u_j\to u$ in $\Ls$ means
\[
\int_{\R^N}\frac{|u_j(y)-u(y)|^{p-1}}{(1+|y|)^{N+sp}}\,dy\longrightarrow0.
\]}
Since $0\notin\Cw$, we have $d(x)\le|x|$ for $x\in\Cw$; in particular $d\le1$ on $\Cw\cap\Sp$. For $u\in\Ls$ and $x\in\R^N$ \mathred{with $u(x)$ finite,} the integral in \eqref{eq:pl} converges away from $x$; near $x$ it is understood {in the principal-value sense.} In the {singular range} $1<p\le\frac{2}{2-s}$ the operator is not defined at critical points of general $C^2$ functions \cite{lindgren}, which motivates the following class: for an open set $D$ and $\alpha>0$, $C^2_\alpha(D)$ is the set of $\varphi\in C^2(D)$ such that
\[
\sup_{x\in \rev{D\setminus\mathcal N_\varphi}}\left[\frac{\min\{\dist(x,\mathcal N_\varphi),1\}^{\alpha-1}}{|\nabla \varphi(x)|}+\frac{|D^2\varphi(x)|}{\dist(x,\mathcal N_\varphi)^{\alpha-2}}\right]<\infty,\qquad \mathcal N_\varphi:=\{\nabla\varphi=0\}.
\]
\mathred{In the coefficient identity from the Introduction, set $A_p(0,0):=0$ for $p\ne2$ and $A_2(0,0):=1$; the integral formula applies when $(a,b)\ne(0,0)$.} Throughout, $q_0:=\max\{2,\frac{sp}{p-1}+1\}$; note that $q_0\ge2$ and $q_0(p-1)>sp$.

\begin{definition}[Viscosity solutions, \cite{lindgren}]\label{def:visc}
Let $D\subset\R^N$ be open and $F:D\times\R\to\R$ continuous. A function $u:\R^N\to[-\infty,\infty]$ is a {viscosity supersolution} of $\pl u=F(x,u)$ in $D$ if $u$ is lower semicontinuous in $D$, $u>-\infty$ in $D$, $u<\infty$ {almost everywhere in} $\R^N$, $u^-\in\Ls$, and the following holds: whenever $x_0\in D$, $B_r(x_0)\subset D$ and $\varphi\in C^2(B_r(x_0))$ satisfy $\varphi(x_0)=u(x_0)$, $\varphi\le u$ in $B_r(x_0)$, and either $\nabla\varphi(x_0)\neq0$, or $p>\frac{2}{2-s}$, or $x_0$ is an isolated critical point of $\varphi$ and $\varphi\in C^2_\alpha(B_r(x_0))$ for some $\alpha>\frac{sp}{p-1}$, then
\[
\pl\varphi_r(x_0)\ge F(x_0,u(x_0)),\qquad \varphi_r:=\varphi\ \text{in }B_r(x_0),\quad \varphi_r:=u\ \text{in }\R^N\setminus B_r(x_0).
\]
Subsolutions are defined symmetrically (upper semicontinuity, $u^+\in\Ls$, test functions touching from above, reversed inequality), and a viscosity solution is both.
\end{definition}

{We say that} a test function $\varphi$ as in Definition~\ref{def:visc} is \emph{admissible} at $x_0$. {We use the notation ``}$\pl u\ge F$ in $D$'' for ``$u$ is a viscosity supersolution of $\pl u=F$ in $D$'', and similarly for subsolutions.

The Sobolev space $W^{s,p}(D)$ consists of functions in $L^p(D)$ with finite Gagliardo seminorm
	\[
	[u]_{W^{s,p}(D)}^p:=\iint_{D\times D}\frac{|u(x)-u(y)|^p}{|x-y|^{N+sp}}\,dx\,dy,
	\]
	and $W^{s,p}_{\mathrm{loc}}(D)$ means that this property holds on every open set compactly contained in $D$. We use $W^{s,p}_0(D)$ for the closure of $C_c^\infty(D)$ in $W^{s,p}(\R^N)$, with functions extended by zero outside $D$.

\begin{definition}[Weak solutions, \cite{lindgren,bm}]
	
	For $u\in W^{s,p}_{\mathrm{loc}}(D)\cap\Ls$, the weak supersolution inequality corresponding to the {{normalization}} in \eqref{eq:pl} is
	\[
	\frac12\iint_{\R^N\times\R^N}
	\frac{\Psi_p(u(x)-u(y))(\phi(x)-\phi(y))}{|x-y|^{N+sp}}\,dx\,dy
	\ge\int_D F(x,u(x))\phi(x)\,dx
	\]
	for every nonnegative $\phi\in C_c^\infty(D)$, whenever the right-hand side is locally integrable. Weak subsolutions satisfy the reversed inequality, and weak solutions satisfy the {corresponding }equality. 
\end{definition}

\par

\medskip
{We also recall the following two elementary inequalities.} For $a,b\in\R$,
\begin{equation}\label{eq:elem}
|\Psi_p(a)-\Psi_p(b)|\le
\begin{cases}
2^{2-p}|a-b|^{p-1}, & 1<p\le2,\\[2pt]
(p-1)|a-b|\,(|a|+|b|)^{p-2}, & p\ge2,
\end{cases}
\end{equation}
and $|a+b|^{p-1}\le 2^{p-1}\big(|a|^{p-1}+|b|^{p-1}\big)$.
In \eqref{eq:elem}, {the second inequality follows from the mean value theorem}, since $|\Psi_p'(\tau)|=(p-1)|\tau|^{p-2}$ and $|\tau|\le|a|+|b|$ on the segment $[a,b]$; {the first is the standard} H\"older continuity of $\Psi_p$ for $p\le2$.

\subsection{\rd{The toolbox: radial powers, regularity and comparison}}
\rd{Here we \changed{collect} the results from the literature on which the proofs rest. The starting point is the action of $\pl$ on radial powers.}

\begin{theorem}[\cite{qp}]\label{thm:qp}
Let $\beta\in\big(-\frac{sp}{p-1},\frac{N}{p-1}\big)$. Then $v_\beta(x)=|x|^{-\beta}$ satisfies, pointwise and in the viscosity sense,
\begin{equation}\label{f.}
\pl v_\beta=C(\beta)\,|x|^{-\beta(p-1)-sp}\qquad\text{in }\R^N\setminus\{0\},
\end{equation}
where $C(\beta)$ is a continuous function \mathred{with $C(0)=0$. If $N\ne sp$, it} vanishes exactly at $\beta=0$ and $\beta=\frac{N-sp}{p-1}$, is positive strictly between these two values and negative otherwise.
\end{theorem}

\mathred{The sign formula for $N\ne sp$ is proved in \cite{qp}. The identity and continuity also hold when $N=sp$. Fix $e\in\Sp$ and set $C(\beta):=\pl v_\beta(e)$. For $\beta\ne0$, the local principal value exists by Lemma~\ref{lem:small}, since $\nabla v_\beta(e)\ne0$, and the remaining integrals converge because $-sp<\beta(p-1)<N$. Rotation and scaling give \eqref{f.}; ordering the patched tests gives the viscosity inequalities. The case $\beta=0$ is constant. On compact exponent intervals, $|\nabla v_\beta(e)|=|\beta|$ and $|D^2v_\beta|\le C|\beta|$ near $e$, so Lemma~\ref{lem:small} bounds the local contribution by $C|\beta|^{p-1}\rho^{p-sp}$, uniformly for sufficiently small $\rho$. Dominated convergence on its complement proves continuity, including at $\beta=0$.}

\rd{In the remaining results} $D$ is a bounded open set, $c\in C(\overline D)$ is nonnegative and $f\in C(\overline D)$, and we consider
\begin{equation}\label{eq:gen}
\pl u+c(x)\Psi_p(u)=f(x)\qquad\text{in }D .
\end{equation}

\begin{itemize}
\item[\textbf{(K1)}] \emph{\bf Equivalence \changed{of weak and viscosity solutions (\cite{lindgren,bm}):}} Viscosity supersolutions (subsolutions) of \eqref{eq:gen} which are bounded and lower (upper) semicontinuous in $\R^N$, and continuous in a neighbourhood of $\overline D$, are weak supersolutions (subsolutions), and continuous weak supersolutions (subsolutions) are viscosity supersolutions (subsolutions). For \eqref{eq:gen}\mathred{, when $N>sp$,} this is \cite[Theorems~1.1 and~1.3]{bm}, which cover every $p>1$, including the singular range with the test functions of Definition~\ref{def:visc}: the right-hand side $f(x)-c(x)\Psi_p(t)$ is nonincreasing in $t$ because $c\ge0$, and, since $u$ is bounded by some $M$, we may replace $\Psi_p(t)$ by $\Psi_p(\max\{-M,\min\{t,M\}\})$, which makes it uniformly continuous without changing the equation for $u$, as required in \cite[Theorem~1.1]{bm}; and the comparison property required in \cite[Theorem~1.3]{bm} holds for \eqref{eq:gen} by the weak comparison principle, again because $c\ge0$. For $f\equiv0$ and $c\equiv0$ see also \cite{lindgren}. The global test functions of \cite{bm} give the same notion as Definition~\ref{def:visc}, since $\varphi_r$ is such a test function and a global test function lies below the corresponding $\varphi_r$. When the right-hand side depends on $u$, as in $\pl u=F(x,u)$, we apply (K1) with the frozen datum $f(x):=F(x,u(x))$, which is continuous and bounded near $\overline D$; this is legitimate because Definition~\ref{def:visc} only evaluates $F$ at $(x_0,u(x_0))$, and it avoids any monotonicity assumption on $F$ in the variable $u$. Boundedness can always be achieved by modifying the functions far from $D$, see Remark~\ref{rem:K1} below; so we use (K1) for functions which are only locally bounded.
\mathred{For $N\le sp$, we have $p>2$ because $N\ge2$ and $s<1$. Freeze $g(x):=f(x)-c(x)\Psi_p(u(x))$. The bounded infimal-convolution approximants in the proof of \cite[Theorem~1.1]{bm} converge almost everywhere, and the local energy estimate \cite[(3.40)]{bm}, which does not require $N>sp$, gives
\[
\sup_\varepsilon\int_K\int_{\R^N}
\frac{|u_\varepsilon(x)-u_\varepsilon(y)|^p}{|x-y|^{N+sp}}\,dy\,dx<\infty
\qquad(K\Subset D).
\]
Bounded convergence replaces the critical Sobolev embedding used there. The {normalised} increments are bounded in $L^p$, and their nonlinear fluxes converge weakly in $L^{p/(p-1)}$ to the flux identified by almost-everywhere convergence. Pairing with smooth test increments gives the weak inequality; the frozen datum is independent of the fractional gradient. Upper convolutions give the subsolution case. The converse uses the dimension-independent weak comparison argument of \cite[Theorem~1.3]{bm}.}
\item[\textbf{(K2)}] \emph{\bf Comparison in bounded domains\changed{:}} Let $f_1\le f_2$ in $D$. If $u$ is a viscosity subsolution of \eqref{eq:gen} with $f=f_1$ and $w$ a viscosity supersolution of \eqref{eq:gen} with $f=f_2$, both continuous in a neighbourhood of $\overline D$ and in $\R^N\setminus Z$ for a closed null set $Z$ disjoint from $\overline D$, $u^+,w^-\in\Ls$, $u$ bounded above and $w$ bounded below in $\R^N$, and $u\le w$ in $\R^N\setminus D$, then $u\le w$ in $D$. \emph{Proof.} Let $M\ge\max\{\sup u,-\inf w\}$ be larger than $|u|,|w|$ in a neighbourhood $D_\delta$ of $\overline D$, and $\tilde u:=\max\{u,-M\}\ge u$, $\tilde w:=\min\{w,M\}\le w$, redefined as $M$ and $-M$ on $Z$ (with $\delta<\dist(Z,\overline D)$), so that they are semicontinuous in $\R^N$; since $Z$ is a null set outside $D$, this affects neither the weak formulations nor the a.e.\ ordering below. They coincide with $u,w$ in $D_\delta$, so by Lemma~\ref{lem:loc}(i) $\tilde u$ is a subsolution and $\tilde w$ a supersolution in $D$; they are bounded, continuous near $\overline D$, and $\tilde u\le\tilde w$ \changed{a.e.\ }in $\R^N\setminus D$ (since $u\le w$, $u\le M$ and $w\ge-M$ there). By (K1) they are weak sub- and supersolutions. For $\tau>0$, $\phi:=(\tilde u-\tilde w-\tau)^+$\mathred{, represented by $0$ on $Z$,} has compact support in $D$ (by continuity near $\partial D$), so it is an admissible test function (by density, since $\tilde u,\tilde w\in W^{s,p}_{\rm loc}(D)$ by (K1) and $\phi\in W^{s,p}$ has compact support in $D$); since $\tilde w+\tau$ is also a weak supersolution (because $c\ge0$), subtracting the two weak inequalities and using the monotonicity of $\Psi_p$ gives $\phi\equiv0$. Letting $\tau\to0$, $u=\tilde u\le\tilde w=w$ in $D$. Note that this truncation, unlike the cut-offs of Remark~\ref{rem:K1}, preserves the order outside $D$.
\item[\textbf{(K3)}] \emph{\bf Dirichlet problem \changed{(\cite{reg}):}} If $D$ is of class $C^{1,1}$ and $c,f\in L^\infty(D)$, there is a unique weak solution \mathred{$u\in W^{s,p}(\R^N)$} of \eqref{eq:gen} with $u=0$ in $\R^N\setminus D$ (monotone operator theory); it is continuous in $\R^N$ and $|u|\le C\,\dist(\cdot,\partial D)^s$ in $D$ \cite[Theorem~4.4]{reg} (for $1<p<2$ see also \cite{IM}). (For $c\not\equiv0$ one first obtains an $L^\infty$ bound by comparison with the solution for $c=0$, $f$ replaced by $|f|$, and then treats $c\Psi_p(u)$ as part of a bounded datum.)
\item[\textbf{(K4)}] \emph{\bf Interior H\"older estimates and weak Harnack inequality \changed{(\cite{dkp,dkp2,reg}):}} Locally bounded weak solutions of \eqref{eq:gen} are locally H\"older continuous (for bounded data see \cite{Cozzi}), with estimates depending on $\|u\|_{L^\infty}$ in a ball, the tail $\int|u|^{p-1}(1+|y|)^{-N-sp}dy$, $\|f\|_\infty$ and $\|c\|_\infty$; by scaling, if $\pl u=f$ in $B_r(x)$ with $f$ bounded, then $[u]_{C^{\gamma_1}(B_{r/2}(x))}\le Cr^{-\gamma_1}\big(\|u\|_{L^\infty(B_r(x))}+\mathrm{Tail}(u;x,r)+r^{\frac{sp}{p-1}}\|f\|_\infty^{\frac1{p-1}}\big)$, with $\mathrm{Tail}(u;x,r):=\big(r^{sp}\int_{\R^N\setminus B_r(x)}|u|^{p-1}|y-x|^{-N-sp}dy\big)^{\frac1{p-1}}$ and $\gamma_1=\gamma_1(N,s,p)\in(0,1)$. If $u\ge0$ in $\R^N$ is a weak supersolution of $\pl u\ge0$ in $B_{4r}(x_0)$, then for some $t>0$ and $c_H>0$ depending only on $N,s,p$,
\[
\inf_{B_{2r}(x_0)}u\ \ge\ c_H\Big(\frac1{|B_{2r}|}\int_{B_{2r}(x_0)}u^t\Big)^{1/t}.
\]
The local H\"older estimate also follows from (K8) with the bounded frozen datum. For the weak Harnack inequality, \cite[Theorem~1.2 and Lemma~4.1]{dkp2} assume $sp<N$, which is used only through the fractional Sobolev embedding with exponent $p^*=\frac{Np}{N-sp}$ (see \cite[(3.9) and (5.10)]{dkp2}). When $N\le sp$, we apply that embedding with any fixed $0<\sigma<\min\{s,N/p\}$ instead of $s$, using $[w]_{W^{\sigma,p}(B_r)}\le(2r)^{s-\sigma}[w]_{W^{s,p}(B_r)}$; the scaling $r^{\sigma p}r^{(s-\sigma)p}=r^{sp}$ is unchanged, the logarithmic and Caccioppoli estimates do not depend on the dimension, and the iteration gives the weak Harnack inequality for some $t>0$, which is all we use.
\item[\textbf{(K5)}] \emph{\bf Stability\changed{:}} If $u_j$ are viscosity solutions of $\pl u_j=F_j(x,u_j)$ in $D$, $u_j\to u$ locally uniformly in $D$ and in $\Ls$, and $F_j\to F$ locally uniformly, then $u$ is a viscosity solution of $\pl u=F(x,u)$. {In all our applications, the} $u_j$ and $u$ are continuous and locally bounded, and $F_j(x,u_j(x))$ are locally uniformly bounded. Then (K5) can be obtained from (K1): \mathred{on each $D'\Subset D$, the modification of Remark~\ref{rem:K1} gives bounded weak solutions with frozen data $F_j(x,u_j)+E_j$. Local uniform convergence and convergence in $\Ls$ imply local uniform convergence of the errors $E_j$ from Lemma~\ref{lem:loc}. The Caccioppoli inequality with tail \cite[Theorem~1.4]{dkp} (whose proof applies to bounded right-hand sides) gives uniform local energy bounds. Almost-everywhere convergence of the {normalised} increments identifies the weak limit of their nonlinear fluxes, so the weak formulation passes to the limit. Applying (K1) and undoing the modification with Lemma~\ref{lem:loc} gives the claim for $u$.}
\item[\textbf{(K6)}] \emph{\bf \changed{Hopf's} lemma \changed{(\cite{hopf}):}} Let $B=B_R(z)$, $\Lambda\ge0$, and let $v\ge0$ in $\R^N$ be a viscosity supersolution of $\pl v\ge -\Lambda \Psi_p(v)$ in $B$ with $v>0$ in $B$ and $v$ continuous in a neighbourhood of $\overline B$. Then $v(x)\ge c\,(R-|x-z|)^s$ in $B$, with $c>0$ depending only on $N,s,p,R,\Lambda$ and $\inf_{B_{\vartheta R}(z)}v$, where $\vartheta=\vartheta(N,s,p)\in(0,1)$ (and $c$ is nondecreasing with respect to the last quantity). This quantitative form follows from the proof of \cite[Lemma~4.1]{hopf}, which gives \cite[Theorem~1.5]{hopf}: there $v$ is compared, in $\{x\in B:\ \dist(x,\partial B)<\rho_H\}$ with $\rho_H=\rho_H(N,s,p,R)$, with $\varepsilon(\delta^s+\alpha\chi_K)$, where $\delta=\dist(\cdot,\R^N\setminus B)$, $\pl\delta^s$ is bounded there by \cite[Theorem~3.6]{reg}, and $K\subset B_{R-\rho_H}(z)$ is a ball. Using (K2) with $c\equiv\Lambda$, the constant $\alpha$ depends only on $N,s,p,R,\Lambda$, and $\varepsilon$ only on these and $\inf_{B_{R-\rho_H}(z)}v$; by scaling, $R-\rho_H=\vartheta R$. The viscosity setting is reduced to the weak one by (K1) and Remark~\ref{rem:K1}.
\item[\textbf{(K7)}] \emph{\bf Exterior barrier \changed{(\cite{reg,hopf}):}} There are $C_*>0$ and a continuous function $w_*$ with $w_*=0$ in $\overline{B_1}$, $0<w_*(x)\le C_*(|x|-1)^s$ for $1<|x|<2$, $1/C_*\le w_*\le C_*$ for $|x|\ge2$, and $\pl w_*\ge 1$ in $B_2\setminus\overline{B_1}$ in the viscosity sense. Such a function is obtained from the solution of $\pl\theta_0=1$ in $B_4\setminus\overline{B_1}$, $\theta_0=0$ outside, given by (K3) (a local version is \cite[Lemma~4.3]{reg}), by raising it to $\max\{\theta_0,m\}$ outside $B_3$ for a constant $m\in(0,\min_{\partial B_3}\theta_0)$ small: since this changes the function by at most $m$ at distance $\ge1$ from $B_2$, by Lemma~\ref{lem:loc} and \eqref{eq:elem} the equation in $B_2\setminus\overline{B_1}$ changes by at most $C(\|\theta_0\|_\infty)\,m^{\min\{1,p-1\}}\le\frac12$, and a multiple $2^{1/(p-1)}$ of the result is the required $w_*$.
Items (K3), (K6) and (K7) hold for every $N\ge2$: \cite{reg} only assumes $N>1$ and does not use Sobolev inequalities, and the proof of \cite[Lemma~4.1]{hopf} only uses the barrier of \cite[Theorem~3.6]{reg} and weak comparison.
\item[\textbf{(K8)}] \emph{\bf Higher H\"older regularity \changed{(\cite[Theorem~1.4]{BLS}; \cite[Theorem~1.2]{GL}):}} Let $\changed{p>1}$ and let $u\in L^\infty(B_1)\cap\Ls$ be a local weak solution of $\pl u=f$ in a neighbourhood of $\overline{B_1}$ with $f\in L^\infty(B_1)$. Then $u\in C^\alpha(B_{1/32})$ for every $\changed{0<}\alpha<\min\{1,\frac{sp}{p-1}\}$, and
\[
\begin{gathered}
[u]_{C^\alpha(B_{1/32})}\le C\Big(\|u\|_{L^\infty(B_{\changed{1/2}})}+\Big(\int_{\R^N\setminus B_{\changed{1/2}}}\frac{|u(y)|^{p-1}}{|y|^{N+sp}}dy\Big)^{\frac1{p-1}}+\|f\|_{L^\infty(B_1)}^{\frac1{p-1}}\Big),\\
C=C(N,s,p,\alpha)
\end{gathered}
\]
\mathred{For $p\ge2$,} this is \cite[Theorem 1.4]{BLS} with $x_0=0$ and $R=\frac14$, so that $B_{4R}=B_1$\mathred{; the tail on $B_{1/2}\setminus B_{1/4}$ is absorbed into the local bound. For $1<p<2$, use \cite[Theorem~1.2]{GL} with the same radius and $q=\infty$.}
\end{itemize}

\medskip
\noindent{\bf Scaling and translations:} If $u_{t,z}(x):=u(z+tx)$ with $t>0$, then $\pl u_{t,z}(x)=t^{sp}(\pl u)(z+tx)$, by the change of variables $y\mapsto z+ty$ in \eqref{eq:pl}; the same holds in the viscosity sense, since $\varphi$ touches $u$ at $z+tx$ if and only if $\varphi(z+t\,\cdot)$ touches $u_{t,z}$ at $x$, and admissibility is preserved. Moreover $\pl(\kappa u)=\kappa^{p-1}\pl u$ and $\pl(u+\kappa)=\pl u$ for constants $\kappa>0$.

\subsection{Elementary lemmas}

\begin{lemma}[Localization]\label{lem:loc}
Let $D\subset\R^N$ be open, $F:D\times\R\to\R$ continuous, and let $u$ be a viscosity subsolution \textup{(}resp.\ supersolution\textup{)} of $\pl u=F(x,u)$ in $D$, \rev{with $u\in\Ls$ finite in $D$}. Let $D'\subset D$ be open and let $\tilde u$ satisfy $\tilde u=u$ in $D'$ and \rev{$\tilde u\in\Ls$}. Set $\Sigma:=\{y\in\R^N\setminus D':\ \tilde u(y)\neq u(y)\}$ and, for $x\in D'$,
\[
E(x):=\int_{\Sigma}\frac{\Psi_p(u(x)-\tilde u(y))-\Psi_p(u(x)-u(y))}{|x-y|^{N+sp}}\,dy .
\]
Then $\tilde u$ is a viscosity subsolution \textup{(}resp.\ supersolution\textup{)} of $\pl\tilde u= F(x,\tilde u)+E(x)$ in $D'$, in the sense that the inequality of Definition~\ref{def:visc} holds with $F(x_0,\tilde u(x_0))+E(x_0)$ on the right-hand side (when $E$ is continuous, this is the usual notion). Moreover:
\begin{itemize}
\item[(i)] if $\tilde u\le u$ in $\R^N$ then $E\ge0$, and if $\tilde u\ge u$ in $\R^N$ then $E\le0$. Consequently, if $\tilde u\ge u$ \textup{(}resp.\ $\tilde u\le u$\textup{)}, then $\tilde u$ is a viscosity subsolution \textup{(}resp.\ supersolution\textup{)} of $\pl\tilde u=F(x,\tilde u)$ in $D'$; \rev{\changed{this} monotone implication also holds under only the one-sided tail assumptions of Definition~\ref{def:visc}, provided the functions are finite in $D'$; no subtraction of integrals is needed in that case.}
\item[(ii)] if $D''\subset\R^N$ is bounded, $u$ is continuous and bounded in $D''$ and $\dist(D'',\Sigma)\ge\delta>0$, then $E$ (defined by the same formula for $x\in D''$) is continuous in $D''$ and
\[
\sup_{D''}|E|\le C\Big(\sup_{D''}|u|^{p-1}+\int_{\Sigma}\frac{|u(y)|^{p-1}+|\tilde u(y)|^{p-1}}{(1+|y|)^{N+sp}}\,dy\Big),
\]
with $C$ depending only on $N,s,p,\delta$ and $R_0:=\sup_{x\in D''}|x|$.
\end{itemize}
\end{lemma}

\begin{proof}
{We consider subsolutions first; the supersolution case is analogous.} Let $x\in D'$, $B_r(x)\subset D'$, and let $\varphi\in C^2(B_r(x))$ be admissible and touch $\tilde u$ from above at $x$ in $B_r(x)$. Since $\tilde u=u$ in $B_r(x)$, $\varphi$ also touches $u$ from above at $x$, and since $B_r(x)\subset D$,
\[
\pl\varphi^{u}_r(x)\le F(x,u(x)),\qquad \varphi^{u}_r:=\varphi\ \text{in }B_r(x),\quad\varphi^{u}_r:=u\ \text{outside}.
\]
Let $\varphi^{\tilde u}_r$ be defined in the same way with $\tilde u$ in place of $u$. The two functions coincide in $B_r(x)$, and outside $B_r(x)$ they differ only on $\Sigma$ (note $\Sigma\cap B_r(x)=\emptyset$). Since $\varphi(x)=u(x)$, subtracting the two integrals defining $\pl$ at $x$ (the principal values in $B_r(x)$ cancel) gives
\[
\pl\varphi^{\tilde u}_r(x)-\pl\varphi^{u}_r(x)=\int_{\Sigma}\frac{\Psi_p(u(x)-\tilde u(y))-\Psi_p(u(x)-u(y))}{|x-y|^{N+sp}}\,dy=E(x),
\]
hence $\pl\varphi^{\tilde u}_r(x)\le F(x,\tilde u(x))+E(x)$, as $\tilde u(x)=u(x)$.

(i) If $\tilde u\le u$, then $u(x)-\tilde u(y)\ge u(x)-u(y)$ and, $\Psi_p$ being increasing, the integrand of $E$ is nonnegative; the case $\tilde u\ge u$ is symmetric. For a subsolution and $\tilde u\ge u$ we get $E\le0$, so $\pl\varphi^{\tilde u}_r(x)\le F(x,\tilde u(x))$; for a supersolution and $\tilde u\le u$ we get $E\ge0$ and $\pl\varphi^{\tilde u}_r(x)\ge F(x,\tilde u(x))$. \rev{Under the one-sided tail assumptions, {the conclusion follows directly} by ordering the two patched test functions outside $B_r(x)$ and using the monotonicity of $\Psi_p$, {without forming their difference}.}

(ii) Let $x\in D''$ and $y\in\Sigma$. Then $|x-y|\ge\delta$, and $|x-y|\ge|y|-R_0\ge\frac12(1+|y|)$ if $|y|\ge 2R_0+1$; if $|y|<2R_0+1$ then $|x-y|\ge\delta\ge\frac{\delta}{2R_0+2}(1+|y|)$. Hence $|x-y|\ge c_\delta(1+|y|)$, with $c_\delta:=\min\{\frac12,\frac\delta{2R_0+2}\}$. By $|\Psi_p(a)-\Psi_p(b)|\le|a|^{p-1}+|b|^{p-1}$ and \eqref{eq:elem},
\[
\frac{|\Psi_p(u(x)-\tilde u(y))-\Psi_p(u(x)-u(y))|}{|x-y|^{N+sp}}\le 2^{p-1}c_\delta^{-N-sp}\,\frac{2|u(x)|^{p-1}+|\tilde u(y)|^{p-1}+|u(y)|^{p-1}}{(1+|y|)^{N+sp}} ,
\]
{which yields the bound}, since $\int_{\R^N}(1+|y|)^{-N-sp}dy<\infty$. If $x_k\to x$ in $D''$, the integrands converge pointwise in $\Sigma$ (because $u$ is continuous at $x$ and $|x-y|\ge\delta$), and they are dominated by the integrable function above with $|u(x)|$ replaced by $\sup_{D''}|u|$; by dominated convergence $E(x_k)\to E(x)$.
\end{proof}

\begin{remark}[Use of (K1) for locally bounded functions]\label{rem:K1}
Let $D$ be bounded with $\overline D\subset\R^N\setminus\{0\}$, let $D_\delta$ be its $\delta$-neighbourhood, with $\overline{D_{2\delta}}\subset\R^N\setminus\{0\}$, and fix $\zeta\in C(\R^N)$ with $0\le\zeta\le1$, $\zeta=1$ on $D_\delta$ and $\zeta=0$ outside $D_{2\delta}$.
\begin{itemize}
\item If $v\ge0$, $v\in\Ls$, is continuous in $\R^N\setminus\{0\}$ and a viscosity supersolution of $\pl v\ge f$ in $D$, let $M'>\max_{\overline{D_\delta}}v$ and $\tilde v:=\min\{v,M'\zeta\}$. Then $\tilde v$ is nonnegative, bounded by $M'$, continuous in $\R^N$ (it vanishes outside $D_{2\delta}$, in particular near $0$), $\tilde v\le v$, and $\tilde v=v$ in $D_\delta$ (where $\zeta=1$ and $v<M'$), so by Lemma~\ref{lem:loc}(i) it is a bounded continuous viscosity supersolution of $\pl\tilde v\ge f$ in $D$, hence a weak supersolution by (K1).
\item If $u$ is continuous in $\R^N\setminus\{0\}$, $u\in\Ls$, and a viscosity solution of $\pl u=f$ in $D$, then $\tilde u:=\zeta u$ is bounded and continuous, $\tilde u=u$ in $D_\delta$, and by Lemma~\ref{lem:loc} (applied twice, as a sub- and as a supersolution) $\tilde u$ is a viscosity solution of $\pl\tilde u=f+E$ in $D$, where $E$ is bounded and continuous in a neighbourhood of $\overline D$ by Lemma~\ref{lem:loc}(ii) (with $D''=D_{\delta/2}$, which is at distance $\ge\delta/2$ from $\Sigma\subset\R^N\setminus D_\delta$).
\end{itemize}
{The same argument applies, after translation, when} the only possible singular point of the functions is a point $\xi_0\notin\overline{D_{2\delta}}$ instead of $0$.
\end{remark}

\begin{lemma}[Strong maximum principle]\label{lem:smp}
Let $D$ be open, $c\in C(D)$, and let $u\ge0$ in $\R^N$ be a viscosity supersolution of $\pl u\ge -c(x)\Psi_p(u)$ in $D$ which is positive on a set of positive measure \textup{(}for instance, $u$ continuous outside a closed null set and $u\not\equiv0$ there\textup{)}. Then $u>0$ in $D$.
\end{lemma}

\begin{proof}
Assume $u(x_0)=0$ for some $x_0\in D$ and let $r>0$ with $B_r(x_0)\subset D$. The function $\varphi(y)=-|y-x_0|^{q_0}$ is of class $C^2$ (as $q_0\ge2$), satisfies $\varphi\le0\le u$ and $\varphi(x_0)=u(x_0)$, and $\nabla\varphi(y)=-q_0|y-x_0|^{q_0-2}(y-x_0)$ vanishes only at $x_0$. It belongs to $C^2_{q_0}(B_r(x_0))$: indeed $\dist(y,\mathcal N_\varphi)=|y-x_0|$, $|\nabla\varphi(y)|=q_0|y-x_0|^{q_0-1}$ and $|D^2\varphi(y)|\le C|y-x_0|^{q_0-2}$; and $q_0>\frac{sp}{p-1}$. So $\varphi$ is admissible, and since $\Psi_p(u(x_0))=0$,
\[
0\le \pl\varphi_r(x_0)=\int_{B_r(x_0)}\frac{|y-x_0|^{q_0(p-1)}}{|x_0-y|^{N+sp}}dy-\int_{\R^N\setminus B_r(x_0)}\frac{u(y)^{p-1}}{|x_0-y|^{N+sp}}dy .
\]
The first integral equals $\mathcal H^{N-1}(\Sp)\,\frac{r^{q_0(p-1)-sp}}{q_0(p-1)-sp}\to0$ as $r\to0$, since $q_0(p-1)>sp$. {By the monotone convergence theorem, }{the second integral tends to} $\int_{\R^N}u^{p-1}|x_0-y|^{-N-sp}dy\in(0,\infty]$, which is positive because $u>0$ on a set of positive measure. {For sufficiently small} $r$, {this is a contradiction.}
\end{proof}

\begin{lemma}[Concave transformations]\label{lem:power}
Let $D$ be open, $F\in C(D)$, and let $v\in C(\R^N\setminus\{0\})\cap\Ls$, $v\ge0$ in $\R^N$, $v>0$ in $D$, be a viscosity supersolution of $\pl v\ge F$ in $D$. Let $a\in(0,1)$, $z:=v^a$, and for $x\in D$, $y\in\R^N$,
\[
\Gamma_a(x,y):=\Psi_p\big(z(x)-z(y)\big)-a^{p-1}v(x)^{(a-1)(p-1)}\,\Psi_p\big(v(x)-v(y)\big).
\]
Then $\Gamma_a\ge0$, $\Gamma_a(x,y)>0$ if $v(y)\neq v(x)$, {and if} $\varphi$ is admissible and touches $z$ from below at $x\in D$ in $B_r(x)\subset D$, then
\begin{equation}\label{eq:power}
\pl\varphi_r(x)\ \ge\ a^{p-1}v(x)^{(a-1)(p-1)}F(x)+\int_{\R^N\setminus B_r(x)}\frac{\Gamma_a(x,y)}{|x-y|^{N+sp}}\,dy .
\end{equation}
\end{lemma}

\begin{proof}
Let $g(t)=t^a$, $t\ge0$. For $t>0$ and $\tau\ge0$, concavity gives
\begin{equation}\label{eq:concave}
g(t)-g(\tau)\ge g'(t)(t-\tau),
\end{equation}
{with strict inequality whenever} $\tau\neq t$ (strict concavity). Since $\Psi_p$ is increasing and $\Psi_p(\lambda\sigma)=\lambda^{p-1}\Psi_p(\sigma)$ for $\lambda>0$, applying \eqref{eq:concave} with $t=v(x)$, $\tau=v(y)$ we get $\Psi_p(z(x)-z(y))\ge g'(v(x))^{p-1}\Psi_p(v(x)-v(y))$, which is $\Gamma_a\ge0$ because $g'(v(x))^{p-1}=a^{p-1}v(x)^{(a-1)(p-1)}$; the inequality is strict if $v(y)\ne v(x)$.

Now let $\varphi$ be as in the statement. Since $\varphi(x)=z(x)>0$, there is $r''\in(0,r]$ with $\varphi>0$ in $B_{r''}(x)$. By Lemma~\ref{lem:grad} below, $\pl\varphi_r(x)\ge\pl\varphi_{r''}(x)$, and since $\Gamma_a\ge0$ the integral in \eqref{eq:power} over $\{|y-x|\ge r''\}$ is not smaller than the one over $\{|y-x|\ge r\}$. {It therefore suffices to prove} \eqref{eq:power} with $r''$ in place of $r$, that is, we may assume $\varphi>0$ in $B_r(x)$. Then $\phi:=\varphi^{1/a}\in C^2(B_r(x))$, $\phi\le z^{1/a}=v$ in $B_r(x)$ and $\phi(x)=v(x)$, i.e.\ $\phi$ touches $v$ from below at $x$. Since $\nabla\phi=\frac1a\varphi^{\frac1a-1}\nabla\varphi$, $\phi$ and $\varphi$ have the same critical points. If $x$ is an isolated critical point and $\varphi\in C^2_\alpha$ with $\alpha>\frac{sp}{p-1}$ (this only matters in the singular range $p\le\frac2{2-s}$, where $\alpha>\frac{sp}{p-1}\ge2$), then $|\nabla\phi|\ge c|\nabla\varphi|$ near $x$, and
\[
D^2\phi=\tfrac1a\varphi^{\frac1a-1}D^2\varphi+\tfrac1a\big(\tfrac1a-1\big)\varphi^{\frac1a-2}\nabla\varphi\otimes\nabla\varphi ,
\]
where $|\nabla\varphi(y)|\le\int_0^1|D^2\varphi(x+\tau(y-x))||y-x|\,d\tau\le C|y-x|^{\alpha-1}$, so that $|\nabla\varphi|^2\le C|y-x|^{2\alpha-2}\le C|y-x|^{\alpha-2}$ near $x$; hence $\phi\in C^2_\alpha$ near $x$. Thus $\phi$ is admissible and, with $\phi_r:=\phi$ in $B_r(x)$, $\phi_r:=v$ outside,
\[
\pl\phi_r(x)\ge F(x).
\]
We have $\varphi_r=g(\phi_r)$ in $\R^N$ (in $B_r(x)$ by definition, outside because $z=g(v)$), and $\phi_r\ge0$. By \eqref{eq:concave} with $t=\phi_r(x)=v(x)$, $\tau=\phi_r(y)$, for every $y$,
\[
\Psi_p(\varphi_r(x)-\varphi_r(y))=g'(v(x))^{p-1}\Psi_p(\phi_r(x)-\phi_r(y))+\Gamma(x,y),\qquad \Gamma(x,y)\ge0,
\]
and $\Gamma(x,y)=\Gamma_a(x,y)$ for $|y-x|\ge r$. We divide by $|x-y|^{N+sp}$ and integrate over $\{|y-x|>\epsilon\}$. As $\epsilon\to0$, {the integrals} of the first two terms converge to $\pl\varphi_r(x)$ and $g'(v(x))^{p-1}\pl\phi_r(x)$ (both principal values exist since $\varphi$, $\phi$ are admissible), and the integral of $\Gamma\ge0$ increases to $\int\Gamma|x-y|^{-N-sp}dy\ge\int_{|y-x|\ge r}\Gamma_a|x-y|^{-N-sp}dy$. Hence \eqref{eq:power}.
\end{proof}

\begin{lemma}[Test functions for homogeneous functions; monotonicity in $r$]\label{lem:grad}
Let $W\in C(\R^N\setminus\{0\})$ satisfy $W(tx)=t^{-\gamma}W(x)$ for $t>0$, with $\gamma\ne0$, and let $x\ne0$ with $W(x)\ne0$. If $\varphi\in C^1$ touches $W$ from below \textup{(}or from above\textup{)} at $x$, then $\nabla\varphi(x)\cdot x=-\gamma W(x)\neq0$. {In particular, every such test function has} $\nabla\varphi(x)\neq0$. Moreover, if $w$ is any function, $\varphi$ touches $w$ from below at $x$ in $B_r(x)$ and $0<r'<r$, then $\pl\varphi_r(x)\ge\pl\varphi_{r'}(x)$\changed{ whenever both operators are defined}.
\end{lemma}

\begin{proof}
The $C^1$ function $\tau\mapsto\varphi(\tau x)-W(\tau x)=\varphi(\tau x)-\tau^{-\gamma}W(x)$ is $\le0$ (or $\ge0$) for $\tau$ near $1$ and vanishes at $\tau=1$, {hence its derivative at} $\tau=1$, which is $\nabla\varphi(x)\cdot x+\gamma W(x)$, vanishes. For the last claim, $\varphi_r=\varphi_{r'}=\varphi$ in $B_{r'}(x)$, $\varphi_r=\varphi\le w=\varphi_{r'}$ in $B_r(x)\setminus B_{r'}(x)$, and $\varphi_r=\varphi_{r'}=w$ outside $B_r(x)$. So $\varphi_r(x)-\varphi_r(y)\ge\varphi_{r'}(x)-\varphi_{r'}(y)$ for all $y$, with equality in $B_{r'}(x)$, and $\Psi_p$ is increasing.
\end{proof}

{The next lemma gives a quantitative bound for the principal value over small balls;} it is used for the bump function below and at the end of Section~\ref{sec:exist}.

\begin{lemma}[Principal values in small balls]\label{lem:small}
Let $x\in\R^N$, $\rho>0$ and $F\in C^2(B_\rho(x))$ with $F(x)=0$, $G:=|\nabla F(x)|>0$, $H:=\sup_{B_\rho(x)}|D^2F|<\infty$ and $\rho H\le G$. Then the limit
\[
P:=\lim_{\epsilon\to0}\int_{\epsilon<|y-x|<\rho}\frac{\Psi_p(F(y))}{|x-y|^{N+sp}}\,dy
\]
exists and $|P|\le C_0\,G^{p-2}H\rho^{p-sp}$, with $C_0$ depending only on $N,s,p$. In particular, if $\varphi\in C^2$ near $x$ and $\nabla\varphi(x)\ne0$, then $\pl\varphi_r(x)$ is well defined \mathred{when $\varphi$ is patched outside the contact ball with a function in $\Ls$}, and the contribution of $B_\rho(x)$ to it tends to $0$ as $\rho\to0$.
\end{lemma}

\begin{proof}
Write $z=y-x$ and {symmetrize} ($z\mapsto-z$):
\[
\int_{\epsilon<|z|<\rho}\frac{\Psi_p(F(x+z))}{|z|^{N+sp}}dz=\frac12\int_{\epsilon<|z|<\rho}\frac{\Psi_p(F(x+z))+\Psi_p(F(x-z))}{|z|^{N+sp}}dz .
\]
By Taylor's formula, $F(x\pm z)=\pm\ell(z)+Q_\pm(z)$ with $\ell(z):=\nabla F(x)\cdot z$ and $|Q_\pm(z)|\le\frac12H|z|^2$. As $\Psi_p$ is odd, the numerator equals $\mathcal D(z):=\Psi_p(a)-\Psi_p(b)$ with $a:=\ell+Q_+$, $b:=\ell-Q_-$, and $|a-b|\le H|z|^2$. We show that $\int_{B_\rho}|\mathcal D(z)||z|^{-N-sp}dz\le 2C_0G^{p-2}H\rho^{p-sp}$; then the limit exists {by the dominated convergence theorem} and hence the bound follows.

\emph{Case $p\ge2$:} Since $|\ell|\le G|z|$ and $H|z|\le H\rho\le G$, we have $|a|+|b|\le 2G|z|+H|z|^2\le3G|z|$, and by \eqref{eq:elem}, $|\mathcal D(z)|\le(p-1)3^{p-2}G^{p-2}H|z|^p$.

\emph{Case $1<p<2$:} Write $z=r\sigma$, $\sigma\in\Sp$, $e:=\nabla F(x)/G$, so $\ell=rG\,(e\cdot\sigma)$. If $|e\cdot\sigma|\ge Hr/G$, then $|\ell|\ge Hr^2$, so $a$ and $b$ have the sign of $\ell$ and $|a|,|b|\ge|\ell|-\frac12Hr^2\ge\frac12|\ell|$; since $|\Psi_p'(\tau)|=(p-1)|\tau|^{p-2}$ is decreasing in $|\tau|$, the mean value theorem gives $|\mathcal D|\le(p-1)(|\ell|/2)^{p-2}Hr^2\le2\,(rG|e\cdot\sigma|)^{p-2}Hr^2$. If $|e\cdot\sigma|<Hr/G$, then by \eqref{eq:elem} $|\mathcal D|\le2^{2-p}(Hr^2)^{p-1}$. Let $c_{N,p}:=\int_{\Sp}|e\cdot\sigma|^{p-2}d\sigma<\infty$ (as $p-2>-1$) and note that $\mathcal H^{N-1}(\{\sigma:|e\cdot\sigma|<\kappa\})\le c_N\kappa$ for $\kappa>0$ (for $N=2$ the measure is $4\arcsin\min\{\kappa,1\}\le2\pi\kappa$; for $N\ge3$ the density of $e\cdot\sigma$ is bounded). Hence
\[
\int_{\Sp}|\mathcal D(r\sigma)|\,d\sigma\le 2c_{N,p}G^{p-2}Hr^p+2c_N(Hr^2)^{p-1}\frac{Hr}{G}\le C\,G^{p-2}Hr^p,
\]
where we used $(Hr^2)^{p-1}\frac{Hr}{G}=(Hr)^{p-1}\frac{Hr^p}{G}\le G^{p-2}Hr^p$, by $Hr\le G$.

In both cases, $\int_{B_\rho}|\mathcal D||z|^{-N-sp}dz\le C\,G^{p-2}H\int_0^\rho r^{p-1-sp}dr=\frac{C}{p-sp}G^{p-2}H\rho^{p-sp}$, and $p>sp$. For {the final claim, apply the result to} $F=\varphi-\varphi(x)$ (note $\Psi_p(\varphi(x)-\varphi(y))=-\Psi_p(F(y))$), with $G=|\nabla\varphi(x)|$ fixed, $H$ bounded and $\rho\to0$.
\end{proof}

\begin{lemma}[A bump function]\label{lem:bump}
Let $k\ge3$ with $k(p-1)>sp$, $r_0\in(0,1]$, and $\eta(y):=\big(1-(|y|/r_0)^{q_0}\big)_+^k$. Then $\eta\in C^2(\R^N)$, $0\le\eta\le1=\eta(0)$, $\nabla\eta\ne0$ in $B_{r_0}\setminus\{0\}$, $\eta\in C^2_{q_0}(B_{r_0/2})$, $\pl\eta(x)$ is well defined for every $x\in B_{r_0}$, and
\begin{align}\label{meta}
M_\eta:=\sup_{x\in B_{r_0}}|\pl\eta(x)|=r_0^{-sp}M_1<\infty,\qquad M_1=M_1(N,s,p,k).
\end{align}
In particular $\eta(\cdot-\hat x)$ is an admissible test function at every point of $B_{r_0}(\hat x)$.
\end{lemma}

\begin{proof}
%\emph{Scaling.}
 If $\eta_1$ denotes the function with $r_0=1$, then $\eta=\eta_1(\cdot/r_0)$ and $\pl\eta(x)=r_0^{-sp}(\pl\eta_1)(x/r_0)$. So we take $r_0=1$ and write $\eta=\eta_1$, $\sigma(y):=1-|y|^{q_0}$.

%\emph{Regularity.}
 $|y|^{q_0}\in C^2(\R^N)$ since $q_0\ge2$, and $\tau\mapsto\tau_+^k\in C^2(\R)$ since $k\ge3$; hence $\eta\in C^2$. For $|y|<1$,
\[
\nabla\eta(y)=-kq_0\,\sigma^{k-1}|y|^{q_0-2}y,\qquad |D^2\eta(y)|\le C\big(\sigma^{k-1}|y|^{q_0-2}+\sigma^{k-2}|y|^{2q_0-2}\big),
\]
with $C=C(N,k,q_0)$. Thus $\nabla\eta(y)\ne0$ for $0<|y|<1$. For $|y|\le\frac12$, $\sigma\ge1-2^{-q_0}$, so $|\nabla\eta(y)|\ge kq_0(1-2^{-q_0})^{k-1}|y|^{q_0-1}$ and $|D^2\eta(y)|\le C|y|^{q_0-2}$; since $\dist(y,\mathcal N_\eta)=|y|$ in $B_{1/2}$, $\eta\in C^2_{q_0}(B_{1/2})$.

%\emph{Two pointwise bounds.}
 For all $y\in\R^N$,
\begin{equation}\label{eq:bump-pw}
0\le1-\eta(y)\le\min\{k|y|^{q_0},1\},\qquad \eta(y)\le q_0^k(1-|y|)_+^k .
\end{equation}
{The first inequality follows from} $1-(1-\tau)^k\le k\tau$ for $\tau\in[0,1]$, {and the second from} $1-|y|^{q_0}\le q_0(1-|y|)$ for $|y|\le1$.

\changed{At} $x=0$, since $q_0(p-1)>sp$,
\[
\pl\eta(0)=\int_{\R^N}\frac{(1-\eta(y))^{p-1}}{|y|^{N+sp}}dy\le\int_{\R^N}\frac{\min\{k|y|^{q_0},1\}^{p-1}}{|y|^{N+sp}}dy<\infty .
\]

For $x\ne0$ we split $\pl\eta(x)$ into the integrals over $B_\rho(x)$ ({inner part}) and over $\R^N\setminus B_\rho(x)$ ({outer part}), with $\rho$ chosen below. The inner part is $-P$, where $P$ is given by Lemma~\ref{lem:small} with $F=\eta-\eta(x)$; {the outer part is bounded by} $\int_{|z|\ge\rho}|\eta(x)-\eta(x+z)|^{p-1}|z|^{-N-sp}dz$.

\emph{Case $0<|x|\le\frac12$:} Let $\rho:=c_1|x|$ with $c_1\in(0,\frac12]$. For $y\in B_\rho(x)$ we have $\frac12|x|\le|y|\le\frac32|x|\le\frac34$, so $H:=\sup_{B_\rho(x)}|D^2\eta|\le C_A|x|^{q_0-2}$, while $G:=|\nabla\eta(x)|\in[g_A|x|^{q_0-1},kq_0|x|^{q_0-1}]$, $g_A:=kq_0(1-2^{-q_0})^{k-1}$. Choosing $c_1:=\min\{\frac12,g_A/C_A\}$ we get $\rho H\le G$. By Lemma~\ref{lem:small} (using the upper bound for $G$ if $p\ge2$ and the lower bound if $p<2$),
\[
|\text{inner part}|\le C\,|x|^{(q_0-1)(p-2)}\,|x|^{q_0-2}\,|x|^{p-sp}=C\,|x|^{q_0(p-1)-sp}\le C .
\]
For $|z|\ge\rho$, by \eqref{eq:bump-pw}, $|\eta(x)-\eta(x+z)|\le k(|x|^{q_0}+|x+z|^{q_0})\le2k(|x|+|z|)^{q_0}\le2k(1+c_1^{-1})^{q_0}|z|^{q_0}$, and also $|\eta(x)-\eta(x+z)|\le1$. Hence {the outer part  is bounded by} $\int_{\R^N}\min\{C|z|^{q_0},1\}^{p-1}|z|^{-N-sp}dz<\infty$.

\emph{Case $\frac12<|x|<1$:} Let $t:=1-|x|\in(0,\frac12)$ and $\rho:=c_2t$ with $c_2\in(0,\frac12]$. For $y\in B_\rho(x)$, $\frac t2\le1-|y|\le\frac{3t}2$ and $|y|\ge\frac14$, so $\sigma(y)\in[1-|y|,q_0(1-|y|)]\subset[\frac t2,\frac{3q_0t}{2}]$ and $H\le C_Bt^{k-2}$. Also $\sigma(x)\in[t,q_0t]$ and $|x|^{q_0-1}\in[2^{1-q_0},1]$, so $G=kq_0\sigma(x)^{k-1}|x|^{q_0-1}\in[g_Bt^{k-1},G_Bt^{k-1}]$. Choosing $c_2:=\min\{\frac12,g_B/C_B\}$ we get $\rho H\le G$, and Lemma~\ref{lem:small} gives
\[
|\text{inner part}|\le C\,t^{(k-1)(p-2)}\,t^{k-2}\,t^{p-sp}=C\,t^{k(p-1)-sp}\le C .
\]
For $|z|\ge\rho$, by \eqref{eq:bump-pw} and $1-|x+z|\le t+|z|$, $|\eta(x)-\eta(x+z)|\le\eta(x)+\eta(x+z)\le2q_0^k(t+|z|)^k\le2q_0^k(1+c_2^{-1})^k|z|^k$, and it is also $\le1$; {hence the outer part} is bounded by $\int_{\R^N}\min\{C|z|^k,1\}^{p-1}|z|^{-N-sp}dz<\infty$, because $k(p-1)>sp$.

{These bounds are independent of} $x$, which proves $M_1<\infty$. The last assertion follows from $\nabla\eta\ne0$ in $B_{r_0}\setminus\{0\}$ and $\eta\in C^2_{q_0}$ near $0$, with $q_0>\frac{sp}{p-1}$.
\end{proof}

\begin{remark}[Touching with the bump function]\label{rem:touch}
The following construction is used several times. Let $\hat x\in\R^N$, $B:=B_{r_0}(\hat x)$, $\tilde\eta:=\eta(\cdot-\hat x)$, and let $w\in\Ls$, $w\ge0$ in $\R^N$, be lower semicontinuous and positive in $\overline B$. Let
\[
c^*:=\sup\{c\ge0:\ c\tilde\eta\le w\ \text{in }B\}=\inf_{y\in B}\frac{w(y)}{\tilde\eta(y)} .
\]
Then $0<\min_{\overline B}w\le c^*\le w(\hat x)$. The function $w/\tilde\eta$ is lower semicontinuous in $B$ and tends to $+\infty$ at $\partial B$ (where $\tilde\eta\to0$ and $w\ge\min_{\overline B}w>0$), so the infimum is attained at some $x_1\in B$. Thus $\phi:=c^*\tilde\eta$ touches $w$ from below at $x_1$, $\phi\le w$ in $\R^N$ (outside $B$, $\phi=0\le w$), and $\phi$ is admissible at $x_1$ by Lemma~\ref{lem:bump}. If $w$ is a viscosity supersolution of $\pl w\ge F\changed{(x,w)}$ in a neighbourhood of $\overline B$, then for small $r$, since $\phi_r\ge\phi$ with equality at $x_1$,
\begin{equation}\label{eq:touch}
F(x_1,w(x_1))\le\pl\phi_r(x_1)\le\pl\phi(x_1)=(c^*)^{p-1}\pl\eta(x_1-\hat x)\le (c^*)^{p-1}M_\eta .
\end{equation}
\end{remark}
\subsection{Homogeneous functions}
For $\beta\in\R$ {we set}
\[
H_\beta:=\big\{u\in C(\R^N\setminus\{0\}):\ u(tx)=t^{-\beta}u(x)\ \text{ for all } t>0,\ x\neq0\big\}.
\]
{Note that every} $u\in H_\beta$ is determined by $\hat u:=u|_{\Sp}$ through $u(x)=|x|^{-\beta}\hat u(x/|x|)$.

 If $u\in H_\beta$ is nonnegative and positive on some open cone $\mathcal C_{\omega'}$, then, in polar coordinates,
\[
\int_{\R^N}\frac{u^{p-1}}{(1+|x|)^{N+sp}}dx=\int_{\Sp}\hat u^{p-1}d\mathcal H^{N-1}\int_0^\infty\frac{t^{N-1-\beta(p-1)}}{(1+t)^{N+sp}}dt ,
\]
{which is finite if and only if} $-\frac{sp}{p-1}<\beta<\frac{N}{p-1}$ (integrability at $0$ and at $\infty$).

%\emph{Scaling.}
 Let $u\in H_\beta\cap\Ls$, $t>0$. If $\varphi$ touches $u$ from below (above) at $x$ in $B_r(x)$, then $\varphi^t(y):=t^{-\beta}\varphi(y/t)$ touches $u$ from below (above) at $tx$ in $B_{tr}(tx)$, admissibility is preserved, and $(\varphi^t)_{tr}=t^{-\beta}\varphi_r(\cdot/t)$, so that $\pl(\varphi^t)_{tr}(tx)=t^{-\beta(p-1)-sp}\pl\varphi_r(x)$. Hence, if $F(tx,t^{-\beta}\xi)=t^{-\beta(p-1)-sp}F(x,\xi)$ for all $t>0$, then $u\in H_\beta$ is a viscosity supersolution (subsolution) of $\pl u=F(x,u)$ in $\Cw$ if and only if it is so in $\Cw\cap\{\frac12<|x|<2\}$, and it is enough to check the defining inequality at points $x\in\Cw\cap\Sp$ (for test functions with small radius, {which suffices by the monotonicity in} $r$ of Lemma~\ref{lem:grad}). This applies to $F=0$, $F=\mu|x|^{-sp}\Psi_p(\xi)$ and $F=\lambda|x|^{-sp}(\psi^{p-1}+|x|^{-\beta(p-1)})$ with $\psi\in H_\beta$\mathred{, $\psi\ge0$}.

We denote by $\mathcal S_\beta(\omega)$ the set of positive homogeneous supersolutions,
\[
\mathcal S_\beta(\omega):=\big\{v\in H_\beta\cap\Ls:\ v\ge0 \text{ in }\R^N,\ v>0\text{ in }\Cw,\ \pl v\ge0\ \text{in }\Cw\big\},
\]
and $\beta^+$ is defined by \eqref{def:beta}. Functions in $\mathcal S_\beta(\omega)$ are only required to be nonnegative outside $\Cw$. Since $\mathcal S_\beta(\omega)\subset\mathcal S_\beta(\omega_0)$ for $\omega_0\subset\omega$, we have $\beta^+(\omega_0)\ge\beta^+(\omega)$ if $\omega_0\subset\omega$.

\section{The exponent \texorpdfstring{$\beta^+$}{beta+} and a nonexistence principle}\label{sec:fund}

\subsection{The range of the exponents}

{Throughout the paper, we fix} a nonempty smooth domain $\omega_0\Subset\omega$, we let $K_0:=\Sp\cap\overline{\mathcal C_{\omega_0}}=\overline{\omega_0}$, and we fix $r_0\in(0,\frac18]$ such that $B_{4r_0}(y)\subset\Cw$ for every $y\in K_0$. We let $\eta$ be the function of Lemma~\ref{lem:bump} with this $r_0$ (and, say, $k=3+\lceil\frac{sp}{p-1}\rceil$). {As} $\omega_0$ is smooth, there is $c_{\omega_0}>0$ with
\begin{equation}\label{eq:cone-measure}
\big|\mathcal C_{\omega_0}\cap B_{\rho}(y)\big|\ge c_{\omega_0}\rho^N\qquad\text{for } y\in K_0,\ 0<\rho\le 2r_0 .
\end{equation}

\begin{lemma}\label{lem:range}
\begin{itemize}
\item[(a)] If $N>sp$, then $|x|^{-\frac{N-sp}{p-1}}\in\mathcal S_{\frac{N-sp}{p-1}}(\omega)$; hence $\beta^+$ is well defined and $\beta^+\ge\frac{N-sp}{p-1}>0$. For every $N\ge2$, $\beta^+$ is well defined and positive by Corollary~\ref{cor:betaplus} below, whose proof {only uses Lemma~\ref{lem:range}(b), Lemmas~\ref{lem:aux}, \ref{lem:Lambda}, \ref{lem:compact}, Proposition~\ref{prop:eig}} and Lemmas~\ref{lem:ext} and~\ref{lem:zero}.
\item[(b)] There is $\delta_0=\delta_0(N,s,p,\omega)>0$ (depending also on the choices of $\omega_0,r_0,k$, which are fixed in terms of $\omega$) such that, if $\beta>0$ and $\mathcal S_\beta(\omega)\neq\emptyset$, then $\beta\le\frac{N}{p-1}-\delta_0$. \changed{Whenever its defining set is nonempty,} $\beta^+<\frac N{p-1}$.
\end{itemize}
\end{lemma}

\begin{proof}
(a) Let $\beta_0:=\frac{N-sp}{p-1}\in(0,\frac N{p-1})$ and $v_0:=|x|^{-\beta_0}$. By Theorem~\ref{thm:qp}, $\pl v_0=C(\beta_0)|x|^{-\beta_0(p-1)-sp}=0$ pointwise in $\R^N\setminus\{0\}$. Since $v_0$ is smooth with nonvanishing gradient in $\R^N\setminus\{0\}$, it is a viscosity solution there: if $\varphi$ touches $v_0$ from below at $x$ in $B_r(x)$, then $\varphi_r\le v_0$ with equality at $x$, so $\pl\varphi_r(x)\ge\pl v_0(x)=0$. Moreover $v_0\in H_{\beta_0}$, $v_0>0$, and $v_0\in\Ls$ by the integrability criterion above. So $v_0\in\mathcal S_{\beta_0}(\omega)$, and the set in \eqref{def:beta} is nonempty.

(b) Let $v\in\mathcal S_\beta(\omega)$ with $\beta>0$. We may assume $\beta\ge\frac{N}{2(p-1)}$ (otherwise {the conclusion is immediate if} $\delta_0\le\frac N{2(p-1)}$), and $\beta(p-1)<N$ because $v\in\Ls$. Let $m_0:=\min_{K_0}v>0$, attained at $\hat x\in K_0$, and $B:=B_{r_0}(\hat x)\subset\Cw$. We apply Remark~\ref{rem:touch} with $w=v$ and $F=0$: there are $c^*\in(0,m_0]$ and $x_1\in B$ such that $\phi:=c^*\eta(\cdot-\hat x)$ touches $v$ from below at $x_1$, $\phi\le v$ in $\R^N$, and for $r$ small
\begin{equation}\label{eq:range1}
0\le\pl\phi_r(x_1)=\pl\phi(x_1)+\int_{\R^N\setminus B_r(x_1)}\frac{\Psi_p(\phi(x_1)-v(y))-\Psi_p(\phi(x_1)-\phi(y))}{|x_1-y|^{N+sp}}dy .
\end{equation}
Since $v\ge\phi$, the last integrand is nonpositive everywhere. Let $r_1:=\min\{2^{-1/\beta},\frac12\}$ and $V:=\mathcal C_{\omega_0}\cap B_{r_1}$. For $y\in V,$ we have the following:
\begin{itemize}
\item $\phi(y)=0$, because $B\subset\{|y|\ge1-r_0\}$ and $r_1\le\frac12<1-r_0$;
\item $y/|y|\in\omega_0\subset K_0$, so $v(y)=|y|^{-\beta}v(y/|y|)\ge m_0|y|^{-\beta}$, and $|y|^{-\beta}\ge r_1^{-\beta}\ge2$;
\item $\phi(x_1)\le c^*\le m_0$, hence $\phi(x_1)-v(y)\le m_0-m_0|y|^{-\beta}\le-\frac12m_0|y|^{-\beta}$, and $\Psi_p(\phi(x_1)-v(y))-\Psi_p(\phi(x_1))\le-2^{1-p}m_0^{p-1}|y|^{-\beta(p-1)}$;
\item $|x_1-y|\le|x_1|+|y|\le\frac98+\frac12\le2$.
\end{itemize}
Also, $\pl\phi(x_1)=(c^*)^{p-1}\pl\eta(x_1-\hat x)\le m_0^{p-1}M_\eta,$ where $M_\eta$ is as in \eqref{meta}. Restricting the integral in \eqref{eq:range1} to $V$ (which does not meet $B_r(x_1)$) and dividing by $m_0^{p-1}$,
\[
M_\eta\ge 2^{1-p-N-sp}\int_{V}|y|^{-\beta(p-1)}dy=2^{1-p-N-sp}\,\mathcal H^{N-1}(\omega_0)\,\frac{r_1^{N-\beta(p-1)}}{N-\beta(p-1)} .
\]
Since $\beta\ge\frac N{2(p-1)}$, $r_1\ge r_*:=\min\{2^{-2(p-1)/N},\frac12\}$, and since $0<N-\beta(p-1)\le N$ and $r_1<1$, $r_1^{N-\beta(p-1)}\ge r_*^{N}$. Therefore
\[
N-\beta(p-1)\ge c_\omega:=\frac{2^{1-p-N-sp}\mathcal H^{N-1}(\omega_0)\,r_*^N}{M_\eta}>0,
\]
and (b) holds with $\delta_0:=\min\{\frac{c_\omega}{p-1},\frac N{2(p-1)}\}$.
\end{proof}

\subsection{Behaviour near the boundary of the cone}

Since $\partial\Cw\setminus\{0\}$ is of class $C^2$, there is $\rho\in(0,\frac1{24})$ such that every $\xi\in\partial\Cw$ with $\frac34\le|\xi|\le\frac43$ admits an interior ball $B_\rho(z_\xi)\subset\Cw$ and an exterior ball $B_\rho(c_\xi)\subset\R^N\setminus\overline{\Cw}$ tangent to $\partial\Cw$ at $\xi$, namely $z_\xi=\xi+\rho\nu_\xi$, $c_\xi=\xi-\rho\nu_\xi$ with $\nu_\xi$ the inner unit normal; and every $x\in\Cw\cap\Sp$ with $d(x)<\rho$ has a unique projection $\xi$ onto $\partial\Cw$ and satisfies $x=\xi+d(x)\nu_\xi$. Note that then $|\xi-x|=d(x)<\rho$, so $\frac{23}{24}<|\xi|<\frac{25}{24}$.

\begin{lemma}\label{lem:boundary}
Let $d(x)=\dist(x,\R^N\setminus\Cw)$.
\begin{itemize}
\item[(a)] \textup{(Hopf's)} Let $0<r_1\le r_2$, $\Lambda\ge0$, and let $v\in C(\R^N\setminus\{0\})$, $v\ge0$ in $\R^N$, $v>0$ in $\Cw\cap\{\frac{r_1}2\le|x|\le2r_2\}$, be a viscosity supersolution of $\pl v\ge-\Lambda\Psi_p(v)$ in $\Cw\cap\{\frac{r_1}2<|x|<2r_2\}$. Then there is $c_v>0$ with $v(x)\ge c_v\,d(x)^s$ for all $x\in\Cw\cap\{r_1\le|x|\le r_2\}$. In particular this holds on $\Cw\cap\Sp$ for every $v\in\mathcal S_\beta(\omega)$.
\item[(b)] \textup{(Upper barrier)} Let $u\in C(\R^N\setminus\{0\})\cap\Ls$ be nonnegative, $u=0$ in $\R^N\setminus\Cw$, and a viscosity subsolution of $\pl u\le F(x)$ in $\Cw\cap\{\frac12<|x|<2\}$, with $F$ continuous and $|F|\le M_F$. Then $u(x)\le C\,d(x)^s$ for $x\in\Cw\cap\Sp$, where $C$ depends only on $N,s,p,\omega$, $M_F$, $M_u:=\max_{\{\frac12\le|x|\le2\}}u$ and $T_u:=\int_{\R^N}u^{p-1}(1+|y|)^{-N-sp}dy$.
\end{itemize}
\end{lemma}

\begin{proof}
(a) {{We first treat the case} $r_1=r_2=1$.} Fix $x\in\Cw\cap\Sp$. Assume first $d(x)<\rho$; let $\xi$ be its projection and $z:=z_\xi$. Then $|x-z|=\rho-d(x)$, and $B_\rho(z)\subset B_{2\rho}(\xi)\subset\{\frac12<|y|<2\}$ because $|\xi|\in(\frac{23}{24},\frac{25}{24})$ and $2\rho<\frac1{12}$. So $v$ is a supersolution of $\pl v\ge-\Lambda\Psi_p(v)$ in $B_\rho(z)$ and $v>0$ there. The set $K_\rho:=\{y:\ d(y)\ge(1-\vartheta)\rho,\ \frac12\le|y|\le2\}$ is a compact subset of $\Cw$, so $m:=\min_{K_\rho}v>0$; since $B_{\vartheta\rho}(z)\subset K_\rho$, $\inf_{B_{\vartheta\rho}(z)}v\ge m$. By (K6), $v(x)\ge c\,(\rho-|x-z|)^s=c\,d(x)^s$, with $c$ depending only on $N,s,p,\rho,\Lambda,m$, hence independent of $x$. If $d(x)\ge\rho$, then $x\in K_\rho$ and $v(x)\ge m\ge m\,d(x)^s$, as $d(x)\le1$.

{Next, {we treat the general case}.} For $t\in[r_1,r_2]$ let $v_t(y):=v(ty)$. By scaling, $v_t$ is a supersolution of $\pl v_t\ge-\Lambda t^{sp}\Psi_p(v_t)$, hence of $\pl v_t\ge-\Lambda r_2^{sp}\Psi_p(v_t)$ (as $v_t\ge0$), in $\Cw\cap\{\frac12<|y|<2\}$. Applying the previous case to $v_t$ (with $m$ replaced by $\min\{v(y): d(y)\ge r_1(1-\vartheta)\rho,\ \frac{r_1}2\le|y|\le2r_2\}>0$, which bounds $\min_{K_\rho}v_t$ from below for all $t$) we get $v_t\ge c\,d^s$ on $\Cw\cap\Sp$ with $c$ independent of $t$. Since $d(ty)=t\,d(y)$, for $|x|=t$ we obtain $v(x)=v_t(x/t)\ge c\,t^{-s}d(x)^s\ge c\,r_2^{-s}d(x)^s$.

If $v\in\mathcal S_\beta(\omega)$, then $v$ is a supersolution of $\pl v\ge0$ in $\Cw$, and {the final assertion follows with} $\Lambda=0$.

(b)  Let $\chi\in C(\R^N)$, $0\le\chi\le1$, with $\chi=1$ in $\{\frac58\le|y|\le\frac32\}$ and $\chi=0$ outside $\{\frac12\le|y|\le2\}$, and $\tilde u:=\chi u$. Then $\tilde u$ is continuous, $0\le\tilde u\le u$, $\tilde u\le M_u$, and $\tilde u=u$ in $D':=\Cw\cap\{\frac58<|y|<\frac32\}$. By Lemma~\ref{lem:loc}, $\tilde u$ is a viscosity subsolution of $\pl\tilde u\le F+E$ in $D'$, where $E\ge0$ by Lemma~\ref{lem:loc}(i). The set $\Sigma$ where $\tilde u\ne u$ is contained in $\Cw\cap(\{|y|<\frac58\}\cup\{|y|>\frac32\})$ (outside $\Cw$ both vanish), which is at distance $\ge\frac1{24}$ from $D'':=\Cw\cap\{\frac23\le|y|\le\frac54\}$. By Lemma~\ref{lem:loc}(ii), $0\le E\le M_E$ on $D''$, with $M_E\le C(M_u^{p-1}+T_u)$.

Now  fix $x\in\Cw\cap\Sp$ with $d(x)<\rho$, let $\xi$ be its projection and $c:=c_\xi$, so that $|x-c|=\rho+d(x)$. Set
\[
W(y):=A\,w_*\Big(\frac{y-c}{\rho}\Big),\qquad D:=\Cw\cap B_{2\rho}(c),
\]
with $w_*$ from (K7) and $A>0$ to be chosen. Since $B_\rho(c)\cap\Cw=\emptyset$ and $B_{2\rho}(c)\subset B_{3\rho}(\xi)$, we have $D\subset\big(B_{2\rho}(c)\setminus\overline{B_\rho(c)}\big)\cap\{\frac34<|y|<\frac54\}\subset D''\cap D'$. By (K7) and scaling, $W$ is a viscosity supersolution of $\pl W\ge A^{p-1}\rho^{-sp}$ in $D$; also $0\le W\le AC_*$, so $W\in\Ls$. Choose
\[
A:=\max\big\{C_*M_u,\ (\rho^{sp}(M_F+M_E))^{1/(p-1)}\big\}.
\]
Then $\pl W\ge M_F+M_E$ in $D$, while $F+E\le M_F+M_E$ in $D$, so that $\tilde u$ is also a viscosity subsolution of $\pl\tilde u\le M_F+M_E$ in $D$. Outside $D$ we have $\tilde u\le W$: outside $\Cw$, $\tilde u=0\le W$; in $\Cw\setminus B_{2\rho}(c)$, $|y-c|\ge2\rho$, so $W(y)\ge A/C_*\ge M_u\ge\tilde u(y)$. By (K2) (with $c\equiv0$ and the constants $f_1=M_F+M_E\le f_2=A^{p-1}\rho^{-sp}$), $\tilde u\le W$ in $D$. Since $x\in D$ and $1<|x-c|/\rho<2$,
\[
u(x)=\tilde u(x)\le A\,C_*\Big(\frac{|x-c|-\rho}{\rho}\Big)^s=A\,C_*\,\rho^{-s}d(x)^s .
\]
If $d(x)\ge\rho$, then $u(x)\le M_u\le M_u\rho^{-s}d(x)^s$. The constant $AC_*\rho^{-s}$ depends only on the quantities in the statement.
\end{proof}

\subsection{A nonexistence principle}

The following proposition replaces the comparison principles for homogeneous functions of \cite{ASS,gab}. It only uses comparison in bounded domains. {The underlying idea, adapted from} \cite[Section 5]{AS}, is that a subsolution with homogeneity $-\beta$ can be placed below any positive supersolution far from the origin, {after its values near the vertex have been replaced by a sufficiently large bump, so that the modified function is a subsolution} of the homogeneous equation outside $B_2$.

\begin{proposition}\label{prop:nonex}
Let $\beta>0$ and let $U\in H_\beta\cap\Ls$, $U\ge0$, $U\not\equiv0$, $U=0$ in $\R^N\setminus\Cw$, be a viscosity subsolution of $\pl U\le0$ in $\Cw$. Let $v\in C(\R^N\setminus\{0\})\cap\Ls$, $v\ge0$ in $\R^N$, $v>0$ in $\Cw$, be a viscosity supersolution of $\pl v\ge0$ in $\Cw$. Then there is $\varepsilon>0$ such that
\[
v\ge\varepsilon\,U\qquad\text{in }\Cw\setminus B_2 .
\]
In particular, $\mathcal S_{\hat\beta}(\omega)=\emptyset$ for every $\hat\beta>\beta$.
\end{proposition}

\begin{proof}
%\emph{Step 0: bounds for $U$.}
Set $L:=\max_{\Sp}U>0$, so $U(y)\le L|y|^{-\beta}$. Applying Lemma~\ref{lem:boundary}(b) with $F=0$ we get $U\le C\,d^s$ on $\Cw\cap\Sp$, and by homogeneity, since $d(y/|y|)=d(y)/|y|$,
\begin{equation}\label{eq:Ud}
U(y)=|y|^{-\beta}U(y/|y|)\le C|y|^{-\beta-s}d(y)^s\qquad\text{for }y\in\Cw .
\end{equation}
Since $U\in\Ls$ and $U>0$ somewhere (hence on an open cone, by continuity and homogeneity), $\beta(p-1)<N$.

%\emph{Step 1: a subsolution outside $B_2$.}
 Let $A':=\mathcal C_{\omega_0}\cap\{\frac14\le|y|\le\frac12\}$, let $\chi\in C_c(\Cw\cap\{\frac18<|y|<\frac58\})$ with $0\le\chi\le1$ and $\chi=1$ on $A'$, and let $\vartheta\in C(\R)$, $0\le\vartheta\le1$, $\vartheta(t)=0$ for $t\le\frac34$, $\vartheta(t)=1$ for $t\ge1$. For a constant $C_b>0$ set
\[
z(y):=\vartheta(|y|)\,U(y)+C_b\chi(y).
\]
Then $z\in C(\R^N)\cap\Ls$ (the first term vanishes in $B_{3/4}$), $z\ge0$, $z=0$ outside $\Cw$, $z=0$ in $B_{1/8}$, and $z=U$ in $\R^N\setminus B_1$ (there $\vartheta=1$ and $\chi=0$). By Lemma~\ref{lem:loc} with $D=\Cw$ and $D'=\Cw\setminus\overline{B_1}$, $z$ is a viscosity subsolution of $\pl z\le E$ in $\Cw\setminus\overline{B_1}$, where, writing $a:=U(x)$,
\[
E(x)=\int_{B_1\cap\Cw}\frac{\Psi_p(a-z(y))-\Psi_p(a-U(y))}{|x-y|^{N+sp}}\,dy ,
\]
because outside $\Cw$ both $z$ and $U$ vanish. Let $|x|\ge2$; then $0\le a\le L2^{-\beta}\le L$, and $|x|-1\le|x-y|\le|x|+1$ for $y\in B_1$. We split $B_1\cap\Cw$ into three sets.
\begin{itemize}
\item On $P:=\{y\in B_1\cap\Cw:\ z(y)<U(y)\}$ the integrand is nonnegative; as $0\le z<U$ there, $|a-z(y)|,|a-U(y)|\le a+U(y)$, and by $|\Psi_p(A)-\Psi_p(B)|\le|A|^{p-1}+|B|^{p-1}$ the integrand is bounded by
\[
2(a+U(y))^{p-1}|x-y|^{-N-sp}\le 2^p\big(L^{p-1}+L^{p-1}|y|^{-\beta(p-1)}\big)(|x|-1)^{-N-sp}.
\]
{Consequently, the contribution of} $P$ is at most $K_1(|x|-1)^{-N-sp}$, $K_1:=2^pL^{p-1}\int_{B_1}(1+|y|^{-\beta(p-1)})dy<\infty$.
\item On $(B_1\cap\Cw)\setminus(P\cup A')$ we have $z\ge U$, and the integrand is $\le0$.
\item Assume $C_b\ge2L4^\beta$. For $y\in A'$ we have $\vartheta(|y|)=0$, so $z(y)=C_b$, and $U(y)\le L|y|^{-\beta}\le4^\beta L<C_b$; thus $A'\cap P=\emptyset$, and, $\Psi_p$ being increasing,
\[
\Psi_p(a-C_b)-\Psi_p(a-U(y))\le\Psi_p(L-C_b)-\Psi_p(-4^\beta L)=-(C_b-L)^{p-1}+(4^\beta L)^{p-1}=:-D(C_b),
\]
where $D(C_b)>0$ because $C_b-L\ge(2\cdot4^\beta-1)L>4^\beta L$. {Consequently, the contribution of} $A'$ is at most $-D(C_b)|A'|(|x|+1)^{-N-sp}$.
\end{itemize}
Using $\frac{|x|+1}{|x|-1}\le3$ for $|x|\ge2$,
\[
E(x)\le K_1(|x|-1)^{-N-sp}-D(C_b)|A'|(|x|+1)^{-N-sp}\le(|x|+1)^{-N-sp}\big[3^{N+sp}K_1-D(C_b)|A'|\big].
\]
Since $D(C_b)\to\infty$ as $C_b\to\infty$, {we may choose} $C_b\ge2L4^\beta$ with $D(C_b)|A'|\ge3^{N+sp}K_1$. Then $E\le0$ in $\Cw\setminus B_2$, and $z$ is a viscosity subsolution of $\pl z\le0$ in $\Cw\setminus\overline{B_2}$.

%\emph{Step 2: ordering in $\overline{B_2}$.}
 By Lemma~\ref{lem:boundary}(a) (with $\Lambda=0$, $r_1=\frac18$, $r_2=2$), $v\ge c\,d^s$ in $\Cw\cap\{\frac18\le|y|\le2\}$; and $v\ge c'>0$ on the compact set $\operatorname{supp}\chi\subset\Cw$. By \eqref{eq:Ud}, $U\le C''d^s$ in $\Cw\cap\{\frac18\le|y|\le2\}$. Hence, in $\Cw\cap\{\frac18\le|y|\le2\}$,
\[
z\le U+C_b\chi\le C''d^s+C_b\chi\le\Big(\frac{C''}{c}+\frac{C_b}{c'}\Big)v,
\]
while $z=0\le v$ in $B_{1/8}$ and outside $\Cw$. Thus $\varepsilon z\le v$ in $\overline{B_2}$, with $\varepsilon:=(\frac{C''}{c}+\frac{C_b}{c'})^{-1}$.

%\emph{Step 3: comparison in annuli.}
 Next, let $\delta>0$, $R>2$ with $\varepsilon LR^{-\beta}\le\delta$, and $D_R:=\Cw\cap(B_R\setminus\overline{B_2})$ be a bounded open set with $\overline{D_R}\subset\R^N\setminus\{0\}$. Thus the function $\varepsilon z-\delta$ is continuous, $(\varepsilon z-\delta)^+\in\Ls$, and it is a viscosity subsolution of $\pl(\varepsilon z-\delta)\le0$ in $D_R$ (as $\pl(\varepsilon z-\delta)=\varepsilon^{p-1}\pl z$). The function $v$ is a continuous supersolution of $\pl v\ge0$ in $D_R$. Moreover $\varepsilon z-\delta\le v$ in $\R^N\setminus D_R$. Indeed, it holds
\begin{itemize}
\item outside $\Cw$, because $\varepsilon z-\delta=-\delta<0\le v$;
\item in $\overline{B_2}$, by \changed{the ordering above};
\item in $\Cw\setminus B_R$, because $\varepsilon z=\varepsilon U\le\varepsilon L|y|^{-\beta}\le\varepsilon LR^{-\beta}\le\delta$.
\end{itemize}
Now by (K2) (with $c\equiv0$, $f_1=f_2=0$), $\varepsilon z-\delta\le v$ in $D_R$. Since this holds for every $R\ge R_\delta:=\max\{3,(\varepsilon L/\delta)^{1/\beta}\}$, we get $\varepsilon U-\delta\le v$ in $\Cw\setminus B_2$, and letting $\delta\to0$, it yields that $v\ge\varepsilon U$ in $\Cw\setminus B_2$.

%\emph{Step 4: different homogeneities.}
 Suppose $v\in\mathcal S_{\hat\beta}(\omega)$ with $\hat\beta>\beta$, and take $y_0\in\Cw\cap\Sp$ with $U(y_0)>0$. Then for $t\ge2$, \changed{the comparison above} gives $\varepsilon t^{-\beta}U(y_0)=\varepsilon U(ty_0)\le v(ty_0)=t^{-\hat\beta}v(y_0)$, that is, $\varepsilon U(y_0)\le t^{\beta-\hat\beta}v(y_0)\to0$ as $t\to\infty$, a contradiction.
\end{proof}

\begin{corollary}\label{cor:unique}
If $\Phi\in H_\beta\cap\Ls$, $\beta>0$, is a nonnegative nontrivial viscosity solution of \eqref{mp}, then $\beta=\beta^+$.
\end{corollary}

\begin{proof}
By Lemma~\ref{lem:smp} (with $c=0$), $\Phi>0$ in $\Cw$, so $\Phi\in\mathcal S_\beta(\omega)$ and $\beta\le\beta^+$ by definition of $\beta^+$. On the other hand $\Phi$ is a subsolution of $\pl\Phi\le0$ in $\Cw$ vanishing outside $\Cw$, so Proposition~\ref{prop:nonex} (with $U=\Phi$) gives $\mathcal S_{\hat\beta}(\omega)=\emptyset$ for all $\hat\beta>\beta$; hence $\beta^+\le\beta$.
\end{proof}

\section{Homogeneous eigenfunctions}\label{sec:eig}

Throughout this section $\beta\in(0,\frac N{p-1})$. For $u\in H_\beta$ we write $\hat u:=u|_{\Sp}$ and $\|u\|:=\max_{\Sp}|u|$. A pair $(\mu,U)$ with $\mu\in\R$, $U\in H_\beta$ continuous, $U\ge0$, $U\not\equiv0$, $U=0$ outside $\Cw$, which is a viscosity solution of
\begin{equation}\label{eq:eig}
\pl U=\mu|x|^{-sp}U^{p-1}\qquad\text{in }\Cw,
\end{equation}
is called a \emph{homogeneous eigenpair} (with exponent $\beta$). By Lemma~\ref{lem:smp} (with $c(x)=\mu^-|x|^{-sp}$), $U>0$ in $\Cw$. Note that $U\in\Ls$ by the integrability criterion of Section~\ref{sec:prel}.

\subsection{An auxiliary problem}
Let $M=M_\beta:=|C(\beta)|+1$, so that $C(\beta)+M\ge1$, where $C(\beta)$ is the constant of Theorem~\ref{thm:qp}. For $\lambda\ge0$ and a continuous $\psi\in H_\beta$ with $\psi\ge0$ and $\psi=0$ outside $\Cw$, consider
\begin{equation}\label{pc}
\pl V+M|x|^{-sp}\Psi_p(V)=\lambda R_\psi,\quad R_\psi:=|x|^{-sp}\big(\psi^{p-1}+|x|^{-\beta(p-1)}\big)\ \text{ in }\Cw,\qquad V=0\ \text{ in }\R^N\setminus\Cw .
\end{equation}
If $V\in H_\beta$, all the terms in \eqref{pc} are homogeneous of degree $-\beta(p-1)-sp$. The term $|x|^{-\beta(p-1)-sp}$ in $R_\psi$ makes the right-hand side positive in $\Cw$ when $\lambda>0$; {this allows the continuation argument below to proceed} {without prior information about the eigenvalue problem}.

\begin{lemma}\label{lem:aux}
Problem \eqref{pc} has a {minimal} nonnegative solution $V=V_{\lambda}[\psi]$: $V\in H_\beta$ is continuous, nonnegative, vanishes outside $\Cw$, is a viscosity solution of \eqref{pc}, and $V\le W$ for every nonnegative $W\in C(\R^N\setminus\{0\})\cap\Ls$ which is a viscosity supersolution of the equation in \eqref{pc}. Moreover:
\begin{itemize}
\item[(i)] $\|V\|^{p-1}\le\lambda\big(\|\psi\|^{p-1}+1\big)$; in particular $V=0$ if $\lambda=0$;
\item[(ii)] $V\le C\,d^s$ on $\Cw\cap\Sp$, with $C$ bounded when $\lambda$ and $\|\psi\|$ are bounded;
\item[(iii)] if $R_2\le\theta^{p-1}R_1$ for some $\theta>0$, where $R_i=\lambda_iR_{\psi_i}$, then $V_{\lambda_2}[\psi_2]\le\theta V_{\lambda_1}[\psi_1]$.
\end{itemize}
\end{lemma}

\begin{proof}
%\emph{Approximating problems.}
 Assume that $\mathcal D$ be the family of bounded open sets $D$ of class $C^{1,1}$ with $\overline D\subset\Cw$. For $D\in\mathcal D$ the coefficient $c(x)=M|x|^{-sp}$ and the datum $f=\lambda R_\psi$ are continuous and bounded in $\overline D$, since $\overline D$ is a compact set not containing $0$. By (K3) there is a unique weak solution $V_D$ of the equation in $D$ with $V_D=0$ outside $D$; it is continuous in $\R^N$ and, by (K1), a viscosity solution in $D$. {We note three consequences of} (K2):
\begin{itemize}
\item $V_D\ge0$: the function $0$ is a subsolution with $f_1=0\le f_2=\lambda R_\psi$, and $0\le V_D$ outside $D$;
\item if $D\subset D'$ in $\mathcal D$, then $V_D\le V_{D'}$: compare in $D$, using $V_D=0\le V_{D'}$ outside $D$;
\item if $W$ is as in the statement, then $V_D\le W$: compare in $D$, using $V_D=0\le W$ outside $D$.
\end{itemize}

%\emph{Upper bound.}
Let $w_K:=K|x|^{-\beta}$, \mathred{$K\ge0$, with $w_0\equiv0$}. Then using Theorem~\ref{thm:qp}, we obtain
\[
\pl w_K+M|x|^{-sp}w_K^{p-1}=K^{p-1}(C(\beta)+M)|x|^{-\beta(p-1)-sp}\ge K^{p-1}|x|^{-\beta(p-1)-sp}
\]
pointwise, hence in the viscosity sense (as in the proof of Lemma~\ref{lem:range}(a)). Since $\psi\le\|\psi\||x|^{-\beta}$, we have $\lambda R_\psi\le\lambda(\|\psi\|^{p-1}+1)|x|^{-\beta(p-1)-sp}$. Hence, for $K^{p-1}:=\lambda(\|\psi\|^{p-1}+1)$, $w_K$ is a nonnegative continuous supersolution in $\Cw$, $w_K\in\Ls$, and $V_D\le w_K$ for all $D\in\mathcal D$.

%\emph{Definition of $V$ and equation.}
Next, set $V:=\sup_{D\in\mathcal D}V_D$. Since $\mathcal D$ is directed by inclusion (two sets of $\mathcal D$ are contained in a third one) and $D\mapsto V_D$ is nondecreasing, $V=\lim_kV_{D_k}$ for any nondecreasing sequence $D_k\in\mathcal D$ exhausting $\Cw$ (every $D\in\mathcal D$ is contained in some $D_k$, by compactness of $\overline D$). Clearly $V\le W$ for every $W$ as in the statement; in particular $V\le w_K$, which gives (i). Fix a compact set $K'\subset\Cw$ and $k_0$ with $K'_{\delta}\subset D_{k_0}$ for some $\delta>0$. For $k\ge k_0$, the functions $V_{D_k}$ are solutions in $D_{k_0}$ with right-hand sides, \changed{local} $L^\infty$ norms and tails bounded independently of $k$ (by $w_K$); by (K4) (after the modification of Remark~\ref{rem:K1}) they are equi-H\"older continuous on $K'$. {A pointwise-convergent, nondecreasing sequence of equicontinuous functions has a continuous limit and converges uniformly on compact sets} (Dini's theorem). So $V$ is continuous in $\Cw$ and $V_{D_k}\to V$ locally uniformly in $\Cw$; also in $\Ls$, by dominated convergence (all functions are bounded by $w_K$ and vanish outside $\Cw$). By (K5), applied in every bounded open $D''\Subset\Cw$ (which is contained in $D_k$ for $k$ large), $V$ is a viscosity solution of the equation in \eqref{pc}.

%\emph{Homogeneity.}
 Let $t>0$ and $D\in\mathcal D$, and set $V^t(x):=t^\beta V_D(tx)$. By scaling, $\pl V^t(x)=t^{\beta(p-1)+sp}(\pl V_D)(tx)$, and since $R_\psi(tx)=t^{-\beta(p-1)-sp}R_\psi(x)$ (because $\psi\in H_\beta$) and $|tx|^{-sp}\Psi_p(V_D(tx))=t^{-\beta(p-1)-sp}|x|^{-sp}\Psi_p(V^t(x))$, the function $V^t$ is a weak solution of the equation in $D/t\in\mathcal D$, vanishing outside $D/t$. {By the uniqueness statement in (K3)}, $V^t=V_{D/t}$. Since $D\mapsto D/t$ is a bijection of $\mathcal D$, $t^\beta V(tx)=\sup_Dt^\beta V_D(tx)=\sup_DV_{D/t}(x)=V(x)$.

%\emph{Boundary behaviour.}
 Observe that each $V_{D_k}$ is a continuous viscosity subsolution of $\pl V_{D_k}\le\lambda R_\psi$ in $\Cw$. In $D_k$ this holds because $M|x|^{-sp}\Psi_p(V_{D_k})\ge0$. At a point $x\in\Cw\setminus D_k$ we have $V_{D_k}(x)=0=\min V_{D_k}$; if \changed{an admissible} $\varphi$ touches $V_{D_k}$ from above at $x$ in $B_r(x)$, then $\varphi_r\ge0=\varphi_r(x)$ in $\R^N$, so every integrand in $\pl\varphi_r(x)$ is $\le0$ and $\pl\varphi_r(x)\le0\le\lambda R_\psi(x)$. Now apply Lemma~\ref{lem:boundary}(b) with $F=\lambda R_\psi$, which is continuous and bounded by $\lambda(\|\psi\|^{p-1}+1)2^{\beta(p-1)+sp}$ in $\{\frac12<|x|<2\}$, $M_u\le\max_{\{1/2\le|x|\le2\}}w_K$ and $T_u\le\int w_K^{p-1}(1+|y|)^{-N-sp}dy$: we get $V_{D_k}\le C\,d^s$ on $\Cw\cap\Sp$, with $C$ independent of $k$ and bounded when $\lambda,\|\psi\|$ are bounded. Letting $k\to\infty$ gives (ii). By (ii) and homogeneity, $V(y)\le C|y|^{-\beta-s}d(y)^s\to0$ as $y\to\xi\in\partial\Cw\setminus\{0\}$, so $V$ (which vanishes outside $\Cw$) is continuous in $\R^N\setminus\{0\}$.

(iii) \mathred{Set $V_i:=V_{\lambda_i}[\psi_i]$, $i=1,2$.} By the $(p-1)$-homogeneity of the left-hand side of \eqref{pc}, $\theta V_{\lambda_1}[\psi_1]$ satisfies $\pl(\theta V_1)+M|x|^{-sp}\Psi_p(\theta V_1)=\theta^{p-1}R_1\ge R_2$, so it is a nonnegative continuous supersolution of the problem with datum $R_2$, and minimality gives $V_{\lambda_2}[\psi_2]\le\theta V_{\lambda_1}[\psi_1]$.
\end{proof}

\subsection{An upper bound for eigenvalues}

\begin{lemma}\label{lem:Lambda}
Let $\bar\beta\in(0,\frac N{p-1})$. There is $\Lambda_0=\Lambda_0(N,s,p,\omega,\bar\beta)$ such that, if $\beta\in(0,\bar\beta]$, $\mu\ge0$, and $\psi\in H_\beta\cap\Ls$, $\psi\ge0$, $\psi>0$ in $\Cw$, is a viscosity supersolution of $\pl\psi\ge\mu|x|^{-sp}\psi^{p-1}$ in $\Cw$, then $\mu\le\Lambda_0$.
\end{lemma}

\begin{proof}
Let $\omega_0$, $K_0$, $r_0$ and $\eta$ be as fixed at the beginning of Section~\ref{sec:fund}. Let $m_0:=\min_{K_0}\psi>0$, attained at $\hat x\in K_0$, and apply Remark~\ref{rem:touch} with $w=\psi$, $B=B_{r_0}(\hat x)$: there are $c^*\in(0,m_0]$ and $x_1\in B$ such that, by \eqref{eq:touch},
\begin{equation}\label{eq:Lambda1}
\mu|x_1|^{-sp}\psi(x_1)^{p-1}\le (c^*)^{p-1}M_\eta\le m_0^{p-1}M_\eta .
\end{equation}
{We next bound} $\psi(x_1)$ from below by the weak Harnack inequality. Since $\mu\ge0$, $\psi$ is a supersolution of $\pl\psi\ge0$ in $B_{4r_0}(\hat x)\subset\Cw\cap\{|y|\ge\frac12\}$. By Remark~\ref{rem:K1} (with $D=B_{4r_0}(\hat x)$), the function $\tilde\psi:=\min\{\psi,M'\zeta\}$ is a bounded, continuous, nonnegative weak supersolution of $\pl\tilde\psi\ge0$ in $B_{4r_0}(\hat x)$ which coincides with $\psi$ there. By (K4),
\[
\psi(x_1)\ge\inf_{B_{2r_0}(\hat x)}\psi\ge c_H\Big(\frac{1}{|B_{2r_0}|}\int_{B_{2r_0}(\hat x)\cap\mathcal C_{\omega_0}}\psi^t\Big)^{1/t}.
\]
For $y\in\mathcal C_{\omega_0}\cap B_{2r_0}(\hat x)$ we have $y/|y|\in\omega_0\subset K_0$ and $|y|\le1+2r_0\le\frac32$, so $\psi(y)=|y|^{-\beta}\psi(y/|y|)\ge(3/2)^{-\bar\beta}m_0$. By \eqref{eq:cone-measure}, $|\mathcal C_{\omega_0}\cap B_{2r_0}(\hat x)|\ge c_{\omega_0}(2r_0)^N$. Hence
\[
\psi(x_1)\ge c_H\Big(\frac{c_{\omega_0}}{|B_1|}\Big)^{1/t}(3/2)^{-\bar\beta}m_0=:c_1m_0 .
\]
Inserting this in \eqref{eq:Lambda1} and using $|x_1|\le1+r_0\le\frac54$, we get $\mu(5/4)^{-sp}c_1^{p-1}m_0^{p-1}\le m_0^{p-1}M_\eta$, that is, $\mu\le(5/4)^{sp}c_1^{1-p}M_\eta=:\Lambda_0$.
\end{proof}

\subsection{Existence of homogeneous eigenpairs}

We use the following global continuation result \cite{Rabinowitz1971}, in its version for maps leaving a cone invariant, which relies on the fixed point index in cones (the index of $0$ for $T(0,\cdot)\equiv0$ is $1$); see \cite{Dancer1973}, \cite[Section~17]{Amann1976} and \cite{ASS}.

\begin{proposition}\label{prop:rab}
Let $X$ be a real Banach space, $\mathcal K\subset X$ a closed convex cone, and $T:[0,\infty)\times\mathcal K\to\mathcal K$ a continuous map which sends bounded sets into relatively compact sets, with $T(0,u)=0$ for all $u\in\mathcal K$. Then there is an unbounded connected set $\mathcal S\subset[0,\infty)\times\mathcal K$ with $(0,0)\in\mathcal S$ such that $T(\lambda,u)=u$ for all $(\lambda,u)\in\mathcal S$.
\end{proposition}

\begin{lemma}\label{lem:compact}
Let $0<b_1\le b_2<\frac N{p-1}$ and $\bar\mu>0$. Let $(\mu_j,U_j)$ be homogeneous eigenpairs with exponents $\beta_j\in[b_1,b_2]$, $|\mu_j|\le\bar\mu$ and $\|U_j\|=1$. Then, up to a subsequence, $\beta_j\to\beta$, $\mu_j\to\mu$, $\hat U_j\to\hat U$ uniformly on $\Sp$, and $(\mu,U)$, $U(x):=|x|^{-\beta}\hat U(x/|x|)$, is a homogeneous eigenpair with exponent $\beta$ and $\|U\|=1$. The same holds if the $U_j$ \textup{(}continuous, nonnegative, $(-\beta_j)$-homogeneous and vanishing outside $\Cw$\textup{)} solve $\pl U_j=\mu_j|x|^{-sp}U_j^{p-1}+\kappa_j|x|^{-\beta_j(p-1)-sp}$ in $\Cw$ with $\kappa_j\ge0$, $\kappa_j\to0$.
\end{lemma}

\begin{proof}
We may assume $\kappa_j\le1$. Since $\|U_j\|=1$, $0\le U_j(x)\le|x|^{-\beta_j}\le|x|^{-b_1}+|x|^{-b_2}$. Let $A_0:=\Cw\cap\{\frac12<|x|<2\}$.
%\emph{Uniform bounds.}
 In $A_0$ the right-hand sides satisfy $\changed{\bigl|}\mu_j|x|^{-sp}U_j^{p-1}+\kappa_j|x|^{-\beta_j(p-1)-sp}\changed{\bigr|}\le(\bar\mu+1)2^{b_2(p-1)+sp}$, the functions $U_j$ are bounded by $2^{b_2}$, and their tails are bounded by $\int(|y|^{-b_1}+|y|^{-b_2})^{p-1}(1+|y|)^{-N-sp}dy<\infty$, since $b_2(p-1)<N$.

%\emph{Interior equicontinuity.}
{Applying (K4)} (via (K1) and Remark~\ref{rem:K1}),  {we find that the} $U_j$ are equi-H\"older continuous on every compact subset of $A_0$.

%\emph{Equicontinuity up to $\partial\omega$.}
 By Lemma~\ref{lem:boundary}(b), applied with $F=\mu_j|x|^{-sp}U_j^{p-1}+\kappa_j|x|^{-\beta_j(p-1)-sp}$ and the bounds above, $U_j\le C\,d^s$ on $\Cw\cap\Sp$ with $C$ independent of $j$. Given $\epsilon>0$, choose $\tau>0$ with $C\tau^s<\epsilon$; on $\{x\in\Sp: d(x)\le\tau\}$ all $\hat U_j$ are $<\epsilon$, and they vanish on $\Sp\setminus\omega$, while on the compact set $\{x\in\Sp:\ d(x)\ge\tau\}$ they are equicontinuous. It follows that $\{\hat U_j\}$ is equicontinuous on $\Sp$.

%\emph{Passage to the limit.}
{Applying Arzel}\`a--Ascoli and using compactness of $[b_1,b_2]\times[-\bar\mu,\bar\mu]$, {we obtain a subsequence for which} $\beta_j\to\beta$, $\mu_j\to\mu$ and $\hat U_j\to\hat U$ uniformly on $\Sp$. Then $U_j(x)=|x|^{-\beta_j}\hat U_j(x/|x|)\to U(x)$ locally uniformly in $\R^N\setminus\{0\}$, and in $\Ls$ by dominated convergence. The nonlinearities $F_j(x,\xi)=\mu_j|x|^{-sp}\Psi_p(\xi)+\kappa_j|x|^{-\beta_j(p-1)-sp}$ converge locally uniformly in $\Cw\times\R$ to $\mu|x|^{-sp}\Psi_p(\xi)$. By (K5), applied in every bounded open $D\Subset\Cw$, $U$ solves \eqref{eq:eig}. Finally $U\in H_\beta$ is continuous, $U\ge0$, $U=0$ outside $\Cw$, and $\|U\|=\lim\|U_j\|=1$, so $U\not\equiv0$.
\end{proof}

\begin{proposition}\label{prop:eig}
For every $\beta\in(0,\frac N{p-1})$ there is a homogeneous eigenpair $(\mu_\beta,U_\beta)$ with exponent $\beta$, $\|U_\beta\|=1$ and $-M_\beta\le\mu_\beta\le\Lambda_0(N,s,p,\omega,\beta)$.
\end{proposition}

\begin{proof}
Consider $X=C(\Sp)$ with the maximum norm, $\mathcal K=\{u\in X:\ u\ge0,\ u=0\text{ on }\Sp\setminus\omega\}$, a closed convex cone, and for $u\in\mathcal K$ let $\psi_u(x):=|x|^{-\beta}u(x/|x|)$; then $\psi_u\in H_\beta$ is continuous, nonnegative and vanishes outside $\Cw$, with $\|\psi_u\|=\|u\|_X$. Define
\[
T(\lambda,u):=\widehat{V_\lambda[\psi_u]},
\]
the restriction to $\Sp$ of the minimal solution of Lemma~\ref{lem:aux}. Since $V_\lambda[\psi_u]$ is continuous, nonnegative and vanishes outside $\Cw$, $T(\lambda,u)\in\mathcal K$; and $T(0,u)=0$ by Lemma~\ref{lem:aux}(i).

%\emph{Compactness.}
 Let $\lambda\le\bar\lambda$ and $\|u\|\le\bar m$. By Lemma~\ref{lem:aux}(i), $\|T(\lambda,u)\|^{p-1}\le\bar\lambda(\bar m^{p-1}+1)$. The functions $V_\lambda[\psi_u]$ solve an equation of the form \eqref{eq:gen} in $\Cw\cap\{\frac12<|x|<2\}$ with uniformly bounded coefficients, data and tails, so they are equi-H\"older continuous on compact subsets by (K4); and $V_\lambda[\psi_u]\le C\,d^s$ on $\Cw\cap\Sp$ by Lemma~\ref{lem:aux}(ii), with $C$ uniform. As in the proof of Lemma~\ref{lem:compact}, the family $\{T(\lambda,u)\}$ is equicontinuous on $\Sp$, and by Arzel\`a--Ascoli it is relatively compact in $X$.

%\emph{Continuity.}
{Now suppose that} $(\lambda_j,u_j)\to(\lambda,u)$ in $[0,\infty)\times\mathcal K$. If $\lambda=0$, then $\|T(\lambda_j,u_j)\|^{p-1}\le\lambda_j(\|u_j\|^{p-1}+1)\to0=\|T(0,u)\|$. Let $\lambda>0$; then $\lambda_j>0$ for large $j$. Let $\varepsilon_j:=\max_{\Sp}|u_j^{p-1}-u^{p-1}|$, which tends to $0$ because $u_j\to u$ uniformly and $\tau\mapsto\tau^{p-1}$ is uniformly continuous on bounded subsets of $[0,\infty)$. Then $\psi_{u_j}^{p-1}\le\psi_u^{p-1}+\varepsilon_j|x|^{-\beta(p-1)}$, and since $R_{\psi_u}\ge|x|^{-\beta(p-1)-sp}$,
{\color{\revisioncolor}
\[
\begin{aligned}
\lambda_jR_{\psi_{u_j}}
&\le\lambda_j R_{\psi_u}+\lambda_j\varepsilon_j|x|^{-\beta(p-1)-sp}\\
&\le\lambda_j(1+\varepsilon_j)R_{\psi_u}
\le\theta_j^{p-1}\lambda R_{\psi_u},\\
\theta_j^{p-1}&:=\max\Big\{\frac{\lambda_j(1+\varepsilon_j)}{\lambda},1\Big\}.
\end{aligned}
\]
}
In the same way, $\lambda R_{\psi_u}\le\tilde\theta_j^{p-1}\lambda_jR_{\psi_{u_j}}$ with $\tilde\theta_j^{p-1}:=\max\{\frac{\lambda(1+\varepsilon_j)}{\lambda_j},1\}$. Clearly $\theta_j,\tilde\theta_j\to1$. By Lemma~\ref{lem:aux}(iii),
\[
\tilde\theta_j^{-1}\,T(\lambda,u)\le T(\lambda_j,u_j)\le\theta_j\,T(\lambda,u)\qquad\text{on }\Sp,
\]
so $\|T(\lambda_j,u_j)-T(\lambda,u)\|\le\max\{\theta_j-1,1-\tilde\theta_j^{-1}\}\|T(\lambda,u)\|\to0$.

%\emph{Continuation.}
 Let $\mathcal S$ be given by Proposition~\ref{prop:rab}. If $(\lambda,u)\in\mathcal S$, then $\psi:=\psi_u$ and $V_\lambda[\psi_u]$ are $(-\beta)$-homogeneous with the same trace on $\Sp$, so $\psi=V_\lambda[\psi]$, that is,
\begin{equation}\label{eq:fixed}
\pl\psi=(\lambda-M)|x|^{-sp}\psi^{p-1}+\lambda|x|^{-\beta(p-1)-sp}\quad\text{in }\Cw,\qquad\psi=0\ \text{ outside }\Cw .
\end{equation}
Assume that $\lambda>M$. Then $\psi\not\equiv0$ (the function $0$ does not solve \eqref{eq:fixed} when $\lambda>0$), $\pl\psi\ge0$ in $\Cw$, so $\psi>0$ in $\Cw$ by Lemma~\ref{lem:smp}, and $\psi$ is a supersolution of $\pl\psi\ge(\lambda-M)|x|^{-sp}\psi^{p-1}$ with $\lambda-M>0$. By Lemma~\ref{lem:Lambda} (with $\bar\beta=\beta$), $\lambda-M\le\Lambda_0$. Hence $0\le\lambda\le M+\Lambda_0$ on $\mathcal S$, and since $\mathcal S$ is unbounded there are $(\lambda_j,u_j)\in\mathcal S$ with $\|u_j\|\to\infty$; {after passing to a subsequence}, $\lambda_j\to\lambda_\infty\in[0,M+\Lambda_0]$. The functions $U_j:=\psi_{u_j}/\|u_j\|$ satisfy $\|U_j\|=1$ and, by the $(p-1)$-homogeneity of $\pl$,
\[
\pl U_j=(\lambda_j-M)|x|^{-sp}U_j^{p-1}+\kappa_j|x|^{-\beta(p-1)-sp},\qquad \kappa_j:=\lambda_j\|u_j\|^{1-p}\to0 .
\]
By Lemma~\ref{lem:compact} (with $\beta_j=\beta$, $\mu_j=\lambda_j-M$, $\bar\mu=M+\Lambda_0$), a subsequence converges to an eigenpair $(\mu_\beta,U_\beta)$ with $\mu_\beta=\lambda_\infty-M\in[-M,\Lambda_0]$ and $\|U_\beta\|=1$.
\end{proof}

\begin{corollary}\label{cor:sign}
Let $(\mu,U)$ be a homogeneous eigenpair with exponent $\beta$. If $\beta<\beta^+$ then $\mu>0$; if $\beta>\beta^+$ then $\mu<0$.
\end{corollary}

\begin{proof}
Let $\beta<\beta^+$ and assume $\mu\le0$. Then $\pl U=\mu|x|^{-sp}U^{p-1}\le0$ in $\Cw$, so $U$ satisfies the hypotheses of Proposition~\ref{prop:nonex}. By definition of $\beta^+$ as a supremum there is $\hat\beta\in(\beta,\beta^+]$ with $\mathcal S_{\hat\beta}(\omega)\ne\emptyset$, contradicting Proposition~\ref{prop:nonex}. Let $\beta>\beta^+$ and assume $\mu\ge0$. Then $\pl U\ge0$ in $\Cw$, $U>0$ in $\Cw$ by Lemma~\ref{lem:smp}, and $U\in H_\beta\cap\Ls$, so $U\in\mathcal S_\beta(\omega)$ with $\beta>\beta^+$, contradicting the definition of $\beta^+$.
\end{proof}

\begin{proposition}\label{prop:limit}
There is a homogeneous eigenpair $(\underline\mu,\underline U)$ with exponent $\beta^+$, $\|\underline U\|=1$ and $\underline\mu\ge0$.
\end{proposition}

\begin{proof}
Recall $0<\beta^+<\frac N{p-1}$ (Lemma~\ref{lem:range} and Corollary~\ref{cor:betaplus}). Take $\beta_j\uparrow\beta^+$ with $\beta_j\ge\beta^+/2$, and let $(\mu_j,U_j)$ be the eigenpairs given by Proposition~\ref{prop:eig}, with $\|U_j\|=1$. By Corollary~\ref{cor:sign}, $\mu_j>0$; then $U_j$ is a positive supersolution of $\pl U_j\ge\mu_j|x|^{-sp}U_j^{p-1}$ with $\beta_j\le\beta^+$, and Lemma~\ref{lem:Lambda} (with $\bar\beta=\beta^+$) gives $\mu_j\le\Lambda_0(N,s,p,\omega,\beta^+)$. Lemma~\ref{lem:compact} (with $b_1=\beta^+/2$, $b_2=\beta^+$) gives a limit eigenpair with exponent $\beta^+$, $\|\underline U\|=1$ and $\underline\mu=\lim\mu_j\ge0$.
\end{proof}

\section{Existence of the fundamental solution $\Phi^+$}\label{sec:exist}

\begin{proof}[Proof of Theorem~\ref{mt}\textup{(a)--(c)}]
(a) is Lemma~\ref{lem:range} and Corollary~\ref{cor:betaplus}; (b) follows from Proposition~\ref{prop:eig} and Corollary~\ref{cor:sign}; (c) is Corollary~\ref{cor:unique}.
\end{proof}

{To establish Theorem~}\ref{mt} (d), we show that the eigenvalue $\underline\mu$ of Proposition~\ref{prop:limit} is $0$; then $\Phi^+:=\underline U$ is a nonnegative $(-\beta^+)$-homogeneous solution of $\pl\Phi^+=0$ in $\Cw$, vanishing outside $\Cw$ and positive in $\Cw$, that is, the required fundamental solution. In this section we assume
\begin{equation}\label{eq:range-p}
p-1>sp ,\qquad\text{that is,}\qquad p>\tfrac1{1-s},
\end{equation}
and, arguing by contradiction, that $\mu:=\underline\mu>0$. We write $U:=\underline U$, so that $\|U\|=1$, $U\in H_{\beta^+}$, $0\le U(y)\le|y|^{-\beta^+}$, $U>0$ in $\Cw$ and $\pl U=\mu|x|^{-sp}U^{p-1}$ in $\Cw$. Since $\mu>0$, $U$ is a supersolution of $\pl U\ge0$, so Lemma~\ref{lem:boundary}(a) applies with $\Lambda=0$; and $U$ is a subsolution of $\pl U\le F$ with $F=\mu|x|^{-sp}U^{p-1}$, which is bounded by $\mu2^{sp+\beta^+(p-1)}$ in $\{\frac12<|x|<2\}$, so Lemma~\ref{lem:boundary}(b) applies. Using also homogeneity ($U(y)=|y|^{-\beta^+}U(y/|y|)$ and $d(y)=|y|\,d(y/|y|)$), there are $0<c_0\le 1\le C_0$ with
\begin{equation}\label{eq:twosided}
c_0\,d(y)^s\le U(y)\le C_0\,d(y)^s\qquad\text{for } y\in\Cw,\ \tfrac12\le|y|\le2 .
\end{equation}
If $p\ge2$ we fix $\alpha\in(0,1)$ such that
\begin{equation}\label{eq:alpha}
\alpha<\frac{sp}{p-1}\qquad\text{and}\qquad 1+\alpha(p-2)>sp .
\end{equation}
{Such a choice is possible. Indeed,  first note that} $\frac{sp}{p-1}<1$ by \eqref{eq:range-p}. If $p=2$, {the second condition becomes} $1>2s$, which is \eqref{eq:range-p}, and any $\alpha\in(0,\frac{sp}{p-1})$ works. If $p>2$, {the second condition becomes} $\alpha>\frac{sp-1}{p-2}$, and
\[
\frac{sp-1}{p-2}<\frac{sp}{p-1}\iff(sp-1)(p-1)<sp(p-2)\iff -(p-1)<-sp\iff sp<p-1 ,
\]
so we may take $\alpha\in\big(\max\{0,\frac{sp-1}{p-2}\},\frac{sp}{p-1}\big)$.

\begin{lemma}\label{lem:gradU}
Let $p\ge2$. There is $C_1$ such that, for every $x\in\Cw\cap\Sp$ and $y\in B_{d(x)/64}(x)$,
\[
|U(y)-U(x)|\le C_1\,d(x)^{s-\alpha}|y-x|^\alpha .
\]
\end{lemma}

\begin{proof}
Let $d:=d(x)\le1$, $r:=d/2\le\frac12$, and $U_r(\xi):=r^{-s}U(x+r\xi)$, $\xi\in\R^N$. Note that $B_{2r}(x)=B_d(x)\subset\Cw$ and $B_r(x)\subset\{\frac12\le|y|\le\frac32\}$.

%\emph{Equation.}
 {Scaling gives}
\[
\begin{aligned}
\pl U_r(\xi)&=r^{-s(p-1)}r^{sp}(\pl U)(x+r\xi)=r^{s}\mu|x+r\xi|^{-sp}U(x+r\xi)^{p-1}\\
&=r^{sp}\mu|x+r\xi|^{-sp}U_r(\xi)^{p-1}=:f_r(\xi)
\end{aligned}
\]
in $B_2$ (in the viscosity sense).

%\emph{Bounds.}
{For} $\xi\in B_1$, $y=x+r\xi$ satisfies $d(y)\le d+r=3r$ and $\frac12\le|y|\le\frac32$, so by \eqref{eq:twosided} $U_r(\xi)\le C_0r^{-s}(3r)^s=3^sC_0$. Hence $\|U_r\|_{L^\infty(B_1)}\le3C_0$ and $\|f_r\|_{L^\infty(B_1)}\le\mu2^{sp}(3C_0)^{p-1}$. Moreover $U_r$ is bounded on $B_{3/2}$ independently of $x$ (as needed for the cut-off below): if $r<\frac13$ this follows from the pointwise bound $U_r(\xi)\le3^sC_0|\xi|^s$ for $1\le|\xi|\le\frac1{2r}$ proved in the tail estimate below, and if $r\ge\frac13$ then $|x+r\xi|\ge1-\frac32r\ge\frac14$ for $\xi\in B_{3/2}$, so $U_r\le r^{-s}4^{\beta^+}\le3\cdot4^{\beta^+}$. \mathred{The equation also bounds $f_r$ uniformly on $B_{5/4}$, since $|x+r\xi|\ge3/8$ there.}

 If $1\le|\xi|\le\frac1{2r}$, then $|y-x|\le\frac12$, so $\frac12\le|y|\le\frac32$ and $d(y)\le d+r|\xi|\le3r|\xi|$, whence $U_r(\xi)\le3^sC_0|\xi|^s$. If $|\xi|>\frac1{2r}$, we change variables $y=x+r\xi$. Therefore
\[
\int_{|\xi|\ge1}\frac{U_r(\xi)^{p-1}}{|\xi|^{N+sp}}d\xi\le(3C_0)^{p-1}\int_{|\xi|\ge1}|\xi|^{s(p-1)-N-sp}d\xi+r^{sp-s(p-1)}\int_{|y-x|>\frac12}\frac{U(y)^{p-1}}{|y-x|^{N+sp}}dy\le C ,
\]
because $s(p-1)<sp$, $r^{s}\le1$, and $|y-x|\ge c(1+|y|)$ for $|y-x|>\frac12$, $|x|=1$, while $U\in\Ls$.\\

{Next, by Remark~}\ref{rem:K1}, applied to $U_r$ with $D=B_{5/4}$ and $\delta=\frac18$ (the only singular point of $U_r$ is $-x/r$, and $|x/r|=2/d\ge2$, so it lies outside $\overline{B_{3/2}}$), cutting $U_r$ off outside $B_{3/2}$ only adds to $f_r$ a bounded continuous term, bounded independently of $x$ by Lemma~\ref{lem:loc}(ii) and the bounds above. By (K1), the cut-off function $\tilde U_r$, which coincides with $U_r$ near $\overline{B_{5/4}}$ and satisfies $0\le\tilde U_r\le U_r$, is a bounded local weak solution in $B_{5/4}$ of an equation $\pl\tilde U_r=\tilde f_r$ with $\|\tilde f_r\|_{L^\infty(B_{5/4})}\le C$, and its tail is bounded by that of $U_r$. By (K8) (the part of the tail in $B_1\setminus B_{\changed{1/2}}$ is bounded by $\|U_r\|_{L^\infty(B_1)}$), $[U_r]_{C^\alpha(B_{1/32})}=[\tilde U_r]_{C^\alpha(B_{1/32})}\le C_1'$ with $C_1'$ independent of $x$. For $|y-x|<\frac d{64}=\frac r{32}$,
\[
|U(y)-U(x)|=r^s\Big|U_r\Big(\frac{y-x}r\Big)-U_r(0)\Big|\le C_1'r^{s-\alpha}|y-x|^\alpha\le 2^{\alpha}C_1'\,d^{s-\alpha}|y-x|^\alpha .
\]
\end{proof}

\subsection{The comparison function}
For $a\in(\frac12,1),$ {set}
\[
\varepsilon:=2\beta^+(1-a),\qquad Z:=U^a,\qquad h(y):=|y|^{-\varepsilon},\qquad W:=hZ,\qquad\gamma:=a\beta^++\varepsilon=\beta^++(1-a)\beta^+ .
\]
Then $Z\in H_{a\beta^+}$, $W\in H_\gamma$ is continuous, nonnegative, positive in $\Cw$ and zero outside $\Cw$, and $W\in\Ls$ as soon as $\gamma<\frac N{p-1}$. By Lemma~\ref{lem:range}(b), $\beta^+(p-1)\le N-\delta_0(p-1)$, so this holds for $a$ close to $1$. We {prove that} $W$ is a viscosity supersolution of $\pl W\ge0$ in $\Cw$ for $a$ close to $1$; then $W\in\mathcal S_\gamma(\omega)$ with $\gamma>\beta^+$, which contradicts the definition of $\beta^+$ and proves Theorem~\ref{mt}(d).

{At a heuristic level,} $\pl W\approx\pl Z$ up to an error of order $\varepsilon$ (Lemma~\ref{lem:error}), while by Lemma~\ref{lem:power} $\pl Z\ge a^{p-1}\mu Z^{p-1}|x|^{-sp}+(\text{concavity gain})$; the first term is bounded below away from $\partial\Cw$, and the gain is large near $\partial\Cw$ (Lemma~\ref{lem:gain}).

For $x\in\Cw\cap\Sp$ we set, for $0<r'<d(x)/64$ (first line) and for $r'>0$ (second line),
\begin{gather*}
J_{r'}(x):=\int_{|y-x|\ge r'}\frac{\Psi_p\big(Z(x)-h(y)Z(y)\big)-\Psi_p\big(Z(x)-Z(y)\big)}{|x-y|^{N+sp}}dy,\\
G_{r'}(x):=\int_{|y-x|\ge r'}\frac{\Gamma_a(x,y)}{|x-y|^{N+sp}}dy,
\end{gather*}
where $\Gamma_a\ge0$ is defined in Lemma~\ref{lem:power} with $v=U$.

Let $a_1\in(\frac12,1)$ be such that, for every $a\in[a_1,1)$,
\begin{equation}\label{eq:a1}
1+as(p-1)-sp\ge\tfrac{1-s}2,\qquad \gamma(p-1)\le N-\tfrac{\delta_0(p-1)}2,\qquad \varepsilon\le\min\big\{1,\tfrac{\delta_0(p-1)}4\big\}.
\end{equation}
This is possible because $1+s(p-1)-sp=1-s>0$, $\gamma\to\beta^+$ and $\varepsilon\to0$ as $a\to1$, and $\beta^+(p-1)\le N-\delta_0(p-1)$. {Throughout the remainder of this section,} $C$ denotes constants independent of $a\in[a_1,1)$, $x\in\Cw\cap\Sp$ and $r'$.

{Since} $U\le1$ on $\Sp$, \eqref{eq:twosided} gives, with $C_Z:=C_0$,
\begin{equation}\label{eq:Zbounds}
Z(y)\le C_Z\,d(y)^{as}\quad\text{for }y\in\Cw,\ \tfrac12\le|y|\le2;\qquad Z(y)\le|y|^{-a\beta^+}\quad\text{for all }y .
\end{equation}
For $\tau>0$, using $|e^u-1|\le|u|\max\{1,e^u\}$ with $u=-\varepsilon\log\tau$, we obtain
\begin{equation}\label{eq:hbounds}
|1-\tau^{-\varepsilon}|\le\varepsilon|\log\tau|\max\{1,\tau^{-\varepsilon}\};\qquad\text{in particular}\quad |1-\tau^{-\varepsilon}|\le4\varepsilon|\tau-1|\ \text{ for }\tfrac12\le\tau\le\tfrac32,
\end{equation}
because then $\tau^{-\varepsilon}\le2^\varepsilon\le2$ and $|\log\tau|\le2|\tau-1|$. Since $\big||y|-1\big|\le|y-x|$ for $|x|=1$, {we obtain} $|1-h(y)|\le4\varepsilon|y-x|$ if $|y-x|\le\frac12$.

\begin{lemma}[Error estimate]\label{lem:error}
There is $C_E$ such that $|J_{r'}(x)|\le C_E\,\varepsilon$ for all $a\in[a_1,1)$, $x\in\Cw\cap\Sp$ and $0<r'<d(x)/64$.
\end{lemma}

\begin{proof}
Fix $x\in\Cw\cap\Sp$, let $d:=d(x)\le1$, and write $z=y-x$,
\[
A:=Z(x)-Z(y),\qquad B:=(1-h(y))Z(y),
\]
so that $Z(x)-h(y)Z(y)=A+B$ and the integrand of $J_{r'}$ is $[\Psi_p(A+B)-\Psi_p(A)]|z|^{-N-sp}$. \mathred{Write $\mathcal E(z):=|\Psi_p(A+B)-\Psi_p(A)|\,|z|^{-N-sp}$. }

\medskip
\noindent\emph{Case $p\ge2$:} By the mean value theorem, since $|\Psi_p'(\tau)|=(p-1)|\tau|^{p-2}$ and $|\tau|\le|A|+|B|$ on the segment between $A$ and $A+B$,
\begin{equation}\label{eq:mvt}
|\Psi_p(A+B)-\Psi_p(A)|\le(p-1)|B|\,(|A|+|B|)^{p-2}.
\end{equation}
{We divide }$\{|z|\ge r'\}$ into three regions.

\emph {Region I $(|z|>\frac12).$} By \eqref{eq:hbounds} and \eqref{eq:Zbounds}, $|B|\le\varepsilon|\log|y||(1+|y|^{-\varepsilon})|y|^{-a\beta^+}\le2\varepsilon|\log|y||(1+|y|^{-\gamma})$ and, since also $|1-h(y)|\le1+|y|^{-\varepsilon}$, $|A|+|B|\le Z(x)+Z(y)+|B|\le1+2|y|^{-a\beta^+}+|y|^{-\gamma}\le4(1+|y|^{-\gamma})$, using $|y|^{-a\beta^+}\le1+|y|^{-\gamma}$. Moreover $|x-y|\ge\frac16(1+|y|)$ when $|x-y|>\frac12$. Hence, by \eqref{eq:mvt},
\[
\int_{|z|>\frac12}\mathred{\mathcal E(z)\,dz}\le C\varepsilon\int_{\R^N}\frac{|\log|y||\,(1+|y|^{-\gamma(p-1)})}{(1+|y|)^{N+sp}}\,dy\le C\varepsilon,
\]
since $\gamma(p-1)\le N-\frac{\delta_0(p-1)}2$ by \eqref{eq:a1}.

\emph{Region II $(\frac d{64}\le|z|\le\frac12).$} Here $\frac12\le|y|\le\frac32$ and $d(y)\le d+|z|\le65|z|$, so by \eqref{eq:Zbounds} $Z(y)\le C_Z65^{as}|z|^{as}\le C|z|^{as}$; also $Z(x)\le C_Zd^{as}\le C|z|^{as}$. Hence $|A|+|B|\le3(Z(x)+Z(y))\le C|z|^{as}$, and $|B|\le4\varepsilon|z|Z(y)\le C\varepsilon|z|^{1+as}$. By \eqref{eq:mvt},
\[
\int_{\frac d{64}\le|z|\le\frac12}\mathred{\mathcal E(z)\,dz}\le C\varepsilon\int_{|z|\le\frac12}|z|^{1+as(p-1)-N-sp}dz=C\varepsilon\,\mathcal H^{N-1}(\Sp)\int_0^{1/2}t^{as(p-1)-sp}dt\le C\varepsilon ,
\]
because $as(p-1)-sp\ge\frac{1-s}2-1>-1$ by \eqref{eq:a1}.

\emph{Region III $(r'\le|z|<\frac d{64})$.} Here $d(y)\ge d-|z|\ge\frac d2$, $d(y)\le2d$ and $\frac12\le|y|\le\frac32$, so $U(y),U(x)\ge c_0(d/2)^s$. By the mean value theorem applied to $\tau\mapsto\tau^a$ (whose derivative $a\tau^{a-1}$ is decreasing) and Lemma~\ref{lem:gradU},
\[
|A|=|U(x)^a-U(y)^a|\le a\big(c_0(d/2)^s\big)^{a-1}|U(x)-U(y)|\le C\,d^{s(a-1)}\,d^{s-\alpha}|z|^\alpha=C\,d^{as-\alpha}|z|^\alpha ,
\]
where $(c_0 2^{-s})^{a-1}\le(c_02^{-s})^{-1/2}$ since $c_02^{-s}\le1$ and $a-1\in(-\frac12,0)$. Also $|B|\le4\varepsilon|z|Z(y)\le C\varepsilon|z|d^{as}$, and since $|z|\le d$, $|z|=|z|^\alpha|z|^{1-\alpha}\le|z|^\alpha d^{1-\alpha}$, so $|B|\le C d^{as-\alpha}|z|^\alpha$ (as $\varepsilon\le1$, $d\le1$). By \eqref{eq:mvt},
\[
\int_{r'\le|z|<\frac d{64}}\mathred{\mathcal E(z)\,dz}\le C\varepsilon\,d^{as}\,d^{(as-\alpha)(p-2)}\int_{|z|<d/64}|z|^{1+\alpha(p-2)-N-sp}dz .
\]
By \eqref{eq:alpha}, $1+\alpha(p-2)-sp>0$, so {the final integral is} $C(d/64)^{1+\alpha(p-2)-sp}$, and the total is at most
\[
C\varepsilon\,d^{as+(as-\alpha)(p-2)+1+\alpha(p-2)-sp}=C\varepsilon\,d^{1+as(p-1)-sp}\le C\varepsilon ,
\]
since the exponent is positive by \eqref{eq:a1} and $d\le1$.

\medskip
\noindent\emph{Case $1<p<2$:} {Here, no regularity of} $U$ is needed; instead we use that $Z$ is strictly monotone along rays. First, for $A\neq0$,
\begin{equation}\label{eq:mvt2}
|\Psi_p(A+B)-\Psi_p(A)|\le 4\,|B|\,|A|^{p-2}.
\end{equation}
Indeed, if $|A|\ge2|B|$, then $|\tau|\ge|A|/2$ on the segment between $A$ and $A+B$, and since $|\Psi_p'(\tau)|=(p-1)|\tau|^{p-2}$ is decreasing in $|\tau|$, the mean value theorem gives the bound $(p-1)2^{2-p}|B||A|^{p-2}\le2|B||A|^{p-2}$. If $|A|<2|B|$, then by \eqref{eq:elem} $|\Psi_p(A+B)-\Psi_p(A)|\le2^{2-p}|B|^{p-1}=2^{2-p}|B||B|^{p-2}\le2^{2(2-p)}|B||A|^{p-2}\le4|B||A|^{p-2}$, because $|B|>|A|/2$ and $p-2<0$.

We use polar coordinates $y=t\theta$, $t>0$, $\theta\in\Sp$, $dy=t^{N-1}dt\,d\theta$. Since $Z(t\theta)=t^{-a\beta^+}Z(\theta)$,
\[
A(t\theta)=Z(x)-t^{-a\beta^+}Z(\theta),\qquad \partial_tA(t\theta)=a\beta^+t^{-a\beta^+-1}Z(\theta)\ge0,\]
\[|B(t\theta)|=|1-t^{-\varepsilon}|\,t^{-a\beta^+}Z(\theta)=\frac{|1-t^{-\varepsilon}|\,t}{a\beta^+}\,\partial_tA(t\theta).
\]
If $Z(\theta)=0$, then $B(t\theta)=0$ and the integrand vanishes. Otherwise $t\mapsto A(t\theta)$ is strictly increasing, it vanishes for at most one $t$, and for any interval $I\subset(0,\infty)$ and any measurable $g\ge0$, \rev{the change of variables $\tau=A(t\theta)$ gives, for every measurable $I'\subset I$,}
{\color{\revisioncolor}
\begin{equation}\label{eq:coarea}
\begin{aligned}
\int_{I'} g(t)\,|A|^{p-2}\partial_tA\,dt
&\le\Big(\sup_{I'}g\Big)\int_{A(I)}|\tau|^{p-2}d\tau\\
&\le\frac{2^{2-p}}{p-1}\Big(\sup_{I'}g\Big)\,|A(I)|^{p-1},
\end{aligned}
\end{equation}
}
where $|A(I)|$ is the length of the interval $A(I)$ (the integral of $|\tau|^{p-2}$ over an interval of length $\ell$ is at most $\frac{2^{2-p}}{p-1}\ell^{p-1}$). By \eqref{eq:mvt2}, in the measure $dt\,d\theta$ the absolute value of the integrand of $J_{r'}$ is bounded by
\[
\frac4{a\beta^+}\,g(t,\theta)\,|A|^{p-2}\partial_tA,\qquad g(t,\theta):=|1-t^{-\varepsilon}|\,t^{N}|x-t\theta|^{-N-sp},
\]
and $\frac4{a\beta^+}\le\frac8{\beta^+}$. In each region below we bound $g$ by its supremum over the region, and then we enlarge the $t$-integration to an interval containing the region, which is allowed since $|A|^{p-2}\partial_tA\ge0$.

{We now introduce} the {dyadic shells $S_\rho:=\{\rho\le|z|<2\rho\}$, with $\rho=2^{-j-2}$, $j\ge0$, $2\rho>r'$.} These shells cover $\{r'\le|z|<\frac12\}$. For $y=t\theta\in S_\rho$ we have $|t-1|\le|z|<2\rho\le\frac12$ and $|\theta-x|\le|\theta-y|+|y-x|=|1-t|+|z|<4\rho$. Hence, by \eqref{eq:hbounds}, $g\le4\varepsilon|t-1|(3/2)^N\rho^{-N-sp}\le C\varepsilon\rho^{1-N-sp}$ on $S_\rho$. For $\theta$ fixed, the $t$-section of $S_\rho$ is contained in $I:=(1-2\rho,1+2\rho)$, and
\[
|A(I)|=Z(\theta)\big((1-2\rho)^{-a\beta^+}-(1+2\rho)^{-a\beta^+}\big)\le C\rho\,Z(\theta)\le C\rho\,(d+4\rho)^{as},
\]
by the mean value theorem ($1-2\rho\ge\frac12$) and \eqref{eq:Zbounds} with $d(\theta)\le d+|\theta-x|$. The set of $\theta$ involved has $\mathcal H^{N-1}$-measure at most $C\rho^{N-1}$. By \eqref{eq:coarea}, the contribution of $S_\rho$ is at most
\[
C\varepsilon\rho^{1-N-sp}\cdot\rho^{N-1}\cdot\big(\rho\,(d+4\rho)^{as}\big)^{p-1}=C\varepsilon\,\rho^{p-1-sp}(d+4\rho)^{as(p-1)} .
\]
{Summing over} $\rho=2^{-j-2}$. The shells with $\rho\le d$ give at most $C\varepsilon(5d)^{as(p-1)}\sum_{\rho\le d}\rho^{p-1-sp}\le C\varepsilon\,d^{as(p-1)+p-1-sp}\le C\varepsilon$, where the geometric sum converges because $p-1-sp>0$; the shells with $\rho>d$ give at most $\mathred{C\varepsilon\sum_{\rho>d}\rho^{p-1-sp}(5\rho)^{as(p-1)}\le C\varepsilon}$. {This is the only point at which} \eqref{eq:range-p} is used in the case $p<2$.

\emph{Region $|z|\ge\frac12$, $t\ge\frac12$.} Here $|x-t\theta|\ge\frac16(1+t)$, and by \eqref{eq:hbounds} $|1-t^{-\varepsilon}|\le2\varepsilon|\log t|$, so $g\le C\varepsilon\sup_{t\ge1/2}|\log t|\,t^N(1+t)^{-N-sp}\le C\varepsilon$. The map $t\mapsto A(t\theta)$ sends $[\frac12,\infty)$ into $[Z(x)-2^{a\beta^+}Z(\theta),Z(x)]$, an interval of length at most $2^{\beta^+}$ (recall $Z\le1$ on $\Sp$). By \eqref{eq:coarea}, for each $\theta$ the $t$-integral is at most $C\varepsilon$, and integrating over $\theta\in\Sp$ gives $C\varepsilon$.

\emph{Region $t<\frac12$} (which is contained in $\{|z|>\frac12\}$). Here $|x-t\theta|\ge\frac12$. Let $\varkappa:=\frac{\delta_0(p-1)}4$. Since $e^u-1\le ue^u$ for $u\ge0$ and $\log\frac1t\le\frac1{e\varkappa}t^{-\varkappa}$,
\[
|1-t^{-\varepsilon}|=t^{-\varepsilon}-1\le\varepsilon\log\tfrac1t\;t^{-\varepsilon}\le\frac{\varepsilon}{e\varkappa}\,t^{-\varepsilon-\varkappa},\qquad\text{so}\qquad g\le C\varepsilon\,t^{\varsigma},\quad \varsigma:=N-\varepsilon-\varkappa .
\]
By \eqref{eq:a1}, $\varepsilon\le\varkappa$, and since $a\beta^+(p-1)\le\beta^+(p-1)\le N-4\varkappa$, we get $\varsigma-a\beta^+(p-1)\ge2\varkappa$. Fix $\theta$ with $Z(\theta)>0$ and change variables $\tau=A(t\theta)$; as $t$ ranges in $(0,\frac12)$, $\tau$ ranges in an interval contained in $(-\infty,1]$. If $\tau\le-2$, then $t^{-a\beta^+}\ge t^{-a\beta^+}Z(\theta)=Z(x)-\tau\ge|\tau|$, i.e.\ $t\le|\tau|^{-1/(a\beta^+)}$, and so $g\le C\varepsilon|\tau|^{-\varsigma/(a\beta^+)}$. Therefore
\[
\int_0^{1/2}g\,|A|^{p-2}\partial_tA\,dt\le C\varepsilon\Big(\int_{-2}^{1}|\tau|^{p-2}d\tau+\int_{2}^\infty\tau^{p-2-\frac{\varsigma}{a\beta^+}}d\tau\Big)\le C\varepsilon,
\]
because $p-2-\frac{\varsigma}{a\beta^+}=-1-\frac{\varsigma-a\beta^+(p-1)}{a\beta^+}\le-1-\frac{2\varkappa}{\beta^+}$. Integrating over $\theta\in\Sp$ we obtain the bound $C\varepsilon$.

{Combining the contributions from all regions in both cases, we obtain} $|J_{r'}(x)|\le C_E\varepsilon$, with $C_E$ independent of $a\in[a_1,1)$, $x$ and $r'$.
\end{proof}

\begin{lemma}[Gain near the boundary]\label{lem:gain}
There are $a_2\in(\frac12,1)$, $\delta_1\in(0,\frac12]$, $\rho_Q>0$ and $c_G,c_3>0$ such that for all $a\in[a_2,1)$, all $x\in\Cw\cap\Sp$ with $U(x)\le\delta_1$ and all $0<r'<\rho_Q$,
\[
G_{r'}(x)\ge\min\Big\{c_G(1-a)\log\frac{\delta_1}{U(x)},\ c_3\Big\}.
\]
\end{lemma}

\begin{proof}
%\emph{Setting.}
 Fix a smooth domain $\omega_1$ with $\omega_0\Subset\omega_1\Subset\omega$, and let $Q:=\mathcal C_{\omega_0}\cap\{\frac12\le|y|\le2\}$, $m_Q:=\min_{\overline Q}U>0$, $M_Q:=\max_{\overline Q}U$, and
\[
\delta_1:=\min\Big\{\tfrac12\min_{\Sp\cap\overline{\mathcal C_{\omega_1}}}U,\ m_Qe^{-2}\Big\}\le\tfrac12 .
\]
If $U(x)\le\delta_1$, then $x\notin\overline{\mathcal C_{\omega_1}}$, so $|x-y|\ge\rho_Q:=\dist(\Sp\setminus\mathcal C_{\omega_1},\overline Q)>0$ for $y\in Q$. For $r'<\rho_Q$ we have $Q\subset\{|y-x|\ge r'\}$, and since $\Gamma_a\ge0$ and $|x-y|\le3$ for $y\in Q$,
\[
G_{r'}(x)\ge3^{-N-sp}\int_Q\Gamma_a(x,y)\,dy .
\]

%\emph{The integrand.}
 Let $y\in Q$, $u_0:=U(x)\le\delta_1$, $u_1:=U(y)\in[m_Q,M_Q]$; so $u_1/u_0\ge e^2$. Since $u_0<u_1$,
\[
\Gamma_a(x,y)=-(u_1^a-u_0^a)^{p-1}+a^{p-1}u_0^{(a-1)(p-1)}(u_1-u_0)^{p-1}=D_2^{p-1}-D_1^{p-1},
\]
with $D_1:=u_1^a-u_0^a$ and $D_2:=a\,u_0^{a-1}(u_1-u_0)\ge D_1$ (by \eqref{eq:concave}). Moreover
\[
c_D:=(1-e^{-1})\min\{1,m_Q\}\le u_1^a(1-e^{-2a})\le D_1\le u_1^a\le\max\{1,M_Q\}=:C_D ,
\]
using $(u_0/u_1)^a\le e^{-2a}\le e^{-1}$ and $u_1^a\ge\min\{1,u_1\}$.

\changed{{We begin by showing that}} $\Gamma_a\ge c_\Gamma\min\{\Delta,1\}$, where $\Delta:=D_2-D_1\ge0$. If $p\ge2$, convexity of $\tau\mapsto\tau^{p-1}$ gives $\Gamma_a\ge(p-1)D_1^{p-2}\Delta\ge(p-1)c_D^{p-2}\Delta$. If $p<2$ and $\Delta\le D_1$, the mean value theorem on $[D_1,D_2]\subset[D_1,2D_1]$ gives $\Gamma_a\ge(p-1)(2D_1)^{p-2}\Delta\ge(p-1)(2C_D)^{p-2}\Delta$. If $p<2$ and $\Delta>D_1$, then $D_2>2D_1$ and $\Gamma_a\ge(2^{p-1}-1)D_1^{p-1}\ge(2^{p-1}-1)c_D^{p-1}$. So {the claim follows with} $c_\Gamma>0$ depending only on $p,m_Q,M_Q$.

%\emph{Step 2: a formula for $\Delta$.}
 With $r:=u_1/u_0\ge e^2$, we have
\[
\Delta=u_0^a\,\zeta_a(r),\qquad \zeta_a(r):=a(r-1)-r^a+1=a\int_1^r(1-t^{a-1})\,dt ,
\]
and, as a function of $u_0\in(0,u_1)$, $\Delta$ is decreasing: $\partial_{u_0}\big[a u_0^{a-1}(u_1-u_0)-u_1^a+u_0^a\big]=a(a-1)u_0^{a-2}(u_1-u_0)<0$.

\changed{{First, suppose that}} $(1-a)\log(M_Q/u_0)\le1$. For $1\le t\le r\le M_Q/u_0$ we have $0\le(1-a)\log t\le1$, and $1-e^{-\sigma}\ge\frac{\sigma}{1+\sigma}\ge\frac\sigma2$ for $\sigma\in[0,1]$, so $1-t^{a-1}\ge\frac12(1-a)\log t$. Hence, using $\int_1^r\log t\,dt=r\log r-r+1\ge\frac12r\log r$ (as $\log r\ge2$), $u_0\le1$ and $a\ge\frac12$,
\[
\Delta\ge\frac a2(1-a)u_0^a\,(r\log r-r+1)\ge\frac a4(1-a)u_0^{a-1}u_1\log\frac{u_1}{u_0}\ge\frac18(1-a)\,m_Q\log\frac{m_Q}{u_0}.
\]

\changed{{Next, suppose that}} $(1-a)\log(M_Q/u_0)>1$. Then $u_0<\bar u:=M_Qe^{-1/(1-a)}$. Let $a_2\in(\frac12,1)$ be such that, for $a\ge a_2$, $\bar u\le\delta_1$ (so that $\bar u\le1$ and $\bar u\le m_Qe^{-2}$) and $(1-a)\log\frac{M_Q}{m_Q}\le\frac12$. \changed{{The preceding lower bound also applies}} at $u_0=\bar u$ (for which $(1-a)\log(M_Q/\bar u)=1$ and $u_1/\bar u\ge e^2$), and by the monotonicity \changed{just established},
\[
\Delta\ge\Delta|_{u_0=\bar u}\ge\frac18(1-a)m_Q\Big(\log\frac{m_Q}{M_Q}+\frac1{1-a}\Big)=\frac{m_Q}8\Big(1-(1-a)\log\frac{M_Q}{m_Q}\Big)\ge\frac{m_Q}{16}.
\]

%\emph{Conclusion.}
 In both cases $\Delta\ge\min\{\frac18(1-a)m_Q\log\frac{m_Q}{u_0},\frac{m_Q}{16}\}$, and $\log\frac{m_Q}{u_0}\ge\log\frac{\delta_1}{u_0}$. By \changed{{the estimate established above}}, for all $y\in Q$,
\[
\Gamma_a(x,y)\ge c_\Gamma\min\Big\{\frac{m_Q}8(1-a)\log\frac{\delta_1}{U(x)},\ \frac{m_Q}{16},\ 1\Big\},
\]
\rev{{Integrating over} $Q$ gives the claim with}
{\color{\revisioncolor}
\[
\begin{aligned}
c_G&:=3^{-N-sp}|Q|c_\Gamma\frac{m_Q}8,\\
c_3&:=3^{-N-sp}|Q|c_\Gamma\min\Big\{\frac{m_Q}{16},1\Big\}.
\end{aligned}
\]
}
\end{proof}

The case $p=2$ of Theorem~\ref{mt}(d), for every $s\in(0,1)$, is proved in Corollary~\ref{cor:p2}; here we assume \eqref{eq:range-p}.

\begin{proof}[Proof of Theorem~\ref{mt}\textup{(d)} when $p-1>sp$]
%\emph{Choice of $a$.}
 Let
\[
\delta_*:=\delta_1\exp\Big(-\frac{4C_E\beta^+}{c_G}\Big)\le\delta_1,
\]
and choose $a\in[\max\{a_1,a_2\},1)$ so close to $1$ that (recall that $\gamma<\frac N{p-1}$ by \eqref{eq:a1})
\begin{equation}\label{eq:choice-a}
2C_E\varepsilon\le\min\big\{c_3,\ 2^{1-p}\mu\,\delta_*^{p-1}\big\}\qquad(\varepsilon=2\beta^+(1-a)).
\end{equation}

%\emph{Reduction.}
 By homogeneity (Section~\ref{sec:prel}) it suffices to verify the supersolution property at points $x\in\Cw\cap\Sp$. Let $\phi$ be admissible and touch $W$ from below at $x$ in $B_r(x)\subset\Cw$. By Lemma~\ref{lem:grad} (with $W\in H_\gamma$, $W(x)>0$), $\nabla\phi(x)\neq0$. Since $h$ is smooth and positive near $x$ with $h(x)=1$, $\varphi:=\phi/h$ is $C^2$ near $x$ and touches $Z$ from below at $x$ in $B_r(x)$, and, again by Lemma~\ref{lem:grad} (with $Z\in H_{a\beta^+}$), $\nabla\varphi(x)\neq0$. Let $0<r'<\min\{r,d(x)/64,\rho_Q\}$. By Lemma~\ref{lem:grad}, $\pl\phi_r(x)\ge\pl\phi_{r'}(x)$.

%\emph{Splitting.}
 Let $\varphi_{r'}:=\varphi$ in $B_{r'}(x)$, $\varphi_{r'}:=Z$ outside. Then $\phi_{r'}=h\varphi_{r'}$ in $\R^N$ and $\phi_{r'}(x)=\varphi(x)=Z(x)$. Comparing the integrals defining $\pl\phi_{r'}(x)$ and $\pl\varphi_{r'}(x)$ inside and outside $B_{r'}(x)$,
\[
\pl\phi_{r'}(x)=\pl\varphi_{r'}(x)+J_{r'}(x)+I_{r'},\qquad I_{r'}:=P_1(r')-P_2(r'),
\]
where $P_1(r'):=\mathrm{P.V.}\int_{B_{r'}(x)}\Psi_p(\varphi(x)-h\varphi(y))|x-y|^{-N-sp}dy$ and $P_2(r')$ is the same expression with $\varphi$ in place of $h\varphi$. The functions $F_1:=h\varphi-\varphi(x)$ and $F_2:=\varphi-\varphi(x)$ are $C^2$ near $x$, vanish at $x$, and $\nabla F_1(x)=\nabla\phi(x)\ne0$, $\nabla F_2(x)=\nabla\varphi(x)\ne0$. By Lemma~\ref{lem:small}, $P_1(r'),P_2(r')\to0$ as $r'\to0$, so $I_{r'}\to0$.

%\emph{Lower bound for $\pl\varphi_{r'}(x)$.}
{By Lemma~}\ref{lem:power} with $v=U$, $D=\Cw$ and $F=\mu|x|^{-sp}U^{p-1}$, since $|x|=1$ and $a^{p-1}U(x)^{(a-1)(p-1)}U(x)^{p-1}=a^{p-1}Z(x)^{p-1}$,
\[
\pl\varphi_{r'}(x)\ge a^{p-1}\mu Z(x)^{p-1}+G_{r'}(x),
\]
and both terms on the right are nonnegative. We claim that their sum is at least $2C_E\varepsilon$:
\begin{itemize}
\item if $U(x)>\delta_*$, then, as $U(x)\le1$ and $a<1$, $Z(x)=U(x)^a\ge U(x)>\delta_*$, so $a^{p-1}\mu Z(x)^{p-1}\ge2^{1-p}\mu\delta_*^{p-1}\ge2C_E\varepsilon$ by \eqref{eq:choice-a};
\item if $U(x)\le\delta_*\le\delta_1$, then $\log\frac{\delta_1}{U(x)}\ge\frac{4C_E\beta^+}{c_G}$, and Lemma~\ref{lem:gain} gives
\[
G_{r'}(x)\ge\min\Big\{c_G(1-a)\frac{4C_E\beta^+}{c_G},\ c_3\Big\}=\min\{2C_E\varepsilon,\ c_3\}=2C_E\varepsilon .
\]
\end{itemize}

%\emph{Conclusion.}
{Using Lemma~}\ref{lem:error}, we deduce
\[
\pl\phi_r(x)\ge\pl\phi_{r'}(x)\ge2C_E\varepsilon-C_E\varepsilon+I_{r'}=C_E\varepsilon+I_{r'}\ \xrightarrow[r'\to0]{}\ C_E\varepsilon>0 .
\]
Hence $\pl\phi_r(x)\ge0$, and $W$ is a viscosity supersolution of $\pl W\ge0$ in $\Cw$. Since $W\in H_\gamma\cap\Ls$, $W\ge0$ and $W>0$ in $\Cw$, we get $W\in\mathcal S_\gamma(\omega)$ with $\gamma>\beta^+$, contradicting the definition of $\beta^+$. Therefore $\underline\mu=0$, and $\Phi^+:=\underline U$ is a fundamental solution with exponent $\beta^+$.
\end{proof}

\section{Liouville\changed{-}type results}\label{sec:liou}

\begin{proof}[Proof of Theorem~\ref{thm:liou}\textup{(a)}]
Assume \rev{$u\not\equiv0$ in $\R^N\setminus\{0\}$}. Since $u^q\ge0$, $u$ is a supersolution of $\pl u\ge0$ in $\Cw$; being continuous, nonnegative and not identically zero, it is positive on a set of positive measure, so $u>0$ in $\Cw$ by Lemma~\ref{lem:smp}.

%\emph{Scaling.}
 For $R>0$, {set} $u_R(x):=u(Rx)$. By scaling, $u_R$ is a viscosity supersolution of $\pl u_R\ge R^{sp}u_R^q$ in $\Cw$ (indeed $\pl u_R(x)=R^{sp}(\pl u)(Rx)\ge R^{sp}u(Rx)^q$, and $\Cw$ is invariant under dilations).

%\emph{The basic inequality.}
 Fix $e\in\omega_0$ (so $e\in K_0$ and $B_{4r_0}(e)\subset\Cw$; recall $r_0\le\frac18$). Let $\eta_2$ be the function of Lemma~\ref{lem:bump} with radius $2r_0$, and $\tilde\eta:=\eta_2(\cdot-e)$. By Remark~\ref{rem:touch} (with $w=u_R$, $B=B_{2r_0}(e)$, $F=R^{sp}u_R^q$) there are $c^*>0$ and $x_1\in B_{2r_0}(e)$ such that $\phi:=c^*\tilde\eta$ touches $u_R$ from below at $x_1$, $u_R\ge c^*\tilde\eta$ in $\R^N$, and
\begin{equation}\label{eq:L1}
R^{sp}u_R(x_1)^q\le(c^*)^{p-1}M_{\eta_2} .
\end{equation}
Next, $u_R$ is a nonnegative supersolution of $\pl u_R\ge0$ in $B_{4r_0}(e)$, so by Remark~\ref{rem:K1} and the weak Harnack inequality (K4),
\[
u_R(x_1)\ge\inf_{B_{2r_0}(e)}u_R\ge c_H\Big(\frac1{|B_{2r_0}|}\int_{B_{2r_0}(e)}u_R^t\Big)^{1/t}\ge c_H\Big(\frac1{|B_{2r_0}|}\int_{B_{2r_0}(e)}(c^*\tilde\eta)^t\Big)^{1/t}=c_4\,c^*,
\]
where $c_4:=c_H\big(\frac1{|B_1|}\int_{B_1}\eta_1^t\big)^{1/t}>0$ is independent of $R$ (by scaling). Hence $c^*\le u_R(x_1)/c_4$, and \eqref{eq:L1} gives $R^{sp}u_R(x_1)^q\le M_{\eta_2}c_4^{1-p}u_R(x_1)^{p-1}$; since $u_R(x_1)>0$,
\begin{equation}\label{eq:L2}
R^{sp}\,u_R(x_1)^{q-p+1}\le C_4:=M_{\eta_2}c_4^{1-p}\qquad\text{for all }R>0 ,
\end{equation}
where $x_1=x_1(R)\in B_{2r_0}(e)$. We also note that
\begin{equation}\label{eq:L3}
c^*\ge\min_{\overline{B_{2r_0}(e)}}u_R\qquad\text{and}\qquad c^*=\frac{u_R(x_1)}{\tilde\eta(x_1)}\ge u_R(x_1),
\end{equation}
the first because $\tilde\eta\le1$, the second because $\phi(x_1)=u_R(x_1)$ and $\tilde\eta\le1$.

{\it Case $q=p-1$:} \eqref{eq:L2} reads $R^{sp}\le C_4$ for all $R>0$, {which is impossible}.

{\it Case $p-1<q<q^*$:} Since $q<q^*=p-1+\frac{sp}{\beta^+}$, we have $\frac{sp}{q-p+1}>\beta^+$, and we can choose $\beta$ with
\[
\beta^+<\beta<\min\Big\{\frac N{p-1},\frac{sp}{q-p+1}\Big\},\qquad\text{so that}\qquad sp-\beta(q-p+1)>0 .
\]
By Proposition~\ref{prop:eig} and Corollary~\ref{cor:sign} there is an eigenpair $(\mu_\beta,U_\beta)$ with $\mu_\beta<0$; thus $U_\beta$ is a subsolution of $\pl U_\beta\le0$ in $\Cw$, nonnegative, nontrivial, vanishing outside $\Cw$. Proposition~\ref{prop:nonex}, applied with $v=u$, gives $\varepsilon>0$ with $u\ge\varepsilon U_\beta$ in $\Cw\setminus B_2$. Let $R\ge4$ and $y\in R\,\overline{B_{2r_0}(e)}$. Then $y/|y|\in K_e:=\{\xi/|\xi|:\ \xi\in\overline{B_{2r_0}(e)}\}$, which is a compact subset of $\Cw\cap\Sp$ (because $\overline{B_{2r_0}(e)}\subset B_{4r_0}(e)\subset\Cw$ and $\Cw$ is a cone), and $2<R(1-2r_0)\le|y|\le R(1+2r_0)\le2R$. Hence
\[
u(y)\ge\varepsilon U_\beta(y)=\varepsilon|y|^{-\beta}U_\beta(y/|y|)\ge\varepsilon(2R)^{-\beta}\min_{K_e}U_\beta=:c_1R^{-\beta},
\]
with $c_1>0$ since $U_\beta>0$ on the compact set $K_e$. By \eqref{eq:L3}, $c^*\ge\min_{\overline{B_{2r_0}(e)}}u_R\ge c_1R^{-\beta}$, and so $u_R(x_1)\ge c_4c^*\ge c_4c_1R^{-\beta}$. Since $q-p+1>0$, \eqref{eq:L2} gives
\[
(c_4c_1)^{q-p+1}R^{sp-\beta(q-p+1)}\le C_4\qquad\text{for all }R\ge4,
\]
{which is impossible} because $sp-\beta(q-p+1)>0$.

{\it Case $0<q<p-1$:} Now $q-p+1<0$, and \eqref{eq:L2} gives $u_R(x_1)\ge(R^{sp}/C_4)^{1/(p-1-q)}$. Let $\rho_1\in(0,2r_0)$ be such that $\tilde\eta\ge\frac12$ on $B_{\rho_1}(e)$. By \eqref{eq:L3}, on $B_{\rho_1}(e)$ we have $u_R\ge c^*\tilde\eta\ge\frac12c^*\ge\frac12u_R(x_1)$. Going back to $u$: for every $R>0$,
\[
u\ge c_5R^{\frac{sp}{p-1-q}}\quad\text{on }B_{R\rho_1}(Re),\qquad c_5:=\tfrac12C_4^{-1/(p-1-q)} .
\]
Take $R=2^k$, $k\ge1$. The balls $B_{2^k\rho_1}(2^ke)$ are pairwise disjoint, since $|2^ke-2^je|\ge2^{\max\{j,k\}-1}$ for $j\ne k$, while the sum of the radii is \changed{at most} $3\cdot2^{\max\{j,k\}-1}\rho_1<2^{\max\{j,k\}-1}$ (as $\rho_1<2r_0<\frac13$). On $B_{2^k\rho_1}(2^ke)$ we have $1+|y|\le1+2^k(1+\rho_1)\le2^{k+2}$. Therefore
\[
\int_{\R^N}\frac{u^{p-1}}{(1+|y|)^{N+sp}}dy\ge\sum_{k\ge1}c_5^{p-1}2^{k\frac{sp(p-1)}{p-1-q}}\,|B_1|\rho_1^N2^{kN}\,2^{-(k+2)(N+sp)}=c\sum_{k\ge1}2^{k\frac{spq}{p-1-q}}=\infty,
\]
because $\frac{sp(p-1)}{p-1-q}-sp=\frac{spq}{p-1-q}>0$. This contradicts $u\in\Ls$. Hence \rev{$u=0$ in $\R^N\setminus\{0\}$} in all cases.
\end{proof}

\begin{proof}[Proof of Theorem~\ref{thm:liou}\textup{(b)}]
Let $q>q^*$ and $\gamma:=\frac{sp}{q-p+1}$. Since $q-p+1>\frac{sp}{\beta^+}>0$, we have $0<\gamma<\beta^+$, and $\gamma q=\gamma(p-1)+\gamma(q-p+1)=\gamma(p-1)+sp$. By definition of $\beta^+$ there are $\hat\beta\in(\gamma,\beta^+]$ and $v\in\mathcal S_{\hat\beta}(\omega)$. Let $a:=\gamma/\hat\beta\in(0,1)$ and $z:=v^a$. Then $z\in H_\gamma$, $z\ge0$, $z>0$ in $\Cw$, and $z\in\Ls$ because $0<\gamma<\frac N{p-1}$.

\emph{Claim: $\pl z\ge c_z|x|^{-\gamma(p-1)-sp}$ in $\Cw$ for some $c_z>0$.} By homogeneity it is enough to check the inequality at points $x\in\Cw\cap\Sp$, for test functions with radius $r<1$ (by Lemma~\ref{lem:grad}, applicable since $\varphi$ touches $z$ from below in $B_r(x)$, decreasing $r$ decreases $\pl\varphi_r(x)$). Let $Q:=\mathcal C_{\omega_0}\cap\{2\le|y|\le4\}$; then $1\le|x-y|\le5$ for $y\in Q$, so $Q\subset\R^N\setminus B_r(x)$. By Lemma~\ref{lem:power} with $D=\Cw$ and $F=0$, and since $\Gamma_a\ge0$,
\[
\pl\varphi_r(x)\ge\int_{\R^N\setminus B_r(x)}\frac{\Gamma_a(x,y)}{|x-y|^{N+sp}}dy\ge g(x):=\int_Q\frac{\Gamma_a(x,y)}{|x-y|^{N+sp}}dy .
\]
We extend $g$ to $\overline{\Cw}\cap\Sp$ by the same formula, with $\Gamma_a(x,y):=+\infty$ when $v(x)=0<v(y)$; this is consistent with the limit of $\Gamma_a(x_k,y)=-\big(v(y)^a-v(x_k)^a\big)^{p-1}+a^{p-1}v(x_k)^{(a-1)(p-1)}(v(y)-v(x_k))^{p-1}$ as $v(x_k)\to0$. Then:
\begin{itemize}
\item $g>0$ on $\overline{\Cw}\cap\Sp$. If $v(x)=0$ this is clear, since $v>0$ on $Q$. If $v(x)>0$, then $\Gamma_a(x,y)>0$ whenever $v(y)\ne v(x)$ (Lemma~\ref{lem:power}), and the set $\{y\in Q:\ v(y)=v(x)\}$ is a null set: in polar coordinates $y=t\theta$, for each $\theta\in\omega_0$ the function $t\mapsto v(t\theta)=t^{-\hat\beta}v(\theta)$ is strictly decreasing, so it takes the value $v(x)$ for at most one $t$.
\item $g$ is lower semicontinuous on $\overline{\Cw}\cap\Sp$: for each $y\in Q$, $x\mapsto\Gamma_a(x,y)|x-y|^{-N-sp}$ is nonnegative and lower semicontinuous (continuous where $v(x)>0$, and tending to $+\infty$ as $v(x)\to0$), so Fatou's lemma applies.
\end{itemize}
A positive lower semicontinuous function on a compact set attains a positive minimum; let $c_z:=\min_{\overline{\Cw}\cap\Sp}g>0$. Then $\pl\varphi_r(x)\ge c_z$ for every admissible $\varphi$ touching $z$ from below at $x\in\Cw\cap\Sp$, which proves the claim.

%\emph{Conclusion.}
 Let $\|z\|:=\max_{\Sp}z>0$ and $w:=\kappa z$ with $\kappa^{q-p+1}:=c_z\|z\|^{-q}$. Then, in the viscosity sense in $\Cw$,
\[
\pl w=\kappa^{p-1}\pl z\ge\kappa^{p-1}c_z|x|^{-\gamma(p-1)-sp}=\kappa^q\|z\|^q|x|^{-\gamma q}\ge\kappa^qz^q=w^q ,
\]
where we used $\gamma(p-1)+sp=\gamma q$ and $z(x)\le\|z\||x|^{-\gamma}$. The function $w$ is nonnegative, continuous in $\R^N\setminus\{0\}$, positive in $\Cw$, $(-\gamma)$-homogeneous and belongs to $\Ls$.
\end{proof}

\begin{remark}\label{rem:crit}
{At the critical exponent} $q=q^*$ the argument above gives no contradiction; \rev{the critical exponent is included in \cite{AS}, and in \cite{gab} for nonlocal Isaacs operators. Their arguments use a quantitative strong maximum principle for differences of subsolutions and supersolutions.}
\end{remark}

\section{The  exponent \texorpdfstring{$\beta^-$}{beta-}, a Hopf\changed{-}type result and \changed{global} Hölder regularity}\label{sec:neg}

In this section we consider homogeneous functions of {positive} degree, that is, $u\in H_\beta$ with $\beta<0$; they vanish at the vertex and grow at infinity. By the integrability criterion of Section~\ref{sec:prel}, a nonnegative $u\in H_\beta$, positive on an open cone, belongs to $\Ls$ if and only if $\beta>-\frac{sp}{p-1}$. Recall that $\beta^-$ and $\alpha^-=-\beta^-$ are defined in \eqref{def:beta}.
Since $\mathcal S_\beta(\omega)\subset\mathcal S_\beta(\omega_0)$ for $\omega_0\subset\omega$, we have $\beta^-(\omega_0)\le\beta^-(\omega)$, i.e.\ $\alpha^-(\omega_0)\ge\alpha^-(\omega)$, if $\omega_0\subset\omega$.

\subsection{Eigenpairs for nonpositive exponents}

\begin{lemma}\label{lem:ext}
The \changed{construction of homogeneous eigenpairs in} Section~\ref{sec:eig} \changed{extends to} exponents in $(-\frac{sp}{p-1},\frac N{p-1})$, with the following modifications:
\begin{itemize}
\item[(i)] Lemma~\ref{lem:aux} holds for every $\beta\in(-\frac{sp}{p-1},\frac N{p-1})$;
\item[(ii)] in Lemma~\ref{lem:Lambda}, if $-\frac{sp}{p-1}<b_1\le b_2<\frac N{p-1}$, there is $\Lambda_0=\Lambda_0(N,s,p,\omega,b_1,b_2)$ such that the conclusion holds for all $\beta\in[b_1,b_2]$;
\item[(iii)] Lemma~\ref{lem:compact} holds for $-\frac{sp}{p-1}<b_1\le b_2<\frac N{p-1}$;
\item[(iv)] for every $\beta\in(-\frac{sp}{p-1},\frac N{p-1})$ there is a homogeneous eigenpair $(\mu_\beta,U_\beta)$ with exponent $\beta$, $\|U_\beta\|=1$ and $-M_\beta\le\mu_\beta\le\Lambda_0(N,s,p,\omega,\beta,\beta)$.
\end{itemize}
\end{lemma}

\begin{proof}
(i) Theorem~\ref{thm:qp} holds for $\beta\in(-\frac{sp}{p-1},\frac N{p-1})$, so $w_K=K|x|^{-\beta}$ is a supersolution as before {(for $\beta=0$, $w_K\equiv K$ is constant, and any patched test function touching it from below satisfies $\varphi_r\le K=\varphi_r(x)$, so the inequality is immediate)}, and $w_K\in\Ls$; the rest of the proof does not use the sign of $\beta$ (the continuity of $V$ at $\partial\Cw\setminus\{0\}$ follows from $V(y)\le C|y|^{-\beta-s}d(y)^s$ as before).
(ii) In the proof of Lemma~\ref{lem:Lambda} the only place where the sign of $\beta$ enters is the lower bound of $\psi$ on $\mathcal C_{\omega_0}\cap B_{2r_0}(\hat x)$, where $\frac34\le|y|\le\frac32$; there $|y|^{-\beta}\ge(3/2)^{-b_2}$ if $\beta\in(0,b_2]$ and $|y|^{-\beta}=|y|^{|\beta|}\ge(3/4)^{|b_1|}$ if $\beta\in[b_1,0]$. So $(3/2)^{-\bar\beta}$ is replaced by $\min\{(3/2)^{-\max\{b_2,0\}},(3/4)^{\max\{-b_1,0\}}\}$.
(iii) The bound $U_j\le|x|^{-b_1}+|x|^{-b_2}$ still holds, and the right-hand side belongs to $\Ls$ because $b_1(p-1)>-sp$ and $b_2(p-1)<N$; the right-hand sides of the equations \mathred{and the functions $U_j$} are bounded in $\{\frac12<|x|<2\}$ uniformly in $j$. {The remainder of the proof is unchanged.}
(iv) The proof of Proposition~\ref{prop:eig} uses only (i)--(iii).
\end{proof}

\begin{lemma}\label{lem:zero}
Every homogeneous eigenpair $(\mu,U)$ with exponent $\beta=0$ has $\mu>0$. In particular $U\in\mathcal S_0(\omega)$, the set in \eqref{def:beta} is nonempty, and \eqref{mp} has no nonnegative nontrivial $0$-homogeneous solution.
\end{lemma}

\begin{proof}
Normalize $\|U\|=1$ and let $x_0\in\Sp$ with $U(x_0)=1$; then $x_0\in\Cw$ and, $U$ being $0$-homogeneous, $U\le1$ in $\R^N$. The function $\varphi(y):=1+|y-x_0|^{q_0}$ touches $U$ from above at $x_0$, has an isolated critical point at $x_0$ and belongs to $C^2_{q_0}$ (as in the proof of Lemma~\ref{lem:smp}), so it is admissible. For $B_r(x_0)\subset\Cw$ the subsolution property gives
\[
\mu=\mu|x_0|^{-sp}U(x_0)^{p-1}\ge\pl\varphi_r(x_0)=-\int_{B_r(x_0)}\frac{|y-x_0|^{q_0(p-1)}}{|x_0-y|^{N+sp}}dy+\int_{\R^N\setminus B_r(x_0)}\frac{(1-U(y))^{p-1}}{|x_0-y|^{N+sp}}dy .
\]
The first integral is $O(r^{q_0(p-1)-sp})\to0$. The second one is at least $\int_{B_3\setminus\Cw}|x_0-y|^{-N-sp}dy\ge4^{-N-sp}|B_3\setminus\Cw|>0$ (since $U=0$ outside $\Cw$ and $\R^N\setminus\Cw$ has nonempty interior). Letting $r\to0$, $\mu\ge4^{-N-sp}|B_3\setminus\Cw|>0$. Hence $\pl U\ge0$, and $U\in\mathcal S_0(\omega)$ by Lemma~\ref{lem:smp}.
\end{proof}

\begin{corollary}\label{cor:betaplus}
For every $N\ge2$ there is $\beta>0$ with $\mathcal S_\beta(\omega)\ne\emptyset$. Hence $\beta^+$ is well defined and $0<\beta^+<\frac N{p-1}$.
\end{corollary}

\begin{proof}
Take $\beta_j\downarrow0$ with $\beta_j\le\frac N{2(p-1)}$, and let $(\mu_j,U_j)$ be homogeneous eigenpairs with exponent $\beta_j$ and $\|U_j\|=1$ (Lemma~\ref{lem:ext}(iv)). If $\mu_j>0$ for some $j$, then $\pl U_j\ge0$ in $\Cw$, so $U_j\in\mathcal S_{\beta_j}(\omega)$. Otherwise $-1-\max_{[0,\frac N{2(p-1)}]}|C|\le\mu_j\le0$ for all $j$, and by Lemma~\ref{lem:ext}(iii) (with $b_1=0$, $b_2=\frac N{2(p-1)}$) a subsequence converges to a homogeneous eigenpair with exponent $0$ and eigenvalue $\lim\mu_j\le0$, contradicting Lemma~\ref{lem:zero}. The upper bound is Lemma~\ref{lem:range}(b).
\end{proof}

\subsection{A nonexistence principle at the vertex}

{The following is the analogue} of Proposition~\ref{prop:nonex}, {with the roles of the origin and infinity interchanged}. {The supersolution is required only near the vertex.}

\begin{proposition}\label{prop:nonexminus}
Let $\beta<0$ and let $U\in H_\beta\cap\Ls$, $U\ge0$, $U\not\equiv0$, $U=0$ in $\R^N\setminus\Cw$, be a viscosity subsolution of $\pl U\le0$ in $\Cw$. Let $v\in C(\R^N\setminus\{0\})\cap\Ls$, $v\ge0$ in $\R^N$, $v>0$ in $\Cw\cap B_6$, be a viscosity supersolution of $\pl v\ge0$ in $\Cw\cap B_6$. Then there is $\varepsilon>0$ such that
\[
v\ge\varepsilon\,U\qquad\text{in }\Cw\cap B_{1/2}.
\]
In particular, $\mathcal S_{\hat\beta}(\omega)=\emptyset$ for every $\hat\beta\in(-\frac{sp}{p-1},\beta)$.
\end{proposition}

\begin{proof}
%\emph{Step 0.}
 Let $L:=\max_{\Sp}U>0$ and $|\beta|=-\beta$. Then $U(y)\le L|y|^{|\beta|}$, so $U$ is continuous in $\R^N$ with $U(0)=0$; $|\beta|(p-1)<sp$ because $U\in\Ls$; and, by Lemma~\ref{lem:boundary}(b) and homogeneity, $U(y)\le C|y|^{|\beta|-s}d(y)^s$ in $\Cw$.

%\emph{Step 1: a subsolution in $\Cw\cap B_{1/2}$.}
 Let $A':=\mathcal C_{\omega_0}\cap\{\frac74\le|y|\le\frac94\}$, $\chi\in C_c(\Cw\cap\{\frac32<|y|<\frac52\})$ with $0\le\chi\le1$, $\chi=1$ on $A'$, and $\vartheta\in C(\R)$, $0\le\vartheta\le1$, $\vartheta(t)=1$ for $t\le1$, $\vartheta(t)=0$ for $t\ge\frac43$. For $C_b>0$ let
\[
z(y):=\vartheta(|y|)U(y)+C_b\chi(y).
\]
Then $z$ is continuous, bounded, compactly supported, $z\ge0$, $z=0$ outside $\Cw$, and $z=U$ in $B_1$. By Lemma~\ref{lem:loc} with $D'=\Cw\cap B_1$, $z$ is a viscosity subsolution of $\pl z\le E$ in $\Cw\cap B_1$, where, with $a:=U(x)$,
\[
E(x)=\int_{\Cw\setminus B_1}\frac{\Psi_p(a-z(y))-\Psi_p(a-U(y))}{|x-y|^{N+sp}}dy .
\]
For $|x|\le\frac12$, {we obtain} $0\le a\le L$, and $\frac12|y|\le|x-y|\le\frac32|y|$ for $|y|\ge1$. On $P:=\{y\in\Cw\setminus B_1:\ z(y)<U(y)\}$ the integrand is nonnegative and, as in Proposition~\ref{prop:nonex}, bounded by $2^p(L^{p-1}+L^{p-1}|y|^{|\beta|(p-1)})(|y|/2)^{-N-sp}$; hence the contribution of $P$ is at most
\[
K_1:=2^{p+N+sp}L^{p-1}\int_{|y|\ge1}\big(1+|y|^{|\beta|(p-1)}\big)|y|^{-N-sp}dy<\infty,
\]
since $|\beta|(p-1)<sp$. {Where} $z\ge U$, {the integrand is} $\le0$. If $C_b\ge2L(9/4)^{|\beta|}$, then on $A'$ we have $z=C_b>U$ (as $\vartheta=0$ there), so $A'\cap P=\emptyset$ and the \changed{numerator} is at most
\[
\Psi_p(L-C_b)-\Psi_p\big(-L(9/4)^{|\beta|}\big)=-D(C_b),\qquad D(C_b):=(C_b-L)^{p-1}-\big(L(9/4)^{|\beta|}\big)^{p-1}>0,
\]
while $|x-y|\le\frac32\cdot\frac94=\frac{27}8$ there. Thus $E(x)\le K_1-D(C_b)|A'|(27/8)^{-N-sp}\le0$ if $C_b$ is large, which we fix. So $z$ is a viscosity subsolution of $\pl z\le0$ in $\Cw\cap B_{1/2}$.

%\emph{Step 2: ordering away from the vertex.}
 {The function} $z$ {now vanishes} outside $\Cw$ and in $\{|y|\ge\frac52\}$. In $\Cw\cap\{\frac12\le|y|\le\frac52\}$ we have $v\ge c\,d^s$ by Lemma~\ref{lem:boundary}(a) (with $r_1=\frac12$, $r_2=\frac52$, whose hypotheses only involve $\Cw\cap\{\frac14\le|y|\le5\}\subset\Cw\cap B_6$), $U\le C''d^s$ by \changed{the boundary bound above}, and $v\ge c'>0$ on the compact set $\operatorname{supp}\chi\subset\Cw\cap B_6$. Hence $\varepsilon z\le v$ in $\R^N\setminus B_{1/2}$, with $\varepsilon:=(\frac{C''}c+\frac{C_b}{c'})^{-1}$.

%\emph{Step 3: comparison near the vertex.}
 Next, let $\delta>0$ and $r\in(0,\frac12)$ with $\varepsilon Lr^{|\beta|}\le\delta$, and $D_r:=\Cw\cap(B_{1/2}\setminus\overline{B_r})$. The function $\varepsilon z-\delta$ is a bounded continuous subsolution of $\pl(\varepsilon z-\delta)\le0$ in $D_r$, and $v\ge0$ is a supersolution of $\pl v\ge0$ in $D_r\subset\Cw\cap B_6$, continuous near $\overline{D_r}$. Outside $D_r$ we have $\varepsilon z-\delta\le v$: outside $\Cw$ trivially; in $\R^N\setminus B_{1/2}$ by \changed{the preceding ordering}; in $\Cw\cap\overline{B_r}$ because there $\varepsilon z=\varepsilon U\le\varepsilon Lr^{|\beta|}\le\delta$. By (K2), $\varepsilon U-\delta\le v$ in $D_r$. Letting $r\to0$ and then $\delta\to0$, $v\ge\varepsilon U$ in $\Cw\cap B_{1/2}$.

%\emph{Step 4.}
 If $v\in\mathcal S_{\hat\beta}(\omega)$ with $-\frac{sp}{p-1}<\hat\beta<\beta$, then it satisfies the hypotheses. Take $y_0\in\Cw\cap\Sp$ with $U(y_0)>0$. Then for $t\in(0,\frac12)$, \changed{the preceding comparison} gives $\varepsilon t^{|\beta|}U(y_0)\le t^{|\hat\beta|}v(y_0)$, i.e.\ $\varepsilon U(y_0)\le t^{|\hat\beta|-|\beta|}v(y_0)\to0$ as $t\to0$, since $|\hat\beta|>|\beta|$. This is a contradiction.
\end{proof}

\subsection{The range of \texorpdfstring{$\beta^-$}{beta-}, signs and uniqueness}

\begin{lemma}\label{lem:rangeminus}
\begin{itemize}
\item[(a)] There is $\delta_0^->0$, depending only on $N,s,p,\omega$, such that if $\beta\in(-\frac{sp}{p-1},0)$ and $\mathcal S_\beta(\omega)\ne\emptyset$, then $\beta\ge-\frac{sp}{p-1}+\delta_0^-$.
\item[(b)] We have $\beta^-<0$.
\end{itemize}
Consequently $-\frac{sp}{p-1}<\beta^-<0$, i.e.\ $0<\alpha^-<\frac{sp}{p-1}$.
\end{lemma}

\begin{proof}
(a) This is the counterpart of Lemma~\ref{lem:range}(b), {with the contribution now coming from infinity}. Let $v\in\mathcal S_\beta(\omega)$, $\alpha:=-\beta\in(0,\frac{sp}{p-1})$, and assume $\alpha\ge\frac{sp}{2(p-1)}$. With $m_0,\hat x,c^*\le m_0,x_1\in B_{r_0}(\hat x)$ and $\phi$ as in the proof of Lemma~\ref{lem:range}(b), \eqref{eq:range1} holds. Let $R_1:=\max\{2^{1/\alpha},2\}$ and $V:=\mathcal C_{\omega_0}\setminus B_{R_1}$. For $y\in V$: $\phi(y)=0$ (as $B_{r_0}(\hat x)\subset B_{9/8}$), $v(y)=|y|^{\alpha}v(y/|y|)\ge m_0|y|^\alpha\ge2m_0$, so $\Psi_p(\phi(x_1)-v(y))-\Psi_p(\phi(x_1))\le-2^{1-p}m_0^{p-1}|y|^{\alpha(p-1)}$, and $|x_1-y|\le|y|+\frac98\le2|y|$. Hence
\[
M_\eta\ge2^{1-p-N-sp}\int_V|y|^{\alpha(p-1)-N-sp}dy=2^{1-p-N-sp}\mathcal H^{N-1}(\omega_0)\frac{R_1^{\alpha(p-1)-sp}}{sp-\alpha(p-1)} .
\]
Since $\alpha\ge\frac{sp}{2(p-1)}$, $R_1\le R_*:=\max\{2^{2(p-1)/(sp)},2\}$, and $R_1^{\alpha(p-1)-sp}\ge R_*^{-sp}$. Thus $sp-\alpha(p-1)\ge c^-:=2^{1-p-N-sp}\mathcal H^{N-1}(\omega_0)R_*^{-sp}/M_\eta$, and (a) holds with $\delta_0^-:=\min\{\frac{c^-}{p-1},\frac{sp}{2(p-1)}\}$.

(b) Assume by contradiction that $\mathcal S_\beta(\omega)=\emptyset$ for all $\beta\in(-\frac{sp}{p-1},0)$. Let $\beta_j\uparrow0$ with $\beta_j\in[-\frac{sp}{2(p-1)},0)$ and let $(\mu_j,U_j)$ be the eigenpairs of Lemma~\ref{lem:ext}(iv). If $\mu_j\ge0$, then $U_j\in\mathcal S_{\beta_j}(\omega)$ (by Lemma~\ref{lem:smp}), which is excluded; so $\mu_j<0$, and $\mu_j\ge-M_{\beta_j}\ge-1-\max_{[-\frac{sp}{2(p-1)},0]}|C|$. By Lemma~\ref{lem:ext}(iii), a subsequence converges to an eigenpair with exponent $0$ and eigenvalue $\le0$, contradicting Lemma~\ref{lem:zero}.
\end{proof}

\begin{corollary}\label{cor:signminus}
\begin{itemize}
\item[(i)] Let $(\mu,U)$ be a homogeneous eigenpair with exponent $\beta\in(-\frac{sp}{p-1},0]$. If $\beta>\beta^-$ then $\mu>0$; if $\beta<\beta^-$ then $\mu<0$.
\item[(ii)] There is a homogeneous eigenpair $(\underline\mu^-,\underline U^-)$ with exponent $\beta^-$, $\|\underline U^-\|=1$ and $\underline\mu^-\ge0$. In particular $\mathcal S_{\beta^-}(\omega)\neq\emptyset$, while $\mathcal S_\beta(\omega)=\emptyset$ for $\beta\in(-\frac{sp}{p-1},\beta^-)$.
\item[(iii)] If $\Phi\in H_\beta\cap\Ls$ with $\beta\le0$ is a nonnegative nontrivial viscosity solution of \eqref{mp}, then $\beta=\beta^-$.
\end{itemize}
\end{corollary}

\begin{proof}
(i) Let $\beta>\beta^-$ and $\mu\le0$. If $\beta=0$ this contradicts Lemma~\ref{lem:zero}. If $\beta<0$, $U$ is a subsolution of $\pl U\le0$; by definition of the infimum there is $\hat\beta\in[\beta^-,\beta)$ with $\mathcal S_{\hat\beta}(\omega)\ne\emptyset$, contradicting Proposition~\ref{prop:nonexminus}. If $\beta<\beta^-$ and $\mu\ge0$, then $U\in\mathcal S_\beta(\omega)$, contradicting the definition of $\beta^-$.

(ii) Take $\beta_j\downarrow\beta^-$ with $\beta_j\in(\beta^-,\frac{\beta^-}2]$ and eigenpairs $(\mu_j,U_j)$ from Lemma~\ref{lem:ext}(iv). By (i), $\mu_j>0$, and by Lemma~\ref{lem:ext}(ii) (with $b_1=\beta^-$, $b_2=\frac{\beta^-}2$), $\mu_j\le\Lambda_0$. By Lemma~\ref{lem:ext}(iii) (note $\beta^->-\frac{sp}{p-1}$ by Lemma~\ref{lem:rangeminus}) a subsequence converges to an eigenpair with exponent $\beta^-$ and $\underline\mu^-=\lim\mu_j\ge0$. Then $\underline U^-\in\mathcal S_{\beta^-}(\omega)$ by Lemma~\ref{lem:smp}. {The final assertion follows from the definition} of $\beta^-$.

(iii) $\beta=0$ is excluded by Lemma~\ref{lem:zero}. Let $\beta<0$. By Lemma~\ref{lem:smp}, $\Phi\in\mathcal S_\beta(\omega)$, so $\beta\ge\beta^-$. By Proposition~\ref{prop:nonexminus} (with $U=\Phi$), $\mathcal S_{\hat\beta}(\omega)=\emptyset$ for $\hat\beta<\beta$, hence $\beta^-\ge\beta$.
\end{proof}

\begin{remark}\label{rem:halfspace}
If $\omega$ is a hemisphere, $\Cw=\{x\cdot e>0\}$, the function $(x\cdot e)_+^s$ is a nonnegative solution of \eqref{mp}: it is $(s,p)$-harmonic in the half-space (a one-dimensional computation, see \cite{IMS20} for $p\ge2$ and \cite{IM} for $1<p<2$), and, being smooth with nonvanishing gradient in the half-space, it is a pointwise and hence a viscosity solution there, as in the proof of Lemma~\ref{lem:range}(a); so by Corollary~\ref{cor:signminus}(iii) $\beta^-=-s$, i.e.\ $\alpha^-=s$, for every $p>1$. By monotonicity, $\alpha^-(\omega)\ge s$ if $\omega$ is contained in a hemisphere, and $\alpha^-(\omega)\le s$ if $\omega$ contains one.
\end{remark}

\subsection{Existence of fundamental solution  \texorpdfstring{$\Phi^-$}{Phi-}}\label{7.4}

\begin{theorem}\label{thm:phiminus}
If $p-1>sp$, then $\underline\mu^-=0$. Hence $\Phi^-:=\underline U^-$ is a fundamental solution of \eqref{mp} with exponent $\beta^-$: $\Phi^-(tx)=t^{\alpha^-}\Phi^-(x)$, $\Phi^->0$ in $\Cw$, $\Phi^-=0$ outside $\Cw$.
\end{theorem}

\begin{proof}
{We follow the proof in Section~\ref{sec:exist}, indicating only the necessary modifications.} Assume $\mu:=\underline\mu^->0$, let $U:=\underline U^-$, so $U\le|y|^{\alpha^-}$ and, by Lemma~\ref{lem:rangeminus}(a), $\alpha^-(p-1)\le sp-\delta_0^-(p-1)$. {The two-sided bound} \eqref{eq:twosided}, {the choice of} $\alpha$ {in} \eqref{eq:alpha} {and} Lemma~\ref{lem:gradU} {apply without modification} (in the case $r\ge\frac13$ of Lemma~\ref{lem:gradU} we use $U\le(7/4)^{\alpha^-}$ on $B_{3r/2}(x)$) {(these arguments use only the equation on} $\{\frac12\le|y|\le2\}$ and $U\in\Ls$).

 For $a\in(\frac12,1)$ let
{\color{\revisioncolor}
\[
\begin{aligned}
\varepsilon&:=2\alpha^-(1-a),
&Z&:=U^a\in H_{a\beta^-},\\
h(y)&:=|y|^{+\varepsilon},
&W&:=hZ\in H_\gamma,\\
\gamma&:=a\beta^--\varepsilon
=\beta^--(1-a)\alpha^-<\beta^-.
\end{aligned}
\]
}
Then $W\in\Ls$ if $|\gamma|(p-1)<sp$, and $W\in\mathcal S_\gamma(\omega)$ would contradict Corollary~\ref{cor:signminus}(ii). Moreover, condition \eqref{eq:a1} is replaced by
\[
1+as(p-1)-sp\ge\tfrac{1-s}2,\qquad|\gamma|(p-1)\le sp-\tfrac{\delta_0^-(p-1)}2,\qquad\varepsilon\le\min\big\{1,\varkappa\big\},\quad \varkappa:=\tfrac{\delta_0^-(p-1)}4 .
\]
Since $|1-\tau^{\varepsilon}|\le\varepsilon|\log\tau|\max\{1,\tau^\varepsilon\}$, \eqref{eq:hbounds} holds with $\tau^{\varepsilon}$ in place of $\tau^{-\varepsilon}$, and $|1-h(y)|\le4\varepsilon|y-x|$ for $|y-x|\le\frac12$.

{We next establish the error estimate.} Lemma~\ref{lem:error} holds with the same statement. For $p\ge2$, Regions II and III are unchanged (they are local), and in Region I we now have $|B|\le2\varepsilon|\log|y||(1+|y|^{|\gamma|})$ and $|A|+|B|\le4(1+|y|^{|\gamma|})$, so the integral is bounded by $C\varepsilon\int|\log|y||(1+|y|^{|\gamma|(p-1)})(1+|y|)^{-N-sp}dy\le C\varepsilon$, because $|\gamma|(p-1)<sp$. For $1<p<2$, now $A(t\theta)=Z(x)-t^{a\alpha^-}Z(\theta)$ is {decreasing} in $t$, $|\partial_tA|=a\alpha^-t^{a\alpha^--1}Z(\theta)$, $|B|=\frac{|1-t^\varepsilon|\,t}{a\alpha^-}|\partial_tA|$, and \eqref{eq:coarea} holds with $|\partial_tA|$. The dyadic shells are treated exactly as before (here $p-1>sp$ is used). The region $t<\frac12$ is now harmless: $|1-t^\varepsilon|\le\varepsilon\log\frac1t$, so $g\le C\varepsilon$, and $A$ maps $(0,\frac12)$ into an interval of length $\le1$. In the region $|z|\ge\frac12$, $t\ge\frac12$, the bound $g\le C\varepsilon$ holds for $t\le2$ and $A([\frac12,2])$ has length $\le2^{sp/(p-1)}$; for $t\ge2$,
\[
g\le C\varepsilon\,\log t\;t^{\varepsilon}\,t^{N}(1+t)^{-N-sp}\le C\varepsilon\,t^{\varepsilon+\varkappa-sp},
\]
and if $\tau=A(t\theta)\le-2$ then $t^{a\alpha^-}\ge t^{a\alpha^-}Z(\theta)=Z(x)-\tau\ge|\tau|$, i.e.\ $t\ge|\tau|^{1/(a\alpha^-)}$, so $g\le C\varepsilon|\tau|^{-(sp-\varepsilon-\varkappa)/(a\alpha^-)}$. Note that $\varepsilon+\varkappa\le2\varkappa\le\frac{sp}2<sp$ (as $\delta_0^-\le\frac{sp}{2(p-1)}$), so $t^{\varepsilon+\varkappa-sp}\le1$ for $t\ge2$ and $g\le C\varepsilon$ also when $\tau\in[-2,1]$. Since $sp-\varepsilon-\varkappa-a\alpha^-(p-1)\ge sp-2\varkappa-(sp-4\varkappa)=2\varkappa$, the change of variables $\tau=A(t\theta)$ gives
\[
\int_2^\infty g\,|A|^{p-2}|\partial_tA|\,dt\le C\varepsilon\Big(\int_{-2}^1|\tau|^{p-2}d\tau+\int_2^\infty\tau^{p-2-\frac{sp-\varepsilon-\varkappa}{a\alpha^-}}d\tau\Big)\le C\varepsilon .
\]

%\emph{Gain and conclusion.}
 Lemma~\ref{lem:gain} {applies without modification} (it only involves $U$ on $Q$ and near $\partial\Cw\cap\Sp$). In place of the choice in the proof of Theorem~\ref{mt}(d) we take $\delta_*:=\delta_1\exp(-4C_E\alpha^-/c_G)$, so that $c_G(1-a)\log\frac{\delta_1}{U(x)}\ge4C_E\alpha^-(1-a)=2C_E\varepsilon$ when $U(x)\le\delta_*$; Lemma~\ref{lem:grad} applies since $\gamma\ne0$ and $a\beta^-\ne0$. The rest of the argument is identical and shows that $W$ is a viscosity supersolution of $\pl W\ge0$ in $\Cw$. Thus $W\in\mathcal S_\gamma(\omega)$ with $-\frac{sp}{p-1}<\gamma<\beta^-$, a contradiction.
\end{proof}

\subsection{A Hopf\changed{-}type lemma at conical points}

{For} $\pl$, {the following result is an analogue of the} {generalised} {Hopf-type lemma at corners} of \cite[Theorem~1.4]{gab}, proved there for nonlocal Isaacs operators, bounded solutions and cones contained in a hemisphere; here {no such restriction on} $\omega$ is needed.

\begin{corollary}\label{cor:hopfcone}
Let $\Omega\subset\R^N$ be open, $x_0\in\partial\Omega$, and assume that $x_0+\mathcal R(\Cw\cap B_\rho)\subset\Omega$ for some rotation $\mathcal R$ and some $\rho>0$. Let $v\in C(\R^N)\cap\Ls$, $v\ge0$ in $\R^N$, $v>0$ in $\Omega$, be a viscosity supersolution of $\pl v\ge0$ in $\Omega\cap B_\rho(x_0)$. Then for every compact set $K\subset\omega$ there is $c>0$ such that, with $\alpha:=\alpha^-(\omega)$,
\begin{equation}\label{eq:hopfcone}
v\big(x_0+\mathcal R x\big)\ge c\,|x|^{\alpha}\qquad\text{for all }x\in\Cw\cap B_{\rho/12}\text{ with }x/|x|\in K .
\end{equation}
This exponent is optimal: if $\Omega=x_0+\mathcal R\Cw$, the function $v(x_0+\mathcal Rx):=\underline U^-(x)$ of Corollary~\ref{cor:signminus}(ii), extended by $0$ at the vertex, satisfies the assumptions and $v(x_0+\mathcal Rx)=|x|^{\alpha^-}\underline U^-(x/|x|)$.
\end{corollary}

\begin{proof}
Let $\rho':=\min\{\rho,1\}$; we first prove \eqref{eq:hopfcone} on $\Cw\cap B_{\rho'/12}$. On the remaining set $\{x\in\Cw:\ \rho'/12\le|x|\le\rho/12,\ x/|x|\in K\}$, whose image under $x\mapsto x_0+\mathcal Rx$ is a compact subset of $\Omega$, $v$ is bounded below by a positive constant and $|x|^\alpha$ is bounded, so \eqref{eq:hopfcone} holds there after decreasing $c$. Thus we may assume $\rho\le1$. Since $\pl$ is invariant under rigid motions and scaling, $\tilde v(x):=v(x_0+\mathcal R(\rho x/6))$ is a nonnegative supersolution of $\pl\tilde v\ge0$ in $\Cw\cap B_6$, positive there, $\tilde v\in C(\R^N\setminus\{0\})\cap\Ls$. By Lemma~\ref{lem:bothsides} there is a homogeneous eigenpair $(\mu_2,U_2)$ with exponent $\beta^-$ and $\mu_2\le0$, so $U_2$ is a subsolution of $\pl U_2\le0$, positive in $\Cw$. By Proposition~\ref{prop:nonexminus}, $\tilde v\ge\varepsilon U_2$ in $\Cw\cap B_{1/2}$, and $U_2(x)=|x|^{\alpha^-}U_2(x/|x|)\ge(\min_KU_2)|x|^{\alpha^-}$ when $x/|x|\in K$. Scaling back gives \eqref{eq:hopfcone}. Lemma~\ref{lem:bothsides} only uses results proved before Corollary~\ref{cor:hopfcone}, so there is no circularity.
\end{proof}

\begin{proof}[Proof of Theorem~\ref{thm:neg}] {Part~(a) follows from} Lemma~\ref{lem:rangeminus} and Corollary~\ref{cor:signminus}(ii). 
{Part~(b) follows from} Lemma~\ref{lem:ext}(iv) and Corollary~\ref{cor:signminus}(i). {Part~(c) follows from} Corollary~\ref{cor:signminus}(iii) and Lemma~\ref{lem:zero}.   {Part~(d) follows from} Theorem~\ref{thm:phiminus} {when} $p-1>sp$ {and from} Corollary~\ref{cor:p2} {when} $p=2$.  {Finally, part~(e) is established by} Corollary~\ref{cor:hopfcone}.
\end{proof}

\begin{remark}\label{rem:lipschitz}
If $\Omega$ is a \changed{bounded} Lipschitz domain, every $x_0\in\partial\Omega$ admits an interior circular cone $x_0+\mathcal R(\mathcal C_{\omega_\theta}\cap B_\rho)\subset\Omega$, where $\omega_\theta$ is a spherical cap of fixed aperture $\theta\in(0,\pi/2)$, depending on the Lipschitz character of $\Omega$, and $\rho$ is independent of $x_0$. Corollary~\ref{cor:hopfcone} then gives the Hopf\changed{-}type lower bound \eqref{eq:hopfcone} along the axis of the cone, with the exponent $\alpha^-(\omega_\theta)$, for every $p>1$; here $\theta$ is an aperture admissible at all boundary points, fixed by the Lipschitz character. Since $\omega_\theta$ is contained in a hemisphere, $\alpha^-(\omega_\theta)\ge s$, and $\alpha^-(\omega_\theta)$ is nonincreasing in $\theta$ (Remark~\ref{rem:halfspace}); for the hemisphere itself $\alpha^-=s$, the rate of the classical Hopf lemma (K6) in $C^{1,1}$ domains.
\end{remark}
\subsection{Global Hölder regularity}
{\begin{theorem}[Boundary decay and global H\"older regularity in Lipschitz domains]\label{cor:lipreg}
		Let $\Omega\subset\R^N$ be a bounded Lipschitz domain. For $\theta\in(0,\frac\pi2)$ let $\Sigma_\theta\subset\Sp$ be the open spherical cap of aperture $\theta$ centred at $-e_N$, and $\omega^c_\theta:=\Sp\setminus\overline{\Sigma_\theta}$. There are $\theta_0\in(0,\frac\pi2)$ and $\rho>0$, depending on the Lipschitz character of $\Omega$, such that every $x_0\in\partial\Omega$ admits a rotation $\mathcal R$ with $x_0+\mathcal R(\mathcal C_{\Sigma_{\theta_0}}\cap B_{2\rho})\subset\R^N\setminus\Omega$. Let $f\in C(\overline\Omega)$ and let $u\in C(\R^N)$, $u=0$ in $\R^N\setminus\Omega$, be a viscosity solution of $\pl u=f$ in $\Omega$. Then, for every $\theta\in(0,\theta_0)$ and every $\alpha\in(0,\alpha^-(\omega^c_\theta))$,
		\begin{equation}\label{eq:lipdecay}
			|u(x)|\le C\,\dist(x,\partial\Omega)^{\alpha}\qquad\text{for }x\in\Omega,
		\end{equation}
		and $u\in C^{\changed{\alpha}}(\R^N)$. The constants depend only on $N,s,p,\Omega,\theta,\alpha,\|f\|_\infty$ and $\|u\|_\infty$. 
\end{theorem}}

\begin{proof} Since $\omega^c_\theta$ contains a hemisphere, we have $0<\alpha^-(\omega^c_\theta)\le s$.\\
	 {Fix} $\theta<\theta_0$, $\alpha<\alpha^-(\omega^c_\theta)$, and let $(\mu,U)$ be a homogeneous eigenpair for $\omega^c_\theta$ with exponent $-\alpha\in(\beta^-(\omega^c_\theta),0)$. {By Theorem~}\ref{thm:neg}(b), {we obtain}  $\mu>0$. The set $K:=\Sp\setminus\Sigma_{\theta_0}$ is a compact subset of $\omega^c_\theta$, so $c_1|x|^\alpha\le U(x)\le C_U|x|^\alpha$ for $x/|x|\in K$, and
	\[
	\pl U=\mu|x|^{-sp}U^{p-1}\ge\mu c_1^{p-1}|x|^{\alpha(p-1)-sp}\ge\mu c_1^{p-1}\rho^{\alpha(p-1)-sp}=:c_2\qquad\text{in }\mathcal C_K\cap B_\rho,
	\]
	because $\alpha(p-1)<sp$. Let $x_0\in\partial\Omega$, $\mathcal R$ as in the statement, so that $\Omega\cap B_{2\rho}(x_0)\subset x_0+\mathcal R\,\mathcal C_K$, and let $D:=\Omega\cap B_\rho(x_0)$, $M:=\|u\|_\infty$. Let $\chi\in C_c(B_{2\rho}(x_0))$, $0\le\chi\le1$, $\chi=1$ in $B_{3\rho/2}(x_0)$, and $\tilde u:=\chi u$. By Lemma~\ref{lem:loc}, $\tilde u$ is a viscosity solution of $\pl\tilde u=f+E$ in $\Omega\cap B_{3\rho/2}(x_0)$, and, by Lemma~\ref{lem:loc}(ii) with $D''=\overline{B_{5\rho/4}(x_0)}$ and $\delta=\rho/4$, $E$ is continuous near $\overline D$ with $|E|\le C_0M^{p-1}$, $C_0=C_0(N,s,p,\rho,\Omega)$. Let $W(x):=A\,U(\mathcal R^{-1}(x-x_0))$, with $A$ {chosen sufficiently large}, {so that} $A^{p-1}c_2\ge\|f\|_\infty+C_0M^{p-1}$ and $Ac_1\rho^\alpha\ge M$. Then $W$ is a viscosity supersolution of $\pl W\ge\|f\|_\infty+C_0M^{p-1}$ in $D$, and $\tilde u\le W$ in $\R^N\setminus D$. Moreover, outside $\Omega\cap B_{2\rho}(x_0)$ we have $\tilde u=0\le W$, and in $(\Omega\cap B_{2\rho}(x_0))\setminus B_\rho(x_0)$, $W\ge Ac_1\rho^\alpha\ge M\ge\tilde u$. By (K2), $u=\tilde u\le W\le AC_U|x-x_0|^\alpha$ in $D$; the same argument applied to $-u$ gives $|u(x)|\le AC_U|x-x_0|^\alpha$ in $D$. Since $U$ and $A$ do not depend on $x_0$, choosing $x_0$ with $|x-x_0|=\dist(x,\partial\Omega)$ gives \eqref{eq:lipdecay} when $\dist(x,\partial\Omega)<\rho$, and \eqref{eq:lipdecay} is trivial otherwise.\\
	
Next, set $\Theta:=\min\{1,sp/(p-1)\}$ and $\eta:=(\alpha+\Theta)/2$. Since $\alpha<\alpha^-(\omega_\theta^c)\le s<\Theta$, we have $\alpha<\eta<\Theta$, and (K8) applies with exponent $\eta$ for every $p>1$. By (K1), $u$ is a local weak solution. Write $d(x):=\dist(x,\partial\Omega)$ in $\Omega$, so that $|u|\le Cd^\alpha$ in $\R^N$ with $d:=0$ outside $\Omega$, and let $x,y\in\R^N$ with $d(x)\ge d(y)$, $x\in\Omega$. The case $x,y\notin\Omega$ is trivial, and if $|x-y|\ge1$ then $|u(x)-u(y)|\le2M|x-y|^{\changed{\alpha}}$. If $d(x)/4\le|x-y|<1$, then $|u(x)-u(y)|\le2Cd(x)^\alpha\le C'|x-y|^\alpha$. If $|x-y|<d(x)/4$, let $r:=d(x)/2$. In $B_r(x)$ we have $|u|\le C(3r)^\alpha$, and $|u(z)|\le C(3|z-x|)^\alpha$ for $|z-x|\ge r$, so that, since $\alpha(p-1)<sp$,
	\[
	\Big(r^{sp}\int_{\R^N\setminus B_r(x)}\frac{|u(z)|^{p-1}}{|z-x|^{N+sp}}dz\Big)^{\frac1{p-1}}\le C r^\alpha .
	\]
	\changed{For $\xi\in B_{r/2}(x)$ and $|z-\xi|\ge {r/4}$, the $1$-Lipschitz property of $d$ gives
		\[
		d(z)\le\tfrac52r+|z-\xi|\le{11}|z-\xi|.
		\]
		Thus the same calculation gives $\mathrm{Tail}(u;\xi,{r/4})\le Cr^\alpha$ and
		$|u|{\le} C(3r)^\alpha$
		in ${B_{r/4}(\xi)}\subset B_r(x).$
		Applying (K8) in $B_{r/2}(\xi)$, rescaled to $B_1$  {(so that $B_{1/2}$ and $B_{1/32}$ correspond to $B_{r/4}(\xi)$ and $B_{r/64}(\xi)$)},
	%	(so that $B_{1/4}$ and $B_{1/32}$ correspond to $B_{r/8}(\xi)$ and $B_{r/64}(\xi)$, respectively)
		yields uniform estimates on $B_{r/64}(\xi)$. The local $L^\infty$ bound {then controls} pairs of
		points separated by a distance of at least $r/64$. Consequently,}
	\[
	[u]_{C^{\changed{\eta}}(B_{r/2}(x))}
	\le Cr^{-\changed{\eta}}\big(r^\alpha+r^{\frac{sp}{p-1}}\|f\|_\infty^{\frac1{p-1}}\big)
	\le Cr^{\alpha-\changed{\eta}},
	\]
	\begin{figure}[htbp]
    \centering
    \includegraphics[width=0.8\textwidth]{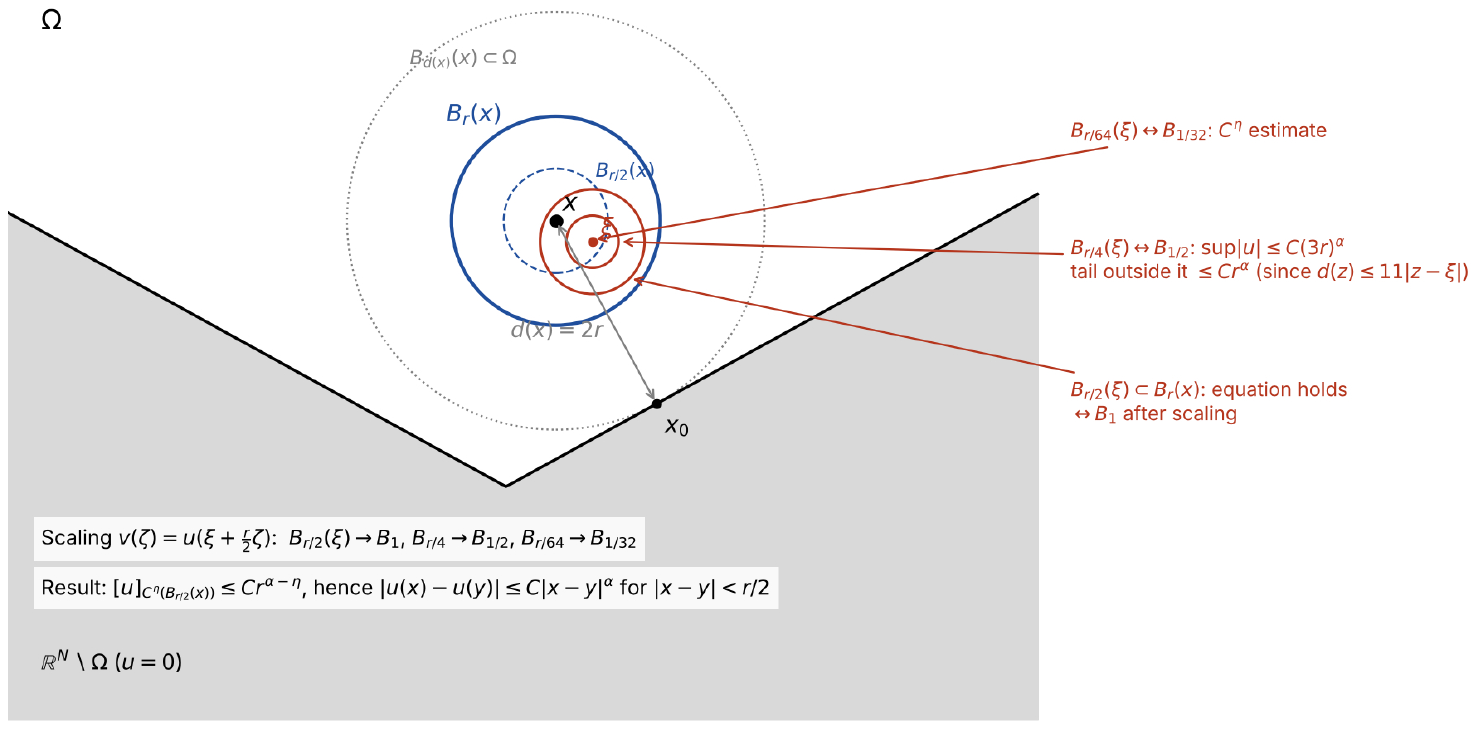}
    \caption{\footnotesize Geometry of the interior rescaling used in the proof of global H\"older regularity. The ball $B_{r/2}(\xi)$ is mapped to $B_1$, while {$B_{r/4}(\xi)$ and $B_{r/64}(\xi)$ correspond to $B_{1/2}$ and $B_{1/32}$}, respectively.}
    \label{fig:holder-rescaling}
\end{figure}

\changed{where {the last constant is uniform because} $\alpha<sp/(p-1)$ and $\Omega$ is bounded. Hence, for $|x-y|<r/2$,}
	\changed{\[
		|u(x)-u(y)|\le Cr^{\alpha-\eta}|x-y|^\eta
		=C|x-y|^\alpha\left(\frac{|x-y|}{r}\right)^{\eta-\alpha}
		\le C|x-y|^\alpha.
		\]}
	{This completes the proof.}
\end{proof}

\section{Uniqueness of the eigenvalue, separation and the restriction \texorpdfstring{$p-1>sp$}{p-1>sp}}\label{sec:sep}

The restriction $p-1>sp$ in Theorems~\ref{mt}(d) and \ref{thm:neg}(d) is used only to show that the eigenvalue at $\beta^\pm$ vanishes, through the error estimate of Lemma~\ref{lem:error}. In this section we show that existence holds for every $p>1$ {provided one can} {separate} a strict supersolution from a strict subsolution with the same homogeneity, {and establish this separation} when $p=2$ for every $s\in(0,1)$.

\begin{definition}\label{def:sp}
Let $\beta\in(-\frac{sp}{p-1},\frac N{p-1})$. We say that the {separation property} $(\mathrm{SP}_\beta)$ holds if there are no two homogeneous eigenpairs $(\mu_1,U_1)$, $(\mu_2,U_2)$ with exponent $\beta$ such that $\mu_2<0<\mu_1$.
\end{definition}

{The separation property follows from uniqueness} of the eigenvalue at $\beta$, i.e.\ by the property that all homogeneous eigenpairs with exponent $\beta$ have the same eigenvalue ({which is weaker than uniqueness of the eigenfunction itself}); by Lemma~\ref{lem:bothsides}, at $\beta=\changed{\beta^-}$ and at $\beta=\beta^+$, this common eigenvalue is necessarily $0$.

\begin{lemma}\label{lem:bothsides}
There are homogeneous eigenpairs $(\mu_1,U_1)$ and $(\mu_2,U_2)$ with exponent $\beta^+$ such that $\mu_1\ge0\ge\mu_2$. {The same conclusion holds with} $\beta^-$ in place of $\beta^+$.
\end{lemma}

\begin{proof}
The pair $(\mu_1,U_1)$ is given by Proposition~\ref{prop:limit}. {For the second pair $(\mu_2,U_2)$, let} $b_2:=\frac12(\beta^++\frac N{p-1})$, take $\beta_j\downarrow\beta^+$ with $\beta_j\in(\beta^+,b_2]$, and let $(\mu_j,U_j)$ be the eigenpairs of Proposition~\ref{prop:eig}, with $\|U_j\|=1$. By Corollary~\ref{cor:sign}, $\mu_j<0$, and $\mu_j\ge-M_{\beta_j}\ge-1-\max_{[\beta^+,b_2]}|C|$. By Lemma~\ref{lem:compact} (with $b_1=\beta^+$), a subsequence converges to an eigenpair with exponent $\beta^+$ and eigenvalue $\mu_2=\lim\mu_j\le0$. For $\beta^-$ {the argument is analogous}, using Corollary~\ref{cor:signminus}(ii) for $\mu_1\ge0$, and eigenpairs with exponents $\beta_j\uparrow\beta^-$, $\beta_j\in[b_1,\beta^-)$ with $b_1:=\frac12(\beta^--\frac{sp}{p-1})$, (which have $\mu_j<0$ by Corollary~\ref{cor:signminus}(i) and $\mu_j\ge-1-\max_{[b_1,\beta^-]}|C|$) together with Lemma~\ref{lem:ext}(iii),(iv) for $\mu_2\le0$.
\end{proof}

\begin{remark}\label{rem:decay}
With $U=U_2$ (exponent $\beta^+$, $\mu_2\le0$) in Proposition~\ref{prop:nonex}, every positive supersolution $v$ as in that proposition satisfies $v\ge\varepsilon U_2$ in $\Cw\setminus B_2$, hence $v(x)\ge c|x|^{-\beta^+}$ for $|x|\ge2$, $x/|x|\in K\Subset\omega$, for every $p>1$. Since $U_1\in\mathcal S_{\beta^+}(\omega)$, this decay rate is optimal; it is the counterpart at infinity of Corollary~\ref{cor:hopfcone}.
\end{remark}

\begin{theorem}\label{thm:sep}
Let $p>1$. If $(\mathrm{SP}_{\beta^+})$ holds, then \eqref{mp} has a fundamental solution with exponent $\beta^+$. If $(\mathrm{SP}_{\beta^-})$ holds, then \eqref{mp} has a fundamental solution with exponent $\beta^-$.
\end{theorem}

\begin{proof}
Let $(\mu_1,U_1)$, $(\mu_2,U_2)$ be given by Lemma~\ref{lem:bothsides}. If $\mu_1=0$ or $\mu_2=0$, the corresponding eigenfunction is a fundamental solution. Otherwise $\mu_2<0<\mu_1$, which is excluded by $(\mathrm{SP}_{\beta^\pm})$.
\end{proof}

\begin{proposition}\label{prop:sp2}
If $p=2$, then $(\mathrm{SP}_\beta)$ holds for every $\beta\in(-2s,N)$ and every $s\in(0,1)$.
\end{proposition}

\begin{proof}
Let $(\mu_1,U_1)$, $(\mu_2,U_2)$ be eigenpairs with exponent $\beta$ and $\mu_2<0<\mu_1$.

%\emph{Step 1: ordering.}
 Since $\mu_1>0$, $U_1$ is a supersolution of $(-\Delta)^sU_1\ge0$; and $U_2$ is a supersolution of $(-\Delta)^sU_2\ge-|\mu_2|2^{2s}U_2$ in $\Cw\cap\{\frac12<|x|<2\}$. By Lemma~\ref{lem:boundary}(a),(b), there are $0<c\le C$ with $c\,d^s\le U_i\le C\,d^s$ on $\Cw\cap\Sp$, $i=1,2$. Since $U_1,U_2$ have the same homogeneity and vanish outside $\Cw$, the set $\{t\ge0:\ tU_2\le U_1\text{ in }\R^N\}$ is an interval $[0,t^*]$ with $c/C\le t^*\le C/c$.

%\emph{Step 2: the difference.}
 {Set }$w:=U_1-t^*U_2\ge0.$ {Thus} $w\in H_\beta,$ {and it is continuous} in $\R^N\setminus\{0\}$, vanishes outside $\Cw$ and belongs to $L^1_{2s}(\R^N)$. Let $D\Subset\Cw$ be a bounded open set. By Remark~\ref{rem:K1} and (K1), suitable cut-offs $\tilde U_i$ of $U_i$, which coincide with $U_i$ near $\overline D$, are weak solutions in $D$ of $(-\Delta)^s\tilde U_i=\mu_i|x|^{-2s}U_i+E_i$ with $E_i$ continuous. Since $(-\Delta)^s$ is linear, $\tilde w:=\tilde U_1-t^*\tilde U_2$ is a continuous weak solution of $(-\Delta)^s\tilde w=f+E_1-t^*E_2$ in $D$, with
\[
f:=|x|^{-2s}\big(\mu_1U_1-t^*\mu_2U_2\big)>0\quad\text{in }\Cw,
\]
hence a viscosity solution by (K1); undoing the cut-off with Lemma~\ref{lem:loc} (the error terms are again $E_1-t^*E_2$, by linearity), $w$ is a viscosity solution of $(-\Delta)^sw=f$ in $D$. As $D$ is arbitrary, $w$ is a viscosity supersolution of $(-\Delta)^sw\ge0$ in $\Cw$.

%\emph{Step 3: separation.}
 If $w\equiv0$, then $\varphi\equiv0$ touches $w$ from below at every $x\in\Cw$, and since $w$ is a viscosity solution of $(-\Delta)^sw=f$ (\changed{as shown above}), $0=(-\Delta)^s\varphi_r(x)\ge f(x)>0$, a contradiction. Thus $w\not\equiv0$, and by Lemma~\ref{lem:smp} $w>0$ in $\Cw$. By Lemma~\ref{lem:boundary}(a) (with $\Lambda=0$), $w\ge c'd^s\ge\frac{c'}{C}U_2$ on $\Cw\cap\Sp$, hence, by homogeneity, $U_1\ge(t^*+\frac{c'}C)U_2$ in $\R^N$. This contradicts the definition of $t^*$. {This proves the proposition.}
\end{proof}

\begin{corollary}\label{cor:p2}
Let $p=2$, $s\in(0,1)$ and $N\ge2$. Then \eqref{mp} has fundamental solutions $\Phi^+\in H_{\beta^+}$ and $\Phi^-\in H_{\beta^-}$, and
\[
\beta^+=N-2s-\beta^-=N-2s+\alpha^- .
\]
In particular $\beta^+>N-2s=\frac{N-sp}{p-1}$, so the lower bound in Theorem~\ref{mt}(a) is strict, and Theorems~\ref{mt}(d) and \ref{thm:neg}(d) hold without the restriction $s<\frac12$.
\end{corollary}

\begin{proof}
Existence of $\Phi^\pm$ follows from Theorem~\ref{thm:sep} and Proposition~\ref{prop:sp2}. Let $\mathcal K u(x):=|x|^{2s-N}u(x/|x|^2)$ be the Kelvin transform; it satisfies $(-\Delta)^s(\mathcal Ku)(x)=|x|^{-N-2s}\big((-\Delta)^su\big)(x/|x|^2)$ for $u\in L^1_{2s}(\R^N)$ which is $C^2$ near $x/|x|^2$ \cite{BZ}, and $\mathcal K$ maps $L^1_{2s}(\R^N)$ into itself, with $\mathcal K\mathcal K u=u$. The identity also holds in the viscosity sense: if $\varphi$ touches $\mathcal Ku$ at $x$ in $B_r(x)$, $r<|x|$, then $\mathcal K\varphi$ touches $u$ at $x/|x|^2$ in $V:=\{y: y/|y|^2\in B_r(x)\}$, and $\mathcal K(\varphi_r)$ equals $\mathcal K\varphi$ in $V$ and $u$ outside; it lies below (above, for subsolutions) the function patched on a ball $B_\rho(x/|x|^2)\subset V$, so the inequality transfers by the pointwise identity. The inversion maps $\Cw$ onto itself. Hence $\mathcal K\Phi^-$ is nonnegative, positive in $\Cw$, zero outside $\Cw$, continuous in $\R^N\setminus\{0\}$, $s$-harmonic in $\Cw$, and $(-(N-2s+\alpha^-))$-homogeneous, with $0<N-2s+\alpha^-<N$ (Lemma~\ref{lem:rangeminus}), so that it belongs to $L^1_{2s}(\R^N)$. By Corollary~\ref{cor:unique}, $\beta^+=N-2s+\alpha^-$. {This completes the proof.}
\end{proof}
\begin{remark}For $p=2$, {these facts are already known} (see  \cite{BB,Michalik,gab}). The proof above only uses the eigenpairs of Section~\ref{sec:eig} and linearity \changed{in the difference equation in} Proposition~\ref{prop:sp2}, {and serves as  a model for} the separation property when $p\ne2$.
\end{remark}

\section{{Acknowledgements}}
R. Biswas {was partially supported by the  VRIP at} Universidad del Bío-Bío, Chile, {through the Proyecto de Investigación Iniciación, grant} IN2647919.\\
 A. Quaas was partially supported by FONDECYT grant 1231585 and
by grant PID2023-148028NB-I00, funded by MICIU/AEI/10.13039/501100011033 and by
ERDF/EU.


{\small
\begin{thebibliography}{99}
\color{\originalrevisioncolor}
\bibitem{ALV}{\color{\originalrevisioncolor}M. Akman, J. Lewis, A. Vogel, Note on an eigenvalue problem for an ODE originating from a homogeneous $p$-harmonic function, \rev{\emph{St. Petersburg Mathematical Journal} \textbf{31} (2020), no.~2, 241--250. %DOI: \doi{10.1090/spmj/1594}.
}}
\bibitem{Amann1976}{\color{\originalrevisioncolor}H. Amann, Fixed point equations and nonlinear eigenvalue problems in ordered Banach spaces, \rev{\emph{SIAM Review} \textbf{18} (1976), no.~4, 620--709. 
	%	DOI: \doi{10.1137/1018114}.
	}}
\bibitem{AS}{\color{\originalrevisioncolor}S. N. Armstrong, B. Sirakov, Nonexistence of positive supersolutions of elliptic equations via the maximum principle, \rev{\emph{Communications in Partial Differential Equations} \textbf{36} (2011), no.~11, 2011--2047.
		% DOI: \doi{10.1080/03605302.2010.534523}.
	 }}
\bibitem{ASS}{\color{\originalrevisioncolor}S. N. Armstrong, B. Sirakov, C. K. Smart, Fundamental solutions of homogeneous fully nonlinear elliptic equations, \rev{\emph{Communications on Pure and Applied Mathematics} \textbf{64} (2011), no.~6, 737--777. 
	%	DOI: \doi{10.1002/cpa.20360}.
	
}}
\bibitem{ASS2}{\color{\originalrevisioncolor}S. N. Armstrong, B. Sirakov, C. K. Smart, Singular solutions of fully nonlinear elliptic equations and applications, \rev{\emph{Archive for Rational Mechanics and Analysis} \textbf{205} (2012), no.~2, 345--394. 
	%	DOI: \doi{10.1007/s00205-012-0505-8}.
	}}
\bibitem{Ataei}{\color{\originalrevisioncolor}A. Ataei, Boundary behavior of solutions to fractional $p$-Laplacian equation, \rev{\emph{Advances in Calculus of Variations} \textbf{18} (2025), no.~2, 255--273. 
	%	DOI: \doi{10.1515/acv-2023-0048}.
	}}
\bibitem{BB}{\color{\originalrevisioncolor}R. Ba\~nuelos, K. Bogdan, Symmetric stable processes in cones, \rev{\emph{Potential Analysis} \textbf{21} (2004), no.~3, 263--288. %DOI: \doi{10.1023/B:POTA.0000033333.72236.dc}.
}}
\bibitem{bm}{\color{\originalrevisioncolor}B. Barrios, M. Medina, Equivalence of weak and viscosity solutions in fractional non-homogeneous problems, \rev{\emph{Mathematische Annalen} \textbf{381} (2021), no.~3--4, 1979--2012. %DOI: \doi{10.1007/s00208-020-02119-w}.
}}
\bibitem{BRS}{\color{\originalrevisioncolor}A. Biswas, S. Roy, A. Sen, Strong comparison principle and symmetry results for the fractional $p$-Laplacian, \rev{preprint (2026), \href{https://arxiv.org/abs/2606.08559}{arXiv:\nolinkurl{2606.08559}}.}}
\bibitem{BT2025}{\color{\originalrevisioncolor}A. Biswas, E. Topp, Lipschitz regularity of fractional $p$-Laplacian, \rev{\emph{Annals of PDE} \textbf{11} (2025), no.~2, Article 27. 
		%DOI: \doi{10.1007/s40818-025-00220-4}.
	}}
\bibitem{BGcone}{\color{\originalrevisioncolor}K. Bogdan, T. Grzywny, Heat kernel of fractional Laplacian in cones, \rev{\emph{Colloquium Mathematicum} \textbf{118} (2010), no.~2, 365--377. 
	%	DOI: \doi{10.4064/cm118-2-1}.
	}}
\bibitem{BKLP}{\color{\originalrevisioncolor}K. Bogdan, P. Knosalla, \L. Le\.{z}aj, D. Pilarczyk, Self-similar solution for fractional Laplacian in cones, \rev{\emph{Electronic Journal of Probability} \textbf{29} (2024), Article 54, 24 pp. 
	%	DOI: \doi{10.1214/24-EJP1111}.
	}}
\bibitem{BSS}{\color{\originalrevisioncolor}K. Bogdan, B. Siudeja, A. St\'os, Martin kernel for fractional Laplacian in narrow cones, \rev{\emph{Potential Analysis} \textbf{42} (2015), no.~4, 839--859.
		% DOI: \doi{10.1007/s11118-014-9459-9}.
	 }}
\bibitem{BZ}{\color{\originalrevisioncolor}K. Bogdan, T. \.Zak, On Kelvin transformation, \rev{\emph{Journal of Theoretical Probability} \textbf{19} (2006), no.~1, 89--120. DOI: \doi{10.1007/s10959-006-0003-8}.}}
\bibitem{BL}{\color{\originalrevisioncolor}L. Brasco, E. Lindgren, Higher Sobolev regularity for the fractional $p$-Laplace equation in the superquadratic case, \rev{\emph{Advances in Mathematics} \textbf{304} (2017), 300--354. %DOI: \doi{10.1016/j.aim.2016.03.039}.
}}
\bibitem{BLS}{\color{\originalrevisioncolor}L. Brasco, E. Lindgren, A. Schikorra, Higher H\"older regularity for the fractional $p$-Laplacian in the superquadratic case, \rev{\emph{Advances in Mathematics} \textbf{338} (2018), 782--846. 
	%	DOI: \doi{10.1016/j.aim.2018.09.009}.
	}}
\bibitem{BMS}{\color{\originalrevisioncolor}L. Brasco, S. Mosconi, M. Squassina, Optimal decay of extremals for the fractional Sobolev inequality, \rev{\emph{Calculus of Variations and Partial Differential Equations} \textbf{55} (2016), no.~2, Article 23. 
	%	DOI: \doi{10.1007/s00526-016-0958-y}.
	}}
\bibitem{BP}{\color{\originalrevisioncolor}L. Brasco, E. Parini, The second eigenvalue of the fractional $p$-Laplacian, \rev{\emph{Advances in Calculus of Variations} \textbf{9} (2016), no.~4, 323--355. 
	%	DOI: \doi{10.1515/acv-2015-0007}.
	}}
\bibitem{Dancer1973}{\color{\originalrevisioncolor}E. N. Dancer, Global solution branches for positive mappings, \rev{\emph{Archive for Rational Mechanics and Analysis} \textbf{52} (1973), no.~2, 181--192. 
	%	DOI: \doi{10.1007/BF00282326}.
	}}
\bibitem{Grube}F. Grube, Boundary regularity and Hopf lemma for nondegenerate stable operators, preprint (2024), arXiv:2410.00829.
\bibitem{hopf}{\color{\originalrevisioncolor}L. M. Del Pezzo, A. Quaas, A Hopf's lemma and a strong minimum principle for the fractional $p$-Laplacian, \rev{\emph{Journal of Differential Equations} \textbf{263} (2017), no.~1, 765--778. 
	%	DOI: \doi{10.1016/j.jde.2017.02.051}.
	}}
\bibitem{qp}{\color{\originalrevisioncolor}L. M. Del Pezzo, A. Quaas, The fundamental solution of the fractional $p$-Laplacian, \rev{\emph{Nonlinear Differential Equations and Applications NoDEA} \textbf{33} (2026), no.~3, Article 62. 
	%	DOI: \doi{10.1007/s00030-026-01202-6}.
	}}
\bibitem{dkp2}{\color{\originalrevisioncolor}A. Di Castro, T. Kuusi, G. Palatucci, Nonlocal Harnack inequalities, \rev{\emph{Journal of Functional Analysis} \textbf{267} (2014), no.~6, 1807--1836. 
	%	DOI: \doi{10.1016/j.jfa.2014.05.023}.
	}}
\bibitem{Cozzi}M. Cozzi, Regularity results and Harnack inequalities for minimizers and solutions of nonlocal problems: a unified approach via fractional De Giorgi classes, \emph{Journal of Functional Analysis} \textbf{272} (2017), no.~11, 4762--4837. 
%DOI: \doi{10.1016/j.jfa.2017.02.016}.
\bibitem{dkp}{\color{\originalrevisioncolor}A. Di Castro, T. Kuusi, G. Palatucci, Local behavior of fractional $p$-minimizers, \rev{\emph{Annales de l'Institut Henri Poincar\'e C, Analyse non lin\'eaire} \textbf{33} (2016), no.~5, 1279--1299. 
	%	DOI: \doi{10.1016/j.anihpc.2015.04.003}.
	}}
\bibitem{Diening}{\color{\originalrevisioncolor}L. Diening, K. Kim, H.-S. Lee, S. Nowak, Higher differentiability for the fractional $p$-Laplacian, \rev{\emph{Mathematische Annalen} \textbf{391} (2025), no.~4, 5631--5693. 
	%	DOI: \doi{10.1007/s00208-024-03057-7}.
	}}
\bibitem{GL}{\color{\originalrevisioncolor}P. Garain, E. Lindgren, Higher H\"older regularity for the fractional $p$-Laplace equation in the subquadratic case, \rev{\emph{Mathematische Annalen} \textbf{390} (2024), no.~4, 5753--5792.
	%	 DOI: \doi{10.1007/s00208-024-02891-z}.
	 }}
\bibitem{GJ2026}\rev{D. Giovagnoli, D. Jesus, H\"older gradient estimates for $(s,p)$-Poisson equations with measurable and H\"older data, preprint (2026), \href{https://arxiv.org/abs/2610.12012}{arXiv:\nolinkurl{2610.12012}}.}
\bibitem{GJS}{\color{\originalrevisioncolor}D. Giovagnoli, D. Jesus, L. Silvestre, $C^{1+\alpha}$ regularity for fractional $p$-harmonic functions, \rev{preprint (2025), \href{https://arxiv.org/abs/2509.26565}{arXiv:\nolinkurl{2509.26565}}.}}
\bibitem{IM}{\color{\originalrevisioncolor}A. Iannizzotto, S. Mosconi, Fine boundary regularity for the singular fractional $p$-Laplacian, \rev{\emph{Journal of Differential Equations} \textbf{412} (2024), 322--379.
		% DOI: \doi{10.1016/j.jde.2024.08.026}.
	 }}
\bibitem{IMP2023}\rev{A. Iannizzotto, S. Mosconi, N. S. Papageorgiou, On the logistic equation for the fractional $p$-Laplacian, \emph{Mathematische Nachrichten} \textbf{296} (2023), no.~4, 1451--1468. 
%	DOI: \doi{10.1002/mana.202100025
%	}.
}
\bibitem{reg}{\color{\originalrevisioncolor}A. Iannizzotto, S. Mosconi, M. Squassina, Global H\"older regularity for the fractional $p$-Laplacian, \rev{\emph{Revista Matem\'atica Iberoamericana} \textbf{32} (2016), no.~4, 1353--1392. 
		%DOI: \doi{10.4171/RMI/921}.
	}}
\bibitem{IMS20}{\color{\originalrevisioncolor}A. Iannizzotto, S. Mosconi, M. Squassina, Fine boundary regularity for the degenerate fractional $p$-Laplacian, \rev{\emph{Journal of Functional Analysis} \textbf{279} (2020), no.~8, Article 108659. 
	%	DOI: \doi{10.1016/j.jfa.2020.108659}.
	}}
\bibitem{jarohs}{\color{\originalrevisioncolor}S. Jarohs, Strong comparison principle for the fractional $p$-Laplacian and applications to starshaped rings, \rev{\emph{Advanced Nonlinear Studies} \textbf{18} (2018), no.~4, 691--704. 
	%	DOI: \doi{10.1515/ans-2017-6039}.
	}}
\bibitem{lindgren}{\color{\originalrevisioncolor}J. Korvenp\"a\"a, T. Kuusi, E. Lindgren, Equivalence of solutions to fractional $p$-Laplace type equations, \rev{\emph{Journal de Math\'ematiques Pures et Appliqu\'ees} \textbf{132} (2019), 1--26. %DOI: \doi{10.1016/j.matpur.2017.10.004}.
}}
\bibitem{CL}A. Cutr\`i, F. Leoni, On the Liouville property for fully nonlinear equations, \emph{Annales de l'Institut Henri Poincar\'e C, Analyse non lin\'eaire} \textbf{17} (2000), no.~2, 219--245.
\bibitem{FQ}P. Felmer, A. Quaas, Fundamental solutions and Liouville type theorems for nonlinear integral operators, \emph{Advances in Mathematics} \textbf{226} (2011), 2712--2738. 
%DOI: \doi{10.1016/j.aim.2010.09.023}.
\bibitem{KKP}J. Korvenp\"a\"a, T. Kuusi, G. Palatucci, The obstacle problem for nonlinear integro-differential operators, \emph{Calculus of Variations and Partial Differential Equations} \textbf{55} (2016), Article 63. 
%DOI: \doi{10.1007/s00526-016-0999-2}.
\bibitem{KM}{\color{\originalrevisioncolor}I. N. Krol', V. G. Maz'ya, The absence of the continuity and H\"older continuity of the solutions of quasilinear elliptic equations near a nonregular boundary, \rev{\emph{Trudy Moskovskogo Matematicheskogo Obshchestva} \textbf{26} (1972), 75--94. (In Russian.) %\url{https://www.mathnet.ru/eng/mmo258}.
}}
\bibitem{KMSdata}{\color{\originalrevisioncolor}T. Kuusi, G. Mingione, Y. Sire, Nonlocal equations with measure data, \rev{\emph{Communications in Mathematical Physics} \textbf{337} (2015), no.~3, 1317--1368. 
	%	DOI: \doi{10.1007/s00220-015-2356-2}.
	}}
\bibitem{KMSself}{\color{\originalrevisioncolor}T. Kuusi, G. Mingione, Y. Sire, Nonlocal self-improving properties, \rev{\emph{Analysis \& PDE} \textbf{8} (2015), no.~1, 57--114. 
		%DOI: \doi{10.2140/apde.2015.8.57}.
	}}
\bibitem{LindgrenHolder}{\color{\originalrevisioncolor}E. Lindgren, H\"older estimates for viscosity solutions of equations of fractional $p$-Laplace type, \rev{\emph{Nonlinear Differential Equations and Applications NoDEA} \textbf{23} (2016), no.~5, Article 55. 
	%	DOI: \doi{10.1007/s00030-016-0406-x}.
	}}
\bibitem{Michalik}{\color{\originalrevisioncolor}K. Michalik, Sharp estimates of the Green function, the Poisson kernel and the Martin kernel of cones for symmetric stable processes, \rev{\emph{Hiroshima Mathematical Journal} \textbf{36} (2006), no.~1, 1--21. 
		%DOI: \doi{10.32917/hmj/1147883392}.
	}}
\bibitem{gab}{\color{\originalrevisioncolor}G. Nornberg, D. dos Prazeres, A. Quaas, Fundamental solutions and critical Lane--Emden exponents for nonlinear integral operators in cones, \rev{\emph{Journal of Functional Analysis} \textbf{287} (2024), no.~4, Article 110487. 
	%	DOI: \doi{10.1016/j.jfa.2024.110487}.
	}}
\bibitem{Nystrom2026}\rev{K. Nystr\"om, Boundary Harnack inequalities for the fractional $p$-Laplacian in Reifenberg flat domains, preprint (2026), \href{https://arxiv.org/abs/2610.21689}{arXiv:\nolinkurl{2610.21689}}.}
\bibitem{OS}{\color{\originalrevisioncolor}P. Ochoa, A. Salort, Hopf's lemmas and boundary point results for the fractional $p$-Laplacian, \rev{\emph{Discrete and Continuous Dynamical Systems} \textbf{45} (2025), no.~3, 717--734. %DOI: \doi{10.3934/dcds.2024109}.
}}
\bibitem{PV}{\color{\originalrevisioncolor}A. Porretta, L. V\'eron, Separable $p$-harmonic functions in a cone and related quasilinear equations on manifolds, \rev{\emph{Journal of the European Mathematical Society} \textbf{11} (2009), no.~6, 1285--1305. 
	%	DOI: \doi{10.4171/JEMS/182}.
	}}
\bibitem{Rabinowitz1971}{\color{\originalrevisioncolor}P. H. Rabinowitz, Some global results for nonlinear eigenvalue problems, \rev{\emph{Journal of Functional Analysis} \textbf{7} (1971), no.~3, 487--513. 
	%	DOI: \doi{10.1016/0022-1236(71)90030-9}.
	}}
\bibitem{Terracini}{\color{\originalrevisioncolor}S. Terracini, G. Tortone, S. Vita, On $s$-harmonic functions on cones, \rev{\emph{Analysis \& PDE} \textbf{11} (2018), no.~7, 1653--1691. 
	%	DOI: \doi{10.2140/apde.2018.11.1653}.
	}}
\bibitem{Tolksdorf}{\color{\originalrevisioncolor}P. Tolksdorf, On the Dirichlet problem for quasilinear equations in domains with conical boundary points, \rev{\emph{Communications in Partial Differential Equations} \textbf{8} (1983), no.~7, 773--817. 
		%DOI: \doi{10.1080/03605308308820285}.
	}}
\bibitem{YZ}{\color{\originalrevisioncolor}D. Ye, W. Zhang, Comparison principles and symmetry for subquadratic fractional $p$-Laplacian equations, \rev{preprint (2026), \href{https://arxiv.org/abs/2609.02626}{arXiv:\nolinkurl{2609.02626}}.}}
\bibitem{ZhangSing2026}\rev{C. Zhang, Gradient H\"older regularity for singular fractional $p$-Laplace equations, preprint (2026), \href{https://arxiv.org/abs/2608.16243v2}{arXiv:\nolinkurl{2608.16243v2}}.}
\end{thebibliography}
}
\end{document}